\documentclass[a4paper,11pt,leqno,openright]{book}
\usepackage[latin1]{inputenc}
\usepackage[T1]{fontenc}
\usepackage{amsmath,amssymb}  % pour utiliser les caracteres R, Z, C
\usepackage{amsfonts} % {\mathbb R} etc...
\usepackage[english,french]{babel}
\usepackage[all,dvips]{xy}
\usepackage[dvips]{graphicx}
\usepackage{latexsym}
\usepackage{dsfont}\let\mathbb\mathds

\usepackage{maple2e}
\def\emptyline{\vspace{12pt}}
\DefineParaStyle{Maple Output}
\DefineParaStyle{Warning}
\DefineCharStyle{2D Math}
\DefineCharStyle{2D Output}

\newenvironment{Dem}{\textit{{D\'emonstration} :}}{\begin{flushright}$\Box$\end{flushright}}

\newenvironment{Rmq}{\emph{\underline{Remarque} :}}{}
\newenvironment{Rmqs}{\emph{\underline{Remarques} :}}{}

\newtheorem{The}{Th\'eor\`eme}[section]
\newtheorem{Teo}{Th\'eor\`eme}
\newtheorem{Pro}{Proposition}[section]
\newtheorem{Lem}{Lemme}[section]

\newtheorem{Def}{D\'efinition}[section]
\newtheorem{Cor}{Corollaire}[section]

\newcommand{\fs}{\mathcal{O}_C}

\newcommand{\im}{\mathop{\rm Im}\nolimits}
\newcommand{\Nm}{\mathop{\rm Nm}\nolimits}

\newcommand{\He}{\mathcal{H}[3]}
\newcommand{\Nmi}{\mathop{\rm N_{f_i}}\nolimits}
\newcommand{\Fix}{\mathop{\rm Fix}\nolimits}
\newcommand{\F}{\mathop{\rm Fix}\nolimits}
\newcommand{\Hom}{\mathop{\rm Hom}\nolimits}
\newcommand{\Det}{\mathop{\rm det}\nolimits}
\newcommand{\Ker}{\mathop{\rm Ker}\nolimits}
\newcommand{\Sing}{\mathop{\rm Sing}\nolimits}
\newcommand{\Deg}{\mathop{\rm deg}\nolimits}
\newcommand{\Rg}{\mathop{\rm rg}\nolimits}
\newcommand{\su}{\mathcal{M}}
\newcommand{\suu}{\mathcal{S}_{2,\lambda}}
\newcommand{\p}{\mathbb{P}^8}
\newcommand{\D}{\mathcal{D}}

\begin{document}
%\maketitle
\begin{titlepage}
\begin{center}
{\bf
UNIVERSIT\'E DE NICE - SOPHIA ANTIPOLIS

FACULT\'E DES SCIENCES

Laboratoire J. A. Dieudonn\'e

U.M.R. du C.N.R.S. No 6621}

\vspace{2.2cm}
{\Large \bf TH\`ESE}

\vspace{3mm}
pr\'esent\'ee pour obtenir le titre de

Docteur en Sciences

sp\'ecialit\'e : Math\'ematiques

par

\vspace{4mm}
{\Large Angela ORTEGA ORTEGA}

\vspace{2cm}
{\LARGE \linespread{1.2} \bf Sur l'espace des modules des fibr\'es vectoriels de rang 3 
sur une courbe de genre 2 et la cubique de Coble. \par}

\vspace{3.8cm}
Soutenue le 15 septembre 2003 devant le jury compos\'e de :
\end{center}

$$\begin{array}{ll}
\mbox{M. Arnaud BEAUVILLE} & \mbox{Professeur \`a l'Universit\'e de Nice}\\
\mbox{M. Andr\'e HIRSCHOWITZ} & \mbox{Professeur \`a l'Universit\'e de Nice}\\
\mbox{M. Yves LASZLO } & \mbox{Professeur \`a l'Universit\'e de Paris VI}\\
\mbox{M. Mudumbai S. NARASIMHAN} & \mbox{Professeur \`a l'International School for}\\ 
& \mbox{Advanced Studies, Trieste, Italie}\\
\mbox{M. Christian PAULY} & \mbox{Charg\'e de recherche CNRS}\\
\mbox{M. Christoph SORGER} & \mbox{Professeur \`a l'Universit\'e de Nantes}\\
\end{array}$$

%\vspace{1,5mm}
%\begin{center}
%\`a 14 heures, en salle de conf\'erences.
%\end{center}
\end{titlepage}

\frontmatter
\tableofcontents
\chapter{Introduction}
\addcontentsline{chapter}{toc}{Introduction}

L'objectif de cette th\`ese est de r\'esoudre une conjecture, d\^ue \`a I. Dolgachev et \'enonc\'ee par Y. Laszlo
dans~\cite{L2}, concernant l'espace des
modules des fibr\'es vectoriels de rang 3 avec d\'eterminant trivial sur une courbe donn\'ee de genre 2. On va d'abord
\'enoncer un r\'esultat de
g\'eom\'etrie alg\'ebrique classique obtenu par Coble au d\'ebut du XX\`eme si\`ecle, ainsi que quelques g\'en\'eralit\'es
et r\'esultats sur les espaces des
modules de fibr\'es vectoriels afin de pr\'eciser le contexte dans lequel se place la dite conjecture. \\

Soit $A$ une vari\'et\'e ab\'elienne complexe de dimension $g$ et soit $L$ un fibr\'e en droites sur $A$
qui d\'efinit une polarisation principale. Pour un entier fixe $d \geq 1$
on note $V_d= H^0(A, L^{\otimes d})$. On consid\`ere $\varphi_d : A \rightarrow \mathbb{P}(V_d)$ \footnote{Ici $\mathbb{P}(V_d)$ d\'esigne l'espace
d'hyperplans de $V_d$} l'application d\'efinie par les sections globales de
$L^{\otimes d}$. Rappelons que pour $d\geq 3$, $\varphi_3$ est un plongement, tandis que $\varphi_2$ induit un plongement 
de la
vari\'et\'e de Kummer $ A/ \{ \pm 1\}$ dans $\mathbb{P} (V_2)$, lorsque $(A,L)$ est ind\'ecomposable.  Soit $A[d]$ le sous-groupe fini des points d'ordre $d$ en $A$; $A[d]$ agit sur $A$ en
preservant le fibr\'e $L^{\otimes d}$, donc agit  sur $\mathbb{P}(V_d)$ de mani\`ere que $\varphi_d$ est $A[d]$-\'equivariant.

\begin{Teo} (A. Coble ) \label{coble} \\
a) Soit $g=2$. Il existe une unique hypersurface cubique $A[3]$-invariante dans $\mathbb{P}(V_3)\simeq \mathbb{P}^8$ singuli\`ere le long de
$\varphi_3(A)$. Les polaires de cette cubique engendrent l'espace des quadriques dans   $\mathbb{P}(V_3)$ contenant $\varphi_3(A)$ ~\cite{C1}.\\
b) Soit $g=3$. Il existe une unique hypersurface quartique $A[2]$-invariante dans $\mathbb{P}(V_2)\simeq \mathbb{P}^7$ singuli\`ere le long de
$\varphi_2(A)$. Les polaires de cette quartique engendrent l'espace des cubiques dans $\mathbb{P}(V_2)$ contenant $\varphi_2(A)$~\cite{C2}.\\
\end{Teo}
Comme on le verra par la suite, il y a un lien \'etroit  entre les observations faites par Coble il y a presque un si\`ecle et quelques r\'esultats, beaucoup
plus r\'ecents, concernant les espaces des modules de fibr\'es vectoriels.\\

Soit  $C$ une courbe lisse, projective et connexe de genre $g \geq 2$. Rappelons que la jacobienne $JC$ param\`etre les
fibr\'es en droites
de degr\'e 0 sur $C$.  On consid\`ere $J^{g-1}$ (vari\'et\'e isomorphe \`a $JC$) la vari\'et\'e qui param\`etre les
fibr\'es en droites de degr\'e $g-1$ sur $C$.
Cette vari\'et\'e contient canoniquement le diviseur th\^eta d\'efini ensemblistement par
\begin{displaymath}
\Theta = \{ M \in J^{g-1}  \mid h^0(C, M) \geq 1 \}.
\end{displaymath}
Soit $\mathcal{SU}_C(r)$ l'espace des modules des fibr\'es vectoriels semi-stables sur $C$ de rang $r$ et de d\'eterminant trivial. C'est une vari\'et\'e
projective irr\'eductible, dont le lieu singulier consiste exactement en les points non-stables, sauf si $r=2, g=2$. Dans ce cas
$\mathcal{SU}_C(2)$ est lisse.
On d\'efinit une application rationnelle
\begin{displaymath}
\theta: \mathcal{SU}_C(r) \dashrightarrow  |r\Theta|
\end{displaymath}
telle que $\theta(E)$ est le diviseur de support
\begin{displaymath}
\Theta_E:= \{ L \in J^{g-1}  \mid h^0(C,E \otimes L) \geq 1 \}.
\end{displaymath}
Pour certaines valeurs de $r$ et de $g$ il existe des \'el\'ements $E$ dans  $\mathcal{SU}_C(r)$ tels que $\Theta_E= J^{g-1}$. M. Raynaud ~\cite{R} donne des exemples de tels fibr\'es
vectoriels et des conditions sous lequelles $\Theta_E$ d\'efinit un diviseur dans $J^{g-1}$. En particulier, les r\'esultats de Raynaud montrent que
l'application $\theta$ est bien d\'efinie pour $r=2$ en genre quelconque et pour $r=3$ et $g=2$, ce qui est le cas qui nous int\'eresse.\\

Lorsque $r=2$ le syst\`eme lin\'eaire $|2\Theta|$ est particuli\`erement int\'eressant parce qu'il contient la vari\'et\'e de Kummer $\mathcal{K}_C$, i.e.
le quotient de la jacobienne $JC$ par l'involution $a \mapsto -a$.
Donc l'application
$$
\begin{array} {rcl}
JC & \rightarrow  & |2\Theta| \\
a & \mapsto & \Theta_a + \Theta_{-a}
\end{array}
$$
factorise par le plongement $\kappa : \mathcal{K}_C \hookrightarrow |2\Theta|$. La partie non-stable de $\mathcal{SU}_C(2)$ est form\'ee des fibr\'es
vectoriels de
la forme $M \oplus M^{-1}$, avec $M \in JC$ et donc elle s'identifie \`a $\mathcal{K}_C$. Pour $g \geq 3$ ces fibr\'es forment le lieu
singulier de $\mathcal{SU}_C(2)$. On obtient ainsi le diagramme commutatif
\begin{eqnarray*}
\shorthandoff{;:!?}
\xymatrix{
	\mathcal{K}_C   \ar@{^{(}->}[d] \ar[rd]^{\kappa} \\
	\mathcal{SU}_C(2)  \ar[r]^{\theta} & |2\Theta|
     }
\end{eqnarray*}
On r\'esume dans le th\'eor\`eme suivant les r\'esultats connus sur l'application $\theta$ pour les fibr\'es de rang 2.
\begin{Teo} \textnormal{    }\\
a) Pour $g=2$, $\theta$ est un isomorphisme de $\mathcal{SU}_C(2)$ sur $|2\Theta| \simeq \mathbb{P}^3$  ~\cite{N-R1}. \\
b) Pour $g\geq 3$ et $C$ hyperelliptique, $\theta$ est fini de degr\'e 2 sur une sous-vari\'et\'e de $|2\Theta|$ qui est d\'ecrite explicitement dans
~\cite{D-R}.\\
c) Pour $g \geq 3$ et $C$ non-hyperelliptique, $\theta$ est un plongement ~\cite{vG-I}.
\end{Teo}
Dans le cas o\`u $C$ est une courbe non-hyperelliptique de genre 3, Narasimhan et Ramanan ~\cite{N-R2} ont d\'emontr\'e que $\theta$ est un
isomorphisme de $\mathcal{SU}_C(2)$ dans une hypersurface quartique $\mathcal{Q}$ dans $|2\Theta|\simeq \mathbb{P}^7$.  Par la remarque ci-dessus,
cette quartique est singuli\`ere le long de la vari\'et\'e de Kummer $\mathcal{K}_C$.
Par ailleurs, par le th\'eor\`eme \ref{coble} {\it a)} il existe une \emph{unique} hypersurface quartique $A[2]$-invariante dans $|2\Theta|$ qui est singuli\`ere le long de
$\mathcal{K}_C$.  Ainsi la quartique $\mathcal{Q}$ est justement la quartique de Coble.
\\

Dans ce travail, on consid\`ere une courbe de genre 2 et  l'espace des modules $\mathcal{SU}_C(3)$.
Par analogie avec le r\'esultat de Narasimhan et Ramanan, le th\'eor\`eme principal de la th\`ese \'etablit le lien entre l'espace de modules
$\mathcal{SU}_C(3)$ et l'hypersurface cubique dans $\mathbb{P}^8$ donn\'ee par le th\'eor\`eme \ref{coble} {\it b)}. On a appel\'e cette hypersurface
\emph{cubique de Coble}.
Dans cette situation l'application
\begin{displaymath}
\theta : \mathcal{SU}_C(3) \rightarrow |3\Theta|\simeq \mathbb{P}^8,
\end{displaymath}
est bien d\'efinie et de degr\'e 2.  On note $\mathfrak{i}$ l'involution dans $\mathcal{SU}_C(3)$ d\'efinie par
\begin{displaymath}
\mathfrak{i} : E \mapsto \iota^*E^*,
\end{displaymath}
o\`u $\iota$ est l'involution hyperelliptique sur $C$. On d\'emontre que le lieu de ramification de $\theta$ est \'egal aux points fix\'es par
cette involution, autrement dit,  il est \'egal \`a
\begin{displaymath}
\{ E \in \mathcal{SU}_C(3) \mid E \sim_{s} \iota^*E^* \}.
\end{displaymath}
Soit $\mathcal{B}$ l'image de cette hypersurface dans $\mathbb{P}^8$. On pose $A=JC$. On consid\`ere 
la vari\'et\'e $J^1=\mathrm{Pic}^1(C)$ et le plongement $\varphi_3 : J^1 \rightarrow |3\Theta|^*$. Le groupe $A[3]$ agit 
sur $J^1$ et $|3\Theta|^*$ de mani\`ere que $\varphi_3$ est $A[3]$-\'equivariant. Par le th\'eor\`eme \ref{coble} {\it b)}
il existe une \emph{unique} cubique $\mathcal{C} \subset \mathbb{P}^{8*}$, $A[3]$-invariante et singuli\`ere le long de
$J^1$. On \'enonce le r\'esultat principal de la th\`ese comme suit
\begin{Teo}
La vari\'et\'e duale de la cubique de Coble $\mathcal{C}$ est l'hypersurface sextique $\mathcal{B}$.
\end{Teo}
\begin{Rmq}
La vari\'et\'e $J^1$ (isomorphe \`a $JC$) n'est pas une vari\'et\'e ab\'elienne, mais elle contient canoniquement le
diviseur th\^eta $\Theta$, c'est pourquoi on a choisi de consid\'erer cette vari\'et\'e au lieu de la jacobienne $JC$. 
\end{Rmq}
\\
\\
{\it Esquisse de la d\'emonstration.}\\

L'outil clef de la preuve est l'action du groupe de Heisenberg $\He$, puisque l'on utilise les restrictions aux plans des points fixes par un \'el\'ement
du groupe de Heisenberg pour comparer les hypersurfaces $\mathcal{C}^*$ et  $\mathcal{B}$. On notera $\mathbb{P}^2_{\tilde{\eta}}$ le plan des
points fix\'es par l'\'el\'ement $\tilde{\eta} \in \He$. \\
Tout d'abord on d\'ecrit le lieu de ramification de l'application $\theta$ en termes de l'involution $\mathfrak{i}$ sur $\mathcal{SU}_C(3)$. On calcule le
degr\'e de la vari\'et\'e duale $\mathcal{C}$, donn\'ee n\'ecessaire pour montrer que $\mathcal{C}^*
\cap \mathbb{P}^2_{\tilde{\eta}}= (\mathcal{C} \cap \mathbb{P}_{\tilde{\eta}}^2)^*$. Les intersections de $\mathcal{C}^*$ avec les plans des
points fixes sont donc des sextiques planes duales des cubiques $\mathcal{C} \cap \mathbb{P}_{\tilde{\eta}}^2$.\\
D'un autre c\^ot\'e, pour montrer  que les sextiques planes $\mathcal{B} \cap \mathbb{P}^2_{\tilde{\eta}}$ sont en fait des courbes duales de certaines cubiques,
on fait appel au th\'eor\`eme sur les vari\'et\'es de Prym  expos\'e dans le chapitre 3 de la th\`ese et prouv\'e dans un cas
plus g\'en\'eral.\\
Le pas suivant consiste \`a d\'efinir deux applications sur les s\'ecantes de $J^1$, une avec image contenue
dans le lieu singulier de $\mathcal{C}$ et l'autre avec image contenue dans le lieu singulier de $\mathcal{B}$.  On d\'emontre que ces deux applications
co\"{\i}ncident. Ceci et le fait que ces sextiques planes sont compl\`etement caract\'eris\'ees par leurs points de rebroussement,
est suffisant pour montrer que, pour tout $\tilde{\eta} \in \He$ non nul, $\mathcal{C}^*\cap \mathbb{P}^2_{\tilde{\eta}}$ et $\mathcal{B} \cap
\mathbb{P}^2_{\tilde{\eta}}$ est la m\^eme sextique. \\
L'id\'ee est de prouver que les intersections d'une sextique $A[3]$-invariante avec les plans $\mathbb{P}^2_{\tilde{\eta}}$
la caract\'erisent compl\`etement.  Pour cela, on commence par donner une base explicite
de $(S^6V)^{A[3]}$, l'espace vectoriel des sextiques $A[3]$-invariantes dans $\mathbb{P}^8$. On consid\`ere
l'application $\nu$ qui envoie une sextique dans $(S^6V)^{A[3]}$ en la somme directe de ses restrictions aux plans $\mathbb{P}_{\tilde{\eta}}^2$.
\`A l'aide du logiciel Maple on calcule le rang de cette application. Il s'av\`ere qu'elle n'est pas injective. Cependant, on prouve que la diff\'erence (\`a scalaire pr\`es) des polyn\^omes
qui d\'efinissent  les sextiques $\mathcal{C}^*$ et $\mathcal{B}$ est invariante par $\iota$, l'involution induite par l'involution hyperelliptique, et que le noyau de $\nu$ est $\iota$-anti-invariant.
On conclut ainsi l'\'egalit\'e de ces hypersurfaces.
\\
\\
{\it Organisation de la th\`ese.}\\

La th\`ese est compos\'ee de trois chapitres. Dans le premier chapitre on rappelle quelques d\'efinitions et  r\'esultats sur le groupe de Heisenberg
et son action sur $H^0(A, \mathcal{O}_A(3\Theta))$, qu'on utilise dans la deuxi\`eme partie. On pr\'esente \'egalement l'\'equation explicite de la
cubique de Coble, utile pour certaines preuves, et on explique son origine.\\
Le chapitre 2 est enti\`erement consacr\'e \`a la d\'emonstration de notre th\'eor\`eme. Finalement, dans le troisi\`eme
chapitre on expose un r\'esultat sur la
d\'ecompositon des vari\'et\'es de Prym associ\'ees aux rev\^etements $n$-cycliques sur une courbe hyperelliptique, 
utilis\'e dans la preuve de la conjecture.\\

\mainmatter

\chapter{ Pr\'eliminaires}

\section{Rappels}

On trouvera plus de d\'etails sur les faits qu'on rappelle ici, dans ~\cite{CAV},~\cite{Mum1}, \cite{Mum2} et
\cite{Mum3}.\\

\subsection{Vari\'et\'es ab\'eliennes}
Soit $V$ un $\mathbb{C}$-espace vectoriel de dim $g$ et $\Lambda \simeq \mathbb{Z}^{2g}$ un r\'eseau dans $V$. Une 
\emph{polarisation} sur le tore complex $A=V/\Lambda$ est par d\'efinition la premi\`ere classe de Chern $H=c_1(L)$
d'un fibr\'e en droites ample $L$ sur $A$. Parfois on consid\`ere le fibr\'e en droites lui-m\^eme comme la polarisation.\\
La polarisation $H$ est une forme hermitienne sur $V$ dont la forme altern\'ee asoci\'ee $E=\im H$ prend des valeurs 
enti\`eres sur $\Lambda$. Il existe une base $\lambda_1, \ldots, \lambda_g,\mu_1, \ldots , \mu_g$ de $\Lambda$, par 
rapport \`a laquelle $E$ est donn\'ee par la matrice
$$\left(
\begin{array}{cc}
0 & D\\
-D & 0
\end{array} \right)
$$
o\`u $D=\mathrm{diag}(d_1, \ldots , d_g)$, avec $d_{\nu}$ entiers positifs v\'erifiant $d_{\nu} \mid d_{\nu +1}$ pour
$1\leq \nu \leq g-1 $. On appelle le vecteur    
$(d_1, \ldots , d_g)$ le \emph{type de la polarisation}. Une polarisation est principale si elle est de type $(1, 
\ldots , 1)$. Une \emph{vari\'et\'e ab\'elienne} est un tore complexe $A$ qui admet une polarisation. Ainsi une 
vari\'et\'e ab\'elienne principalement polaris\'ee (v.a.p.p.) est un couple $(A,L)$ o\`u $L$ est une polarisation
principale.
\\
Soit $C$ une courbe lisse non singuli\`ere de genre $g$. On note $H^0(\omega_C)$ le $\mathbb{C}$-espace vectoriel des
1-formes holomorphes sur $C$. Le groupe d'homologie $H_1(C, \mathbb{Z})$ est un groupe libre abelien de rang $2g$.
L'application injective $H_1(C, \mathbb{Z}) \rightarrow H^0(\omega_C)^*$ donn\'ee par $\gamma \mapsto \int_{\gamma}$ nous
permet de consid\'erer $H_1(C, \mathbb{Z})$ comme un r\'eseau de $H^0(\omega_C)^*$. La \emph{vari\'et\'e jacobienne de
$C$}, d\'efinie par 
\begin{displaymath}
JC:= H^0(\omega_C)^*/ H_1(C, \mathbb{Z}),
\end{displaymath}
est le premier exemple d'une vari\'et\'e ab\'elienne qui admet naturellement une polarisation principale $\mathcal{O}_{JC}(\Theta)$. 

\subsection{Vari\'et\'es de Prym}

Soient $\widetilde{C}$ et $C$ des courbes compl\`etes non-singuli\`eres avec des jacobiennes $\widetilde{J}$ et $J$ 
respectivement. Soit $\pi : \widetilde{C} \rightarrow C$ un rev\^etement double et $i: \widetilde{C} \rightarrow
 \widetilde{C}$ l'involution qui \'echange les feuilles au-dessus de chaque point de $C$. On d\'efinit l'application 
\emph{Norme de $\pi$}, en termes de diviseurs, par 
$$
\begin{array}{rccl}
\Nm_{\pi} : & \widetilde{J} & \rightarrow & J \\
&\sum n_i p_i & \mapsto & \sum n_i \pi (p_i),
\end{array}
$$
o\`u $p_i\in \widetilde{C}$ et $n_i \in \mathbb{Z}$.
Pour tout diviseur $\mathfrak{D}$ dans $\widetilde{C}$ on a 
$\pi^{-1}(\pi \mathfrak{D}) = \mathfrak{D} + i(\mathfrak{D})$,
d'o\`u 
\begin{displaymath}
\pi^{*}(\Nm_{\pi} x) = x + i(x), \qquad \forall x \in \widetilde{J}.
\end{displaymath}
Comme l'application $\Nm_{\pi}$ est surjective, ceci montre que
\begin{displaymath}
i(\pi^{*} y) = \pi^{*} y, \qquad \forall y \in J.
\end{displaymath}
Donc $i_{|_{\pi^* J}} = +1$ et $i_{|_{\Ker \Nm_{\pi}}} = -1$. On d\'efinit la \emph{vari\'et\'e de Prym} de $\widetilde{C}$
sur $C$ par 
\begin{displaymath}
P:= (\Ker \Nm_{\pi})^0 = \im (1_{\widetilde{J}}-i).
\end{displaymath}
Cette vari\'et\'e est une sous-vari\'et\'e ab\'elienne de $\widetilde{J}$ (la partie "impaire" de $\widetilde{J}$).
Dans ~\cite{Mum2} D. Mumford \'etablit les cas o\`u $P$ admet une polarisation principale.\\

Dans le chapitre 3 on consid\'erera une notion de vari\'et\'es de Prym \'elargie: on consid\`ere un rev\^etement \'etale
$n$-cyclique $\pi: \widetilde{C} \rightarrow C$ et on d\'efinit la vari\'et\'e de Prym comme la composante connexe du noyau
de la Norme $\Nm_{\pi}$, qui contient $0 \in \widetilde{J}$. \`A diff\'erence des vari\'et\'es consid\'er\'ees par 
Mumford, celles-ci ne sont pas principalement polaris\'ees.

\subsection{Groupe de Heisenberg}

Soit $(A,L)$ une vari\'et\'e ab\'elienne principalement polaris\'ee.
Le \emph{groupe th\^eta de $L$} de niveau 3 est le groupe
$$\mathcal{G}(L)= \{ (\varphi, \eta ) \mid \eta \in A,\ \varphi: t_{\eta}^*(L^3) \stackrel{\sim}{\rightarrow} L^3 \}, $$
avec l'op\'eration du groupe donn\'ee par
\begin{displaymath}
(\varphi, \eta )\cdot (\varphi', \eta' ) = (t^*_{\eta'} \varphi \circ \varphi', \eta + \eta') .
\end{displaymath}
On a l'extension centrale de groupes
\begin{displaymath}
1 \longrightarrow \mathbb{C}^* \stackrel{i}{\longrightarrow}   \mathcal{G}(L)  \stackrel{p}{\longrightarrow} A[3] \longrightarrow 0,
\end{displaymath}
avec $i(\alpha)=(\alpha ,0)$ et $p(\varphi , \eta )=\eta$ et $A[3]$ le sous-groupe fini des points d'ordre 3.
Le commutateur $[(\varphi, \eta ),(\varphi', \eta') ]$ de deux \'el\'ements dans $\mathcal{G}(L)$ appartient au centre, et induit la \emph{ forme de Weyl}
\begin{displaymath}
e^L : A[3] \times A[3]  \rightarrow  \mathbb{C}^*.
\end{displaymath}
Comme groupe abstrait $\mathcal{G}(L)$ est isomorphe au groupe de Heisenberg
\begin{displaymath}
\mathcal{H}(3) := \mathbb{C}^* \times (\mathbb{Z}/3)^g \times (\widehat{\mathbb{Z}/3})^g \simeq \mathbb{C}^* \times 
(\mathbb{Z}/3)^{2g},
\end{displaymath}
o\`u $(\widehat{\mathbb{Z}/3})^g:= \Hom((\mathbb{Z}/3)^g, \mathbb{C}^*)$,  avec la loi de groupe
\begin{displaymath}
(t,x,x^*)\cdot(s,y,y^*) = (st (y^*(x)- x^*(y) ), x+y, x^* y^* ).
\end{displaymath}
On note $\langle (x,x^*),(y,y^*) \rangle : = y^*(x) - x^*(y)$ la forme bilin\'eaire sur  $(\mathbb{Z}/3)^g \times 
(\widehat{\mathbb{Z}/3})^g$.
La projection $(t,x,x^*)\mapsto (x,x^*)$ d\'efinit une extension  centrale de groupes
\begin{displaymath}
1 \longrightarrow \mathbb{C}^* \longrightarrow   \mathcal{H}(3)  \longrightarrow (\mathbb{Z}/3)^{2g} \longrightarrow 0.
\end{displaymath}
Une \emph{structure th\^eta} de niveau 3 pour $L$ est un isomorphisme
\begin{displaymath}
\alpha: \mathcal{H}(3) \stackrel{\sim}{\rightarrow} \mathcal{G}(L).
\end{displaymath}
Par projection sur $(\mathbb{Z}/3)^{2g}$, une structure th\^eta $\alpha$ (de niveau 3) induit un isomorphisme
\begin{displaymath}
\tilde{\alpha}: (\mathbb{Z}/3)^{2g} \stackrel{\sim}{\rightarrow} A[3],
\end{displaymath}
tel que le diagramme suivant commute
$$
\shorthandoff{;:!?}
\xymatrix @!0 @R=1.3cm @C=1.8cm {
        1  \ar[r] & \mathbb{C}^* \ar[r] \ar@{=}[d] &  \mathcal{G}(L) \ar[r]  & A[3] \ar[r] & 0  \\
         1  \ar[r] & \mathbb{C}^* \ar[r]  &  \mathcal{H}(3) \ar[r] \ar[u]_{\alpha}  & (\mathbb{Z}/3)^{2g}  \ar[r]  \ar[u]_{\tilde{\alpha}} & 0  \\
     }
$$
De plus, l'isomorphisme $\tilde{\alpha}$ est un isomorphisme symplectique par rapport aux formes bilin\'eaires $e^L$ et  
$\langle \cdot, \cdot \rangle$.
On appelle \emph{structure de niveau 3} un isomorphisme symplectique $\tilde{\alpha}: (\mathbb{Z}/3)^{2g} \rightarrow  A[3] $.
\\
Consid\'erons le demi-espace de Siegel 
$$\mathcal{S}_g=\{ \tau \in M_g(\mathbb{C}) \mid \tau {}^t\tau, \im \tau >0 \}$$ 
et le groupe modulaire de Siegel
$$  \Gamma_g =  \left\{ M={A \ B \choose C \ D}  \in
    GL(2g,\mathbb{Z}) \mid {}^{t}M {0 \ -\mathrm{I}_{g} \choose \mathrm{I}_{g}  \ \ 0} M  = {0 \ -\mathrm{I}_{g} \choose \mathrm{I}_{g} \ \ 0} \right\}, $$
qui agit sur $\mathcal{S}_g$ par 
\begin{displaymath}
M: \tau \mapsto M_{\tau}:= (A\tau +B)(C\tau +D)^{-1}.
\end{displaymath}
Soit  
$$\Gamma_{g}(3) = \{ M \in \Gamma_g  \mid  M  \equiv  I_{2g} \textrm{ mod } 3 \} $$
un sous-groupe de $\Gamma_g$.
Deux \'el\'ements $\tau$ et $\tau'$ dans $\mathcal{S}_g$ d\'efinissent des v.a.p.p. avec structure de niveau 3 si et seulement s'il existe $M \in \Gamma_g(3)$
tel que  $\tau' = M_{\tau}$. La vari\'et\'e complexe ${\textit{\LARGE a}}_g(3) := \mathcal{S}_g / \Gamma_g(3)$  
param\`etre donc les classes
d'isomorphismes de v.a.p.p. de dimension $g$ avec structure de niveau 3.
\\
\\
Le groupe $\mathcal{G}(L)$ a une repr\'esentation naturelle dans l'espace $V_A:= H^0(A, L^3)$, donn\'ee par
\begin{displaymath}
((\varphi,\eta)s) (a) = \varphi_a(s(a+\eta)) ,
\end{displaymath}
o\`u $a\in A$ et $\varphi_a : t_{\eta}^*(L^3)_a \simeq (L^3)_{a+\eta} \rightarrow (L^3)_a $.
\\
D'un autre c\^ot\'e, si on note $V(g)$ le $\mathbb{C}$-espace vectoriel de fonctions sur $(\mathbb{Z}/3)^{g}$ \`a valeurs complexes,
le groupe $\mathcal{H}(3)$ agit lin\'eairement sur $V(g)$ par
\begin{displaymath}
((t,x,x^*) f ) (v) = t x^*(v-x) f(v-x), \quad v\in (\mathbb{Z}/3)^g.
\end{displaymath}
Ceci d\'efinit une repr\'esentation de $\mathcal{H}(3)$ connue dans la litt\'erature comme la repr\'esentation de Schr\"odinger.
Toutes les repr\'esentations irr\'eductibles du groupe de Heisenberg
o\`u le centre agit par multiplication sont isomorphes (proposition 3 ~\cite{Mum1}).
Il existe donc un isomorphisme de repr\'esentations lin\'eaires
\begin{displaymath}
    \psi_{\alpha}: V(g) \rightarrow V_{A}.
\end{displaymath}
Par le lemme de Schur cet isomorphisme est unique \`a scalaire
pr\`es.\\
On consid\`ere le sous-groupe de $\mathcal{H}(3)$
\begin{displaymath}
    \He:= \{ (t,x,x^*) \in \mathcal{H}(3) \mid
    t^{3}=1 \}.
\end{displaymath}
C'est une extension centrale non-triviale
\begin{displaymath}
    1 \longrightarrow {\bf\mu_3}\longrightarrow \He
    \longrightarrow (\mathbb{Z}/3)^{2g} \longrightarrow 1,
\end{displaymath}
o\`u ${\bf\mu_{3}}$ est le groupe des racines cubiques de l'unit\'e.
Par restriction, $V(g)$ est une repr\'esentation lin\'eaire de
$\He$. Pour chaque $v \in (\mathbb{Z}/3)^g$ soit
$\delta_{v} \in V(g)$ la fonction caract\'eristique
$$ \delta_{v}= \left\{
\begin{array}{cc} 0 & \textnormal{ si } \ x\neq v \\ 1 &\textnormal{ si } \ x=v. \end{array}
\right.$$
On a donc
\begin{displaymath}
  (t,x,x^{*}) \delta_{v}= tx^{*}(v-x)\delta_{v-x},
\end{displaymath}
pour tout $ (t,x,x^{*}) \in \He$ . Apr\`es avoir fix\'e
un ordre dans $(\mathbb{Z}/3)^{g}$, on obtient une base canonique $\{
\delta_{v}\}_{v\in (\mathbb{Z}/3)^g}$  de $V(g)$. L'image de cette
base par l'isomorphisme $\psi_{\alpha}$ correspond \`a la base de
fonctions th\^eta 
\begin{displaymath}
   \left\{ X_{b} := \vartheta \left[ \begin{array}{c} b/3 \\ 0 \end{array}
    \right] (3z; 3\tau), \quad b\in (\mathbb{Z}/3)^{g}, \ z\in
    \mathbb{C}^g \right\}
\end{displaymath}
de l'espace $V_A $, avec  $\tau \in \mathcal{S}$ fix\'e. On consid\`ere $\varphi=\varphi_{\tau} : A
\rightarrow \mathbb{P}(V_{A}^*)$ l'application d\'efinie par
\begin{displaymath}
  z \mapsto (\ldots , X_{b}(z), \ldots)
\end{displaymath}
correspondant au syst\`eme lin\'eaire $|L^{3}|$. Soit
$\varphi_{\alpha}$ la compos\'ee des applications
\begin{displaymath}
 A \stackrel{\varphi}{\longrightarrow} \mathbb{P}(V_{A}^*)
 \stackrel{{}^t\psi_{\alpha}}{\longrightarrow} \mathbb{P}(V(g)^*).
\end{displaymath}
On observe que le groupe $A[3] = \mathcal{G}(L)/\mathbb{C}^*$ agit par
translations sur $A$ et le groupe
\begin{displaymath}
 \He / {\bf\mu_{3}} = \mathcal{H}(3) / \mathbb{C}^* =
 (\mathbb{Z}/3)^{2g}
\end{displaymath}
agit sur $ \mathbb{P}(V(g)^*)$ par le projectivis\'e de la
repr\'esentation duale de $\He$ sur $V(g)^*$.  Ces
actions sont compatibles dans le sens o\`u elles v\'erifient la relation
\begin{displaymath}
 \varphi_{\alpha}(a+ \eta) = \tilde{\alpha}^{-1}(\eta)
 \varphi_{\alpha} (a), \quad\forall \eta \in A[3], \quad a\in A.
\end{displaymath}
Explicitement, l'action de $\He$ sur la base $\{ X_b \}$ est donn\'ee par
\begin{displaymath}
(t,x,x^*) X_b = tx^*(b-x)X_{b-x},
\end{displaymath}
pour tout $(t,x,x^*) \in \He$ et $b \in (\mathbb{Z}/3)^2$.
\\

Dans la deuxi\`eme partie on utilisera seulement l'action du groupe de
Heisenberg \emph{fini} $\He$ et on le nommera simplement  groupe de
Heisenberg.

\section{La cubique de Coble} \label{coble's_cubic}

Soit $g=2$. \'Etant donn\'ee une surface ab\'elienne $(A,L,\alpha)$,
principalement polaris\'ee avec une structure th\^eta $\alpha$,
le r\'esultat de Coble assure l'existence d'une unique hypersurface
cubique $A[3]$-invariante dans $\mathbb{P}(V_{A}^*) \simeq
\mathbb{P}^8$, singuli\`ere le long de $\varphi(A)$. On cherche \`a
expliciter l'\'equation de cette cubique. Pour cela on utilisera la
base de fonctions th\^eta $\{ X_{b}\}$ de
$V_{A}$. \\

Dans ~\cite{B3} (proposition 1), A. Beauville d\'eduit que la dimension
de l'espace des formes de degr\'e 3 $\He$-invariantes
dans $V_{A}$ est \'egale \`a 5. Soient $ K := \{ (1,x,0) \mid x \in
(\mathbb{Z}/3)^2 \}$ et $\widehat{K} :=\{(1,0,x^{*}) \mid x^* \in
(\widehat{\mathbb{Z}/3})^2 \}$ sous-groupes maximaux de
$\He$. Un mon\^ome de la forme
$X_{\sigma_{1}}X_{\sigma_{2}}X_{\sigma_{3}}$, avec $\sigma_{i}\in
(\mathbb{Z}/3)^2$ est $\widehat{K}$-invariant si et seulement si
$\sigma_{1}+\sigma_{2}+\sigma_{3}=0$. En faisant agir les \'el\'ements
de $K$ sur ce type de mon\^omes et en prenant la somme, on obtient une
base de formes cubiques $\He$-invariantes comme suit
\begin{eqnarray*}
    F_0 &:= &\sum_{b \in (\mathbb{Z}/3)^2} X_{b}^3 \\
    F_1 &:= &\sum_{b} X_{b}X_{01 + b}X_{02+b} \\
    F_2 &:= &\sum_{b} X_{b}X_{10 + b}X_{20+b} \\
    F_3 &:= &\sum_{b} X_{b}X_{11 + b}X_{22+b} \\
    F_4 &:= &\sum_{b} X_{b}X_{12 + b}X_{21+b}.
\end{eqnarray*}
Une cubique g\'en\'erale dans $(S^3V_{A})^{\He}$ est
donc de la forme
\begin{equation} \label{cubique}
F_{\beta}:= \beta_{0}F_{0} + \beta_{1}F_{1} + \beta_{2}F_{2}
+\beta_{3}F_{3} +\beta_{4}F_{4},
\end{equation}
avec $\beta= (\beta_{0},\beta_{1},\beta_{2},\beta_{3},\beta_{4})$ des
coefficients complexes.\\

Par ailleurs, comme $\varphi(A)$ est projectivement normal dans $\mathbb{P}(V_A^*)\simeq \mathbb{P}^8$ ~\cite{K}, l'application 
$H^0(\mathbb{P}^8,\mathcal{O}(2))\rightarrow H^0(A,\mathcal{O}_A(2)) $ est surjective et donc la dimension de l'espace de 
quadriques dans $\mathbb{P}^{8}$ qui contiennent $\varphi(A)$
est \'egale \`a
\begin{eqnarray*}
    \dim H^0(\mathbb{P}^{8}, \mathcal{I}_{A}(2)) &=& \dim
    H^0(\mathbb{P}^8, \mathcal{O}(2)) - \dim H^0(A,
    \mathcal{O}_A(6\Theta))\\
    &=& {10 \choose 2} - 6^2\\
    &=& 9
\end{eqnarray*}
o\`u $\mathcal{I}_{A}$ d\'esigne le faisceau d'id\'eaux correspondant au
plongement $\varphi: A \hookrightarrow \mathbb{P}^8$. Ainsi, par la
proposition 1~\cite{B3}, les d\'eriv\'ees partielles
$\left\{ \frac{\partial F_{\beta}}{ \partial X_b}\right\}_{b \in
(\mathbb{Z}/3)^2}$ forment une base pour l'espace des hypersufaces
quadriques contenant $A$. Coble avait conjectur\'e que ces quadriques
suffisaient pour d\'efinir la surface $A$. Ce fait a \'et\'e
d\'emontr\'e par W. Barth dans ~\cite{Barth} lorsque $A$ est ind\'ecomposable. Explicitement, $A$ est
d\'efinie comme l'intersection des quadriques suivantes ~\cite{Barth}
\begin{eqnarray*}
    Q_{00} &:= &\beta_{0}X_{00}^2 + \beta_{1}X_{01}X_{02} +
    \beta_{2}X_{10}X_{20} + \beta_{3}X_{11}X_{22} + \beta_{4}X_{12}X_{21}   \\
    Q_{01} &:= &\beta_{0}X_{01}^2 + \beta_{1}X_{02}X_{00} +
    \beta_{2}X_{11}X_{21} + \beta_{3}X_{12}X_{20} +
    \beta_{4}X_{10}X_{22}\\
    Q_{02} &:= &\beta_{0}X_{02}^2 + \beta_{1}X_{00}X_{01} +
    \beta_{2}X_{12}X_{22} + \beta_{3}X_{10}X_{21} +
    \beta_{4}X_{11}X_{20} \\
    Q_{10} &:= &\beta_{0}X_{10}^2 + \beta_{1}X_{11}X_{12} +
    \beta_{2}X_{20}X_{00} + \beta_{3}X_{21}X_{02} + \beta_{4}X_{22}X_{01} \\
    Q_{11} &:= &\beta_{0}X_{11}^2 + \beta_{1}X_{12}X_{10} +
    \beta_{2}X_{21}X_{01} + \beta_{3}X_{22}X_{00} +
    \beta_{4}X_{20}X_{02}\\
     Q_{12} &:= &\beta_{0}X_{12}^2 + \beta_{1}X_{10}X_{11} +
    \beta_{2}X_{22}X_{02} + \beta_{3}X_{20}X_{01} +
    \beta_{4}X_{21}X_{00}\\
    Q_{20} &:= &\beta_{0}X_{20}^2 + \beta_{1}X_{21}X_{22} +
    \beta_{2}X_{00}X_{10} + \beta_{3}X_{01}X_{12} + \beta_{4}X_{02}X_{11}   \\
    Q_{21} &:= &\beta_{0}X_{21}^2 + \beta_{1}X_{22}X_{20} +
    \beta_{2}X_{01}X_{11} + \beta_{3}X_{02}X_{10} +
    \beta_{4}X_{00}X_{12}\\
     Q_{22} &:= &\beta_{0}X_{22}^2 + \beta_{1}X_{20}X_{21} +
    \beta_{2}X_{02}X_{12} + \beta_{3}X_{00}X_{11} +
    \beta_{4}X_{01}X_{10}.
\end{eqnarray*}
On donnera un bref aper\c cu du lien entre les param\`etres $\beta=(\beta_{0},\beta_{1},\beta_{2},\beta_{3},\beta_{4})$
et la surface ab\'elienne $A$ qui est l'intersection des quadriques ci-dessus. On consid\`ere l'involution 
$\iota : a \mapsto -a$ sur $A$ qui agit sur les coordonn\'ees de $V_A$ par 
\begin{displaymath}
\iota \cdot X_b = X_{-b}, \quad b \in (\mathbb{Z}/3)^2.  
\end{displaymath}
Cette action induit une d\'ecomposition de $\mathbb{P}^8$ dans un espace $\mathbb{P}^4_+$ $\iota$-invariante et un espace
$\mathbb{P}^3_-$ $\iota$-anti-invariante. Soient
\begin{eqnarray*}
Y_0 := X_{00}  \ \quad \quad \qquad  & & \\
Y_1:= \frac{1}{2}(X_{01}+ X_{02}) & & Z_1:= \frac{1}{2}(X_{01}- X_{02})\\
Y_2:= \frac{1}{2}(X_{10}+ X_{20}) & & Z_2:= \frac{1}{2}(X_{10}- X_{20})\\
Y_3:= \frac{1}{2}(X_{11}+ X_{22}) & & Z_3:= \frac{1}{2}(X_{11}- X_{22})\\
Y_4:= \frac{1}{2}(X_{12}+ X_{21}) & & Z_4:= \frac{1}{2}(X_{12}- X_{21})\\
\end{eqnarray*}
des coordonn\'ees de $\mathbb{P}^4_+$ et  $\mathbb{P}^3_-$ respectivement. Les 16 points de 2-torsion de $A$ sont les seuls points fixes par l'involution $\iota$, donc
\begin{eqnarray*}
A \cap \mathbb{P}^3_+ &=& \{ \textnormal{6 points de 2-torsion impairs}\} \\
A \cap \mathbb{P}^4_+ &=& \{ \textnormal{10 points de 2-torsion pairs} \}.
\end{eqnarray*}
Dans les coordonn\'es $Y_i$ et $Z_j$ les quadriques $Q_b$ prennent la forme ~\cite{Barth}
\begin{flushleft}
$ \left(
\begin{array}{c} Q_1 \\  Q_2 \\ Q_3 \\ Q_4 \\ Q_5  \end{array} \right)= 
\left( \left( \begin{array}{ccccc} 
Y_0^2  & Y_1^2  & Y_2^2  & Y_3^2  & Y_4^2  \\ 
Y_1^2  & Y_0Y_1  & Y_3Y_4  & Y_2Y_4  & Y_2Y_3 \\
Y_2^2  & Y_3Y_4  & Y_0Y_2  & Y_1Y_4  & Y_1Y_3 \\
Y_3^2  & Y_2Y_4  & Y_1Y_4  & Y_0Y_3  & Y_1Y_2 \\
Y_4^2  & Y_2Y_3  & Y_1Y_3  & Y_1Y_2  & Y_0Y_4 \\
\end{array} \right) \right. +  $
\end{flushleft} 
\begin{flushright}
$
+ \left. \left(
\begin{array}{rrrrr} 
 0 & -Z_1^2 & -Z_2^2 & -Z_3^2 & -Z_4^2 \\
Z_1^2 & 0 & -Z_3Z_4 & -Z_2Z_4 & -Z_2Z_3 \\
Z_2^2 & Z_3Z_4 &  0 & Z_1Z_4 & -Z_1Z_3  \\
Z_3^2 & Z_2Z_4 & -Z_1Z_4 &  0 & Z_1Z_2  \\
Z_4^2 & Z_2Z_3 & Z_1Z_3 & - Z_1Z_2 & 0 
\end{array} \right) \right ) \cdot \beta
$
\end{flushright}
\begin{flushleft}
$ \left(
\begin{array}{c} Q_6 \\  Q_7 \\ Q_8 \\ Q_9  \end{array} \right)=$
\end{flushleft}
$$
\left(
\begin{array}{rrrrr} 
 -\beta_1Z_1 & 2\beta_0Z_1 & \beta_3Z_4-\beta_4Z_3 &  \beta_4Z_2-\beta_2Z_4 & \beta_2Z_3-\beta_3Z_2   \\
 -\beta_2Z_2 & -\beta_3Z_4-\beta_4Z_3 &  2\beta_0Z_2 & \beta_1Z_4+\beta_4Z_1 & \beta_1Z_3-\beta_3Z_1   \\
 -\beta_3Z_3 & -\beta_2Z_4-\beta_4Z_2 &  \beta_1Z_4-\beta_4Z_1 & 2\beta_0Z_3 & \beta_1Z_2+\beta_2Z_1  \\
 -\beta_4Z_4 & -\beta_2Z_3-\beta_3Z_2 &  \beta_1Z_3+\beta_3Z_1 &\beta_1Z_2-\beta_2Z_1 & 2\beta_0Z_4   \\
\end{array} \right) \cdot 
\left( \begin{array}{c} Y_0 \\  Y_1 \\ Y_2 \\ Y_3 \\ Y_4  \end{array} \right).
$$
Les restrictions des quadriques $Q_i$ \`a $\mathbb{P}^3_-$ est donn\'ees par 
$$ \left(
\begin{array}{c} q_1 \\  q_2 \\ q_3 \\ q_4 \\ q_5  \end{array} \right)= 
\left( 
\begin{array}{rrrrr} 
 0 & -Z_1^2 & -Z_2^2 & -Z_3^2 & -Z_4^2 \\
Z_1^2 & 0 & -Z_3Z_4 & -Z_2Z_4 & -Z_2Z_3 \\
Z_2^2 & Z_3Z_4 &  0 & Z_1Z_4 & -Z_1Z_3  \\
Z_3^2 & Z_2Z_4 & -Z_1Z_4 &  0 & Z_1Z_2  \\
Z_4^2 & Z_2Z_3 & Z_1Z_3 & - Z_1Z_2 & 0 
\end{array}  \right) \cdot \left(
\begin{array}{c} \beta_0 \\   \beta_1 \\  \beta_2 \\  \beta_3 \\  \beta_4  \end{array} \right) = (q_{ij}(Z))\cdot \beta .
$$
Les quadriques $q_i$ d\'efinent l'application de \emph{Steiner} $\beta :\mathbb{P}^3_- \rightarrow \mathbb{P}^4 $ sur 
les points $z \in \mathbb{P}^3_-$ tels que la matrice $q_{ij}$ est de rang 4. L'application $\beta$ envoie $z$ sur le noyau de 
la matrice $q_{ij}$. Les 6 points de 2-torsion de l'intersection $A \cap\mathbb{P}^3_- $ ont tous
la m\^eme image par $\beta$, autrement dit  $\beta :\mathbb{P}^3_- \rightarrow \mathbb{P}^4 $ est de degr\'e 6 sur
son image. Coble a montr\'e (~\cite{H} pag. 190) que l'image de cette application est la \emph{ quartique de Burkhardt}. 
La quartique de Burkhardt est la seule quartique dans $\mathbb{P}^4$ invariante par l'action du groupe unitaire de reflexions 
$\mathbb{P}Sp(4, \mathbb{Z}/3)$. \\
Le th\'eor\`eme suivant (5.3.4 ~\cite{H}) r\'esume cette situation
\begin{The} 
Il existe une application $\mathbb{P}Sp(4, \mathbb{Z}/3)$-\'equivariante et bir\'eguliere d'un ouvert de Zariski de la 
quartique de Burkhardt dans un ouvert de Zariski de l'espace de modules ${\textit{\LARGE a}}_{(3)}$    
\end{The}
Pour plus de d\'etails sur les propri\'et\'es de la quartique de Burkhardt voir ~\cite{H}.

\chapter{Preuve de la conjecture de Dolgachev}

\section{L'application $\theta$}

Dans tout ce travail on consid\`ere une courbe projective $C$ complexe non singuli\`ere de genre 2. Soit $JC$ la
jacobienne de $C$ et $J^1$ la vari\'et\'e (isomorphe \`a $JC$) qui param\`etre les fibr\'es en droites de
degr\'e 1 sur $C$. Dans la suite on notera $h^0(\xi)$ au lieu de
$h^0(C,\xi)$ pour tout fibr\'e $\xi$ sur $C$. Soit $\Theta$ le diviseur th\^eta canonique dans $J^1$ d\'efini
ensemblistement par
\begin{eqnarray*}
\Theta = \{\xi \in J^1 \mid h^0(\xi) \geq 1 \}.
\end{eqnarray*}
On note
$\su$ l'espace des modules des fibr\'es vectoriels semi-stables de rang 3 avec d\'eterminant trival.
C'est une vari\'et\'e projective, dont les points peuvent \^etre identifi\'es avec les classes d'isomorphisme de fibr\'es vectoriels
qui sont des sommes directes de fibr\'es vectoriels stables de degr\'e 0.
\\
On consid\`ere l'application
\begin{eqnarray*}
\theta :\su \longrightarrow |3\Theta|\simeq \p
\end{eqnarray*}
qui associe \`a un fibr\'e $E$ le diviseur de support 
\begin{eqnarray*}
\Theta_E :=\{\xi\in J^1 \mid H^0(C, \xi\otimes E) \neq 0\}.
\end{eqnarray*}
On v\'erifie que $\Theta_E \in |3\Theta|$ (c.f. \S 2 ~\cite{B1} ) et par le corollaire 1.7.4. ~\cite{R}, l'application $\theta$ est bien d\'efinie.
Cette application a \'et\'e \'etudi\'ee pour des courbes de genre $g\geq 2$ et de rang $r=2$ dans ~\cite{B1} et ~\cite{N-R2} et pour
un rang quelconque dans ~\cite{B-N-R} et  ~\cite{B2}.
\\

Pour tout fibr\'e en droites $L\in J^1$  on d\'efinit
\begin{displaymath}
\Theta_L = \{ E \in \mathcal{M} \mid h^0(E\otimes L)\geq 1 \}.
\end{displaymath}
Il s'av\`ere que $\Theta_L$ est un diviseur de Cartier dans $\mathcal{M}$ ~\cite{D-N}. Le fibr\'e en droites associ\'e $\mathcal{L}:=
\mathcal{O}_{\mathcal{M}}(\Theta_L)$ ne d\'epend pas du choix de $L$; il est appel\'e le \emph{fibr\'e d\'eterminant}. Le partie \emph{i)} du
th\'eor\`eme suivant  montre qu'il est canonique.
\begin{The} \label{picard} \textnormal{~\cite{D-N}} \\
i) Le groupe de Picard $Pic(\mathcal{M})$ est isomorphe \`a $\mathbb{Z}$, engendr\'e par $\mathcal{O}(\Theta_L)$.\\
ii) Le faisceau canonique de $\mathcal{M}$ est isomorphe \`a $\mathcal{O}(-6\Theta_L)$.
\end{The}
\begin{Rmq}
Le diviseur $\Theta_L$ est connu dans la litt\'erature comme le \emph{diviseur th\^eta g\'en\'eralis\'e} puisque dans le cas $n=1$ la d\'efinition
correspond au diviseur th\^eta dans $J^1$. De m\^eme, par analogie avec les fibr\'es de rang 1, les sections du fibr\'e $\mathcal{L}^k$ sont
souvent appel\'ees \emph{fonctions th\^eta g\'en\'eralis\'ees.}
\end{Rmq}
\begin{The} \label{iso-can}
Il existe un isomorphisme canonique (\`a un scalaire pr\`es)
\begin{displaymath}
H^0(\mathcal{M},\mathcal{L}) \stackrel{\sim}{\longrightarrow} H^0(J^1, \mathcal{O}_{J^1}(3\Theta) )^*,
\end{displaymath}
qui rend commutatif le diagramme suivant
$$
\shorthandoff{;:!?}
\xymatrix {
        & |\mathcal{L}|^* \ar[dd]^{\wr} \\
        \mathcal{M} \ar[rd]_{\theta}  \ar[ru]^{\varphi_{\mathcal{L}}}  & \\
        & |3\Theta|
     }
$$
\end{The}
Le th\'eor\`eme ~\ref{iso-can} est d\'emontr\'e en g\'en\'eral dans~\cite{B-N-R}  pour $g \geq 2$ et pour un rang quelconque.  Il en d\'ecoule que
$\mathcal{L}=\theta^*\mathcal{O}(1)$, o\`u $\mathcal{O}(1)$ est le fibr\'e tautologique dans $|3\Theta|\simeq \mathbb{P}^8$.
\begin{Lem}   \label{degre}
L'application $\theta$ est de degr\'e 2.
\end{Lem}
\begin{Dem}
Le degr\'e de $\theta $ est \'egal au nombre d'intersection $c_1(\mathcal{L})^8$. Pour calculer
ce nombre d'intersection on utilisera la formule de Verlinde qui donne la dimension des espaces vectoriels
$H^0(\su,\mathcal{L}^k)$.
\begin{Pro} \textnormal{~\cite{V}} \textnormal{~\cite{Z}}
Formule de Verlinde pour $n=3$, $g=2$.
\begin{displaymath}
\dim H^0(\su, \mathcal{L}^k)=3 \left( \frac{k+3}{8} \right) ^2 V_{1,1,1}(k+3)
\end{displaymath}
o\`u
\begin{displaymath}
V_{l,m,n}(k)= \sum_{\substack {a,b\geq 0\\ a+b\leq k}} \left( \sin\frac{\pi a}{k} \right) ^{-2l}
\left( \sin\frac{\pi b}{k} \right) ^{-2m}\left( \sin\frac{\pi (a+b)}{k}\right) ^{-2n}
\end{displaymath}
pour $k,l,m,n \in \mathbb{N}$.
\end{Pro}
Le terme de droite de cette formule est un polyn\^ome en $k$ et son
coefficient dominant est $\frac{c_1(\mathcal{L})^8}{8!}$. Dans son article (c.f. \S 2~\cite{Z}) Zagier donne une formule, en termes de nombres
de Bernoulli, pour le coefficient dominant de $V_{h,h,h}$, avec $h=g-1$:
\begin{displaymath}
v_{h,h,h} =(-1)^h 2^{6h} \sum_{r=0}^h \binom{4h-2r-1}{2h-1}\frac{B_{2r}}{(2r)!}\frac{B_{6h-2r}}{(6h-2r)!}.
\end{displaymath}
Pour $h=1$ on obtient $v_{1,1,1}= \frac{2^4}{3(7)!}$ et donc $c_1(\mathcal{L})=2$.
\end{Dem}
Le lemme ~\ref{degre} a \'et\'e prouv\'e d'une autre fa\c con par Y. Laszlo dans  ~\cite{L2} (Section V.).
Notons $\mathcal{B}\subset \su$ le lieu de ramification de l'application $\theta$.
\begin{Lem}
L'hypersurface $\mathcal{B}$ est de degr\'e 6.
\end{Lem}
\begin{Dem}
Comme $\theta$ est un rev\^etement double sur $\mathbb{P}^8$ il existe un fibr\'e en droites $\mathcal{O}
(\delta)$ sur $\mathbb{P}^8$ tel que
\begin{displaymath}
\theta_*\mathcal{O}_{\mathcal{M}} = \mathcal{O}_{\mathbb{P}^8} \oplus \mathcal{O}_{\mathbb{P}^8}(-\delta),
\end{displaymath}
et donc $\mathcal{B}\in |\mathcal{O}_{\mathbb{P}^8}(2\delta)|$. On sait que dans cette situation (Lemme 17.1~\cite{B-P-V}) on a la relation
suivante
\begin{displaymath}
K_{\mathcal{M}} = \theta^*(K_{\mathbb{P}^8} + \delta),
\end{displaymath}
o\`u $K_{\mathcal{M}}$ et $K_{\mathbb{P}^8}$ d\'esignent les diviseurs canoniques de $\mathcal{M}$ et ${\mathbb{P}^8}$ .  Par le th\'eor\`eme
~\ref{picard} \emph{i)} on a
\begin{displaymath}
\theta^*(\mathcal{O}_{\mathbb{P}^8}(K_{\mathbb{P}^8} + \delta)) \simeq \theta^*(\mathcal{O}_{\mathbb{P}^8}(\delta-9))
\simeq \mathcal{O}_{\mathcal{M}}(-6\Theta_L),
\end{displaymath}
d'o\`u on en d\'eduit $\delta=3$ et donc que $\mathcal{B}$ est de degr\'e 6.
\end{Dem}

\section{Le lieu de ramification de l'application $\theta$}

Le but de cette section est de d\'ecrire le lieu de ramification $\mathcal{B}$ en termes d'une certaine involution sur $\su$. Soit $\iota : C \rightarrow C$
l'involution hyperelliptique. On note $\mathfrak{i}$ l'involution sur $\su$ d\'efinie par
\begin{displaymath}
\mathfrak{i}: E \mapsto \iota^*E^*,
\end{displaymath}
o\`u $E$ est un repr\'esentant d'une classe d'isomorphisme de fibr\'es semi-stables et $E^*$ son fibr\'e dual.
\begin{Pro}   \label{theta}
L'involution $\mathfrak{i}$ commute avec l'application $\theta$.
\end{Pro}
\begin{Dem}
Il suffit de d\'emontrer, au niveau des supports de diviseurs, que $\Theta_{\iota^*E^*}=\Theta_{E}$. D'abord,
observons que
\begin{eqnarray*}
\iota^*\Theta_E & = & \{ \iota^*\xi \in J^1 \mid  h^0(E \otimes \xi) \neq 0 \} \\
& = & \{ \xi\in J^1 \mid h^0(E \otimes \iota_{*}\xi) \neq 0 \} \\
& = & \{ \xi\in J^1 \mid h^0(\iota^*E \otimes \xi) \neq 0 \}\\
& = & \Theta_{\iota^*E} .
\end{eqnarray*}
Rappelons que l'involution $\iota^*$ restreinte \`a $J^1$ est donn\'ee par $\xi \mapsto \xi^{-1}\otimes
\omega_C$, o\`u $\omega_C $ est le fibr\'e canonique sur $C$. Comme $\chi (E^* \otimes\xi)=0$ on a
\begin{eqnarray*}
\Theta_{E^*} & = & \{ \xi \in J^1 \mid  h^0(E^* \otimes \xi) \neq 0 \} \\
& = & \{ \xi\in J^1 \mid h^0(E \otimes \xi^{-1} \otimes \omega_C) \neq 0 \} \\
& = & \{ \xi\in J^1 \mid h^0(E \otimes \iota^*\xi) \neq 0 \}\\
& = & \iota^*\Theta_E .
\end{eqnarray*}
Ainsi on obtient que
\begin{displaymath}
\Theta_{\iota^*E^*}= \iota^*\Theta_{E^*}= (\iota^*\circ\iota^*)\Theta_E =\Theta_E.
\end{displaymath}
\end{Dem}
Pour conclure que $\mathfrak{i}$ est l'involution associ\'ee au rev\^etement double $\theta$, il faut montrer que cette involution n'est pas triviale.
Pour cela on va exhiber un fibr\'e $E \in \su$ tel que $\mathfrak{i} (E) \ncong E$. \\
On consid\`ere un fibr\'e vectoriel $F$ sur $C$ de rang 2 et de degr\'e 1. On note $W=\underline{End} F$ le fibr\'e
d'endomorphismes de $F$ et $W_0 = \underline{End}_0 F$ le noyau de l'application trace
\begin{displaymath}
Tr:W \rightarrow \fs,
\end{displaymath}
i.e. $W_0$ est le fibr\'e des endomorphismes de $F$ de trace nulle. Donc $\Rg W_0=3$ et $\Det W_0 \simeq \fs$.
\\
On cherche un fibr\'e $F$ tel que $W_0$ soit stable parce que dans ce cas la S-\'equivalence
dans $\su$ co\" incide avec l'isomorphisme de fibr\'es. Soit $\suu$ l'espace des modules de classes d'\'equivalence de fibr\'es stables de
rang 2, avec d\'eterminant isomorphe \`a $\lambda$, suppos\'e de degr\'e 1.
\begin{Lem}
Soit $F \in \suu$, pour que $\underline{End}_0(F)$ soit stable il faut et il suffit que $F$ ne soit pas isomorphe \`a
$ F\otimes\alpha$ pour tout $\alpha \in JC[2]$ non nul.
\end{Lem}
\begin{Dem}
Soit $F \in \suu$. Supposons que $\underline{End}_0(F)$ n'est pas stable. Donc, il contient
un sous-fibr\'e $L$ tel que $\Deg L \geq 0$ et $\Rg L < 3$ . Supposons $\Rg L= 2$. Observons que
\begin{displaymath}
\underline{End}(F)^* \simeq \underline{End}(F^*) \simeq \underline{End}(F),
\end{displaymath}
d'o\`u il d\'ecoule que $\underline{End}_0(F)^* \simeq \underline{End}_0(F)$.  En prenant le dual de la suite exacte :
\begin{displaymath}
0 \longrightarrow L \longrightarrow \underline{End}_0(F) \longrightarrow Q \longrightarrow 0,
\end{displaymath}
o\`u $Q$ est le fibr\'e quotient de rang 1, on obtient la suite exacte
\begin{displaymath}
0 \longrightarrow Q^* \longrightarrow \underline{End}_0(F) \longrightarrow L^* \longrightarrow 0 .
\end{displaymath}
Comme $\Deg Q^* =-\Deg Q =\Deg L$ on peut donc supposer que L est un fibr\'e en droites.
Le morphisme injectif non nul $L \hookrightarrow \underline{End}_0(F)$ induit un homomorphisme
\begin{displaymath}
u: F \longrightarrow F\otimes L^{-1},
\end{displaymath}
non nul de trace nulle. Comme $F$ est stable, le morphisme $u$
est n\'ecessairement un isomorphisme; alors $F\simeq F\otimes L^{-1}$ et on a
$$
\lambda\simeq \Lambda^2 F \simeq \Lambda^2 (F \otimes L^{-1}) \simeq \lambda \otimes L^{-2} ,
$$
donc $L^2 \simeq \fs$, i.e.,  $L\in JC[2]$. L'hypoth\`ese sur $F$ entra\^{\i}ne $L \simeq \fs$, donc $u$ est un endomorphisme de $F$. Par ailleurs,
les endomorphismes d'un fibr\'e stable sont des homoth\'eties; or, comme $u$ est de trace nulle, il est n\'ecessairement
l'homomorphisme nul, ce qui nous donne une contradiction. \\
Supposons qu'il existe $\alpha \in JC[2]$ non nul tel que $u: F\otimes \alpha \stackrel{\sim}{\rightarrow} F$. Donc $u$
induit un morphisme injectif $\alpha \rightarrow \underline{End}(F) = \fs \oplus \underline{End}_0(F)$, dont la projection sur $\fs$ est nulle car $\alpha$ est
non trivial. Donc $\alpha$ est un sous-fibr\'e de $\underline{End}_0(F)$, ce qui contredit la stabilit\'e de $\underline{End}_0(F)$.
\end{Dem}
\'Etant donn\'e un groupe alg\'ebrique  $G$ on note $\mathcal{O}G$ le faisceau des fonctions holomorphes sur $C$ \`a valeurs dans $G$.
Alors les classes d'isomorphisme de $G$-fibr\'es sur $C$ sont param\'etris\'ees par le groupe de cohomologie de C\v ech $H^1(C,
\mathcal{O}G)$. Dans notre situation on a la suite exacte (2.7 ~\cite{O} )
\begin{equation} \label{suite}
0 \longrightarrow \textnormal{Pic}(C) \longrightarrow  H^1(C, \mathcal{O}GL_2)  \stackrel{\pi}{\longrightarrow}
H^1(C, \mathcal{O}SO_3) \rightarrow 0,
\end{equation}
o\`u l'application $\pi$ est donn\'ee par $F \mapsto  \underline{End}_0(F)$.
Prenons $F \in \suu$. On a donc
\begin{displaymath}
\iota^*\underline{End}_0(F) \simeq  \underline{End}_0(\iota^*F).
\end{displaymath}
Supposons  $\underline{End}_0(\iota^*F) \simeq  \underline{End}_0(F)$. Donc, par la suite (\ref{suite}), il existe $\alpha \in \mathrm{Pic}(C)$ tel
que
\begin{displaymath}
\iota^*F \simeq F\otimes \alpha.
\end{displaymath}
On a donc
$$
\iota^* \lambda \simeq \Lambda^2 \iota^*F  \simeq  \Lambda^2 (F\otimes \alpha) \simeq \lambda \otimes \alpha^2.
$$
En choissisant $\lambda$ tel que $\iota^*\lambda\simeq\lambda$ on en d\'eduit que $\alpha \in JC[2]$.
Il suffit donc, pour que l'involution $\mathfrak{i}$ soit non triviale, de d\'emontrer qu'il existe un fibr\'e $F \in
\suu$ tel que $\iota^*F\not\simeq F\otimes\alpha$, pour tout $\alpha \in JC[2]$. Pour montrer l'existence d'un tel fibr\'e dans $\suu$ on
va utiliser la description du groupe de automorphismes de $\suu$ donn\'ee par Newstead dans son article ~\cite{New}.

\section{Les automorphismes de $\suu$}

Newstead utilise une description g\'eom\'etrique de l'espace de modules $\suu$ pour expliciter ses automorphismes. On consid\`ere le
groupe
\begin{displaymath}
G:= \{ (\rho,M) \mid \rho \in Aut(C),  M \in \textrm{Pic}(C) \textrm{ \ avec \ } M^2 \simeq
\lambda\otimes(\rho^*\lambda)^* \}
\end{displaymath}
On munit  $G$ d'une structure de groupe donn\'ee par
\begin{displaymath}
(\rho_1, M_1)\cdot(\rho_2, M_2)= (\rho_1 \circ\rho_2, M_2 \otimes\rho_2^*(M_1)).
\end{displaymath}
Le groupe $G$ agit sur $\suu$ \`a droite par
\begin{displaymath}
F \cdot(\rho, M) = (\rho^*(F)\otimes M),
\end{displaymath}
pour tout $F \in \suu$. On observe que
\begin{displaymath}
\Det (\rho^*(F)\otimes M) \simeq \rho^*(\lambda)\otimes M^2 \simeq \lambda.
\end{displaymath}
On peut donc d\'efinir un homomorphisme naturel $\Phi: G \rightarrow Aut( \suu )$ par
\begin{displaymath}
(\rho, M) \mapsto \{ F \mapsto \rho^*(F)\otimes M \}.
\end{displaymath}
\begin{The} \label{Newstead}
\textnormal{~\cite{New}} L'application $\Phi$ est un isomorphisme de groupes.
\end{The}
Par le th\'eor\`eme ~\ref{Newstead} l'automorphisme de $\suu$ donn\'e par $F \mapsto \iota^* F \otimes \alpha$ n'est pas l'identit\'e.
Par cons\'equent, il existe un fibr\'e
$F \in \suu$ tel que $\iota^*F\not\simeq F\otimes\alpha$ pour tout $\alpha \in JC[2]$.\\
On a donc d\'emontr\'e la proposition suivante
\begin{Pro}
L'involution associ\'ee au rev\^etement  $\theta$ est  $\mathfrak{i}: E \mapsto \iota^*E^*$.
\end{Pro}
\begin{Cor}
L'hypersurface sextique de ramification $\mathcal{B}\subset \su$ est d\'efinie par
\begin{displaymath}
\mathcal{B}= \{ E \in \su \mid E^* \sim_{s} \iota^*E \} .
\end{displaymath}
\end{Cor}

\section{Calcul du degr\'e de la vari\'et\'e duale $\mathcal{C}^*$} \label{degre_duale}

On cherche \`a calculer le degr\'e de  l'hypersurface $\mathcal{C^*}$, l'hypersurface duale de la cubique de Coble.
Soit $\mathcal{D}: |3\Theta|^*  \dashrightarrow |3\Theta| \cong \mathbb{P}^{8*}$ l'application rationnelle d\'efinie par les d\'eriv\'ees partielles
de $\mathcal{C}$,  i.e.,
\begin{displaymath}
\mathcal{D}: x \mapsto \left[ \frac{\partial F}{\partial X_0}(x), \ldots, \frac{\partial F}{\partial X_8}(x) \right],
\end{displaymath}
o\`u $F=0$ est l'\'equation qui d\'efinit la cubique de Coble. 
La vari\'et\'e duale de $\mathcal{C}$ est par d\'efinition l'adh\'erence de l'image de la partie lisse de $\mathcal{C}$ 
par l'application
$\mathcal{D}$.\\
Soit $(A,L)$ une surface ab\'elienne complexe principalement polaris\'ee. Il est connu que l'id\'eal de la surface $A$
plong\'ee dans
$\mathbb{P}^8$ par le fibr\'e tr\`es ample $L^{\otimes 3}$ est engendr\'e par des quadriques et des cubiques.
Barth a d\'emontr\'e ~\cite{Barth} que le faisceau d'id\'eaux $\mathcal{I}_{A/\mathbb{P}^8}$ est engendr\'e par les
quadriques qui s'annulent en $A$ lorsque $A$ n'est pas le produit de courbes elliptiques.
\\
Dans notre cas,$J^1$ ($\simeq JC$) est d\'efini sch\'ematiquement par les d\'eriv\'ees partielles de $\mathcal{C}$, i.e.
$J^1$ est le lieu singulier de la cubique de Coble o\`u l'application $\mathcal{D}$
n'est pas d\'efinie. On consid\`ere le diagramme commutatif suivant

$$
\shorthandoff{;:!?}
\xymatrix {
        E  \ar@{^{(}->}[r]^{j}  \ar[d] & \widetilde{\mathbb{P}}^8 \ar[d]^{\pi} \ar[rd]^f  \\
        J^1 \ar@{^{(}->}[r] ^{i} & \mathbb{P}^8  \ar@{-->}[r]^{\mathcal{D}}  & \mathbb{P}^{8*}
     }
$$
o\`u $\pi: \widetilde{\mathbb{P}}^8 \rightarrow \mathbb{P}^8$ est l'\'eclatement de $\mathbb{P}^8$  le long de $J^1$ et $E$ est le diviseur
exceptionnel. Ici on a pris le
plongement $i$ comme celui d\'efini par le fibr\'e en droites $\mathcal{O}(3\Theta)$ sur $J^1$.  Soit $h=i^*\mathcal{O}(1)$, on a donc
$h^2=9 (\Theta.\Theta)=18$. \\
Soit $e=[E]$ la classe du diviseur exceptionnel dans  $ \widetilde{\mathbb{P}}^8$ et $H$ la classe d'un hyperplan dans $\mathbb{P}^8$; on
note $\widetilde{H}:=\pi^*H$.  Par la proposition 4.4 ~\cite{F} on a
\begin{displaymath}
d\cdot\Deg (\mathcal{C}^*) = \int_{\pi^{-1}(\mathcal{C})}  c_1(f^*\mathcal{O}(1))^7,
\end{displaymath}
o\`u $d$ est le degr\'e g\'en\'erique de l'application $\mathcal{D}$. 
Soit $N= N_{J^1}\mathbb{P}^8$ le fibr\'e normal. \'Etant donn\'e que le fibr\'e tangent \`a $J^1$ est trivial, les classes de Chern totales de $N$ et de $i^*T_{\mathbb{P}^8}$ sont \'egales. On a donc
\begin{displaymath}
c(N)= c(i^*T_{\mathbb{P}^8}) = (1 + h)^9 = 1 +9h + 36h^2.
\end{displaymath}
On note $\xi := c_1(\mathcal{O}_E (-1))=e_{|_E}$ et  $\tilde{h}:= \pi^*(h)$. Par le lemme 3.2.~\cite{F} on a
\begin{displaymath}
c_6(\pi^*N \otimes \mathcal{O}_E (1))= \sum_{i=0}^6 c_1 (\mathcal{O}_E(1))^i c_{6-i}(\pi^* N) =0,
\end{displaymath}
d'o\`u
\begin{equation} \label{relation}
\xi^6 - 9\tilde{h}\xi^5+36\tilde{h}^2\xi^4=0.
\end{equation}
\`A l'aide de cette relation on va calculer les produits $\widetilde{H}^r \cdot e^{8-r}$ pour $0\leq r \leq 8$, dont on aura besoin pour le
calcul du degr\'e. \\
On observe d'abord que $\tilde{h}^2\xi^5 = -18$ car $h^2=18$. Donc, en multipliant (\ref{relation}) par $\tilde{h}$, on obtient
\begin{displaymath}
\tilde{h}\xi^6-9\tilde{h}^2\xi^5=0,
\end{displaymath}
et donc $\xi^6\tilde{h}=-162$. En multipliant (1) par $\xi$ on a
\begin{displaymath}
\xi^7-9\tilde{h}\xi^6 + 36\tilde{h}^2\xi^5=0,
\end{displaymath}
d'o\`u $\xi^7= 36(18)-9(162)=-810$.  On a
\begin{displaymath}
\widetilde{H}^r \cdot  e^{8-r} = (\tilde{h}^r \cdot \xi^{8-r-1})_E.
\end{displaymath}
Observons que la transform\'ee stricte de $\mathcal{C}$ dans $\widetilde{\mathbb{P}}^8$ correspond au diviseur 
$3\widetilde{H} - 2e$ puisque le lieu singulier de $\mathcal{C}$ est de multiplicit\'e 2 dans $\mathcal{C}$. 
En utilisant le fait que $\tilde{h}^r=0$ pour $r\geq3$ il en r\'esulte
\begin{eqnarray*}
\int_{\pi^{-1}(\mathcal{C})} c_1(f^*\mathcal{O}(1))^7 &=&  (3\widetilde{H}-2e)(2\widetilde{H}-e)^7 \\
&=&  -3\widetilde{H}e^7 +210\widetilde{H}^2e^6 +2e^8 -28\widetilde{H}e^7 +384\widetilde{H} ^8\\
&=&  -31 \tilde{h}\xi^6 + 210 \tilde{h}^2 \xi^5  +2 \xi^7 +384 \\
&=& 6.
\end{eqnarray*}
Ceci montre que $\Deg (\mathcal{C}^*) \leq 6 $. On v\'erifiera ult\`erieurement (c.f. section 2.7  (\ref{grado})) que
$d=1$, i.e.  $ \Deg (\mathcal{C}^*) = 6 $.

\section{L'application duale sur les s\'ecantes de $J^1$}

Consid\'erons $\D: \p \dashrightarrow \mathbb{P}^{8*}$ l'application duale, qui n'est pas d\'efinie sur le lieu singulier
de $\mathcal{C}$. Consid\'erons $J^1= \Sing \mathcal{C}$ plong\'ee dans $\p$ par le syst\`eme lin\'eaire
$|3\Theta|$. 
On note $Q_i:=\frac{\partial F}{\partial X_i}$ et $\widetilde{Q}_i$ les formes bilin\'eaires associ\'ees, avec $\ i=0,
\ldots ,8 $.
On prend $a,b \in J^1$ et on consid\`ere la s\'ecante $\ell$ de $J^1$ qui passe par $a$ et
$b$. Puisque $\mathcal{C}$ est singuli\`ere le long de $J^1$, la multiplicit\'e de l'intersection $ \ \ell \cap 
\mathcal{C} \ $ est au moins 4, ce qui
entra\^{\i}ne que la droite $\ell$ est contenue dans $\mathcal{C}$.
\\
\begin{Lem}    \label{singularite}
L'application duale $\D$ est constante le long de la droite $\ell$, i.e. tous les points de $\ell$ autres que $a,b$ ont 
le m\^eme hyperplan tangent
\`a la cubique $\mathcal{C}$.
\end{Lem}
\begin{Dem}
On \'ecrit  $p=a+tb \in l$ pour $t\in \mathbb{C}$.  On a donc  $\D(p)= \D(a+tb) =  [\widetilde{Q}_0(a+tb,a+tb), \ldots,
\widetilde{Q}_8(a+tb,a+tb)]$.  En
d\'eveloppant la forme $\widetilde{Q}_i(a+tb,a+tb)$ on obtient
\begin{eqnarray*}
\widetilde{Q}_i(a+tb,a+tb) & = &  \widetilde{Q}_i(a,a) +t^2\widetilde{Q}_i(b,b) + 2\widetilde{Q}_i(a,tb) \\
                    & =&  2t\widetilde{Q}_i(a,b)
\end{eqnarray*}
pour tout $i=0,\ldots,8 $. Donc l'application $\D$ ne d\'epend que de $a$ et $b$ sur la droite s\'ecante.
\end{Dem}
\begin{Rmq}
Observons que par le lemme pr\'ec\'edent l'application duale $\mathcal{D}$ n'est pas bijective sur les points lisses de la droite $\ell$. Or, comme
$\mathcal{C}^{**} \simeq \mathcal{C}$, l'application duale sur $\mathcal{C}^*$ n'est pas d\'efinie sur l'image de $\ell$ par $\mathcal{D}$, ce qui
implique que $\mathcal{D}$ envoie $\ell$ sur un point singulier de $\mathcal{C}^*$.
\end{Rmq}
\\

\begin{figure}[!h]
\begin{center}
\includegraphics[width=6cm, height=6cm]{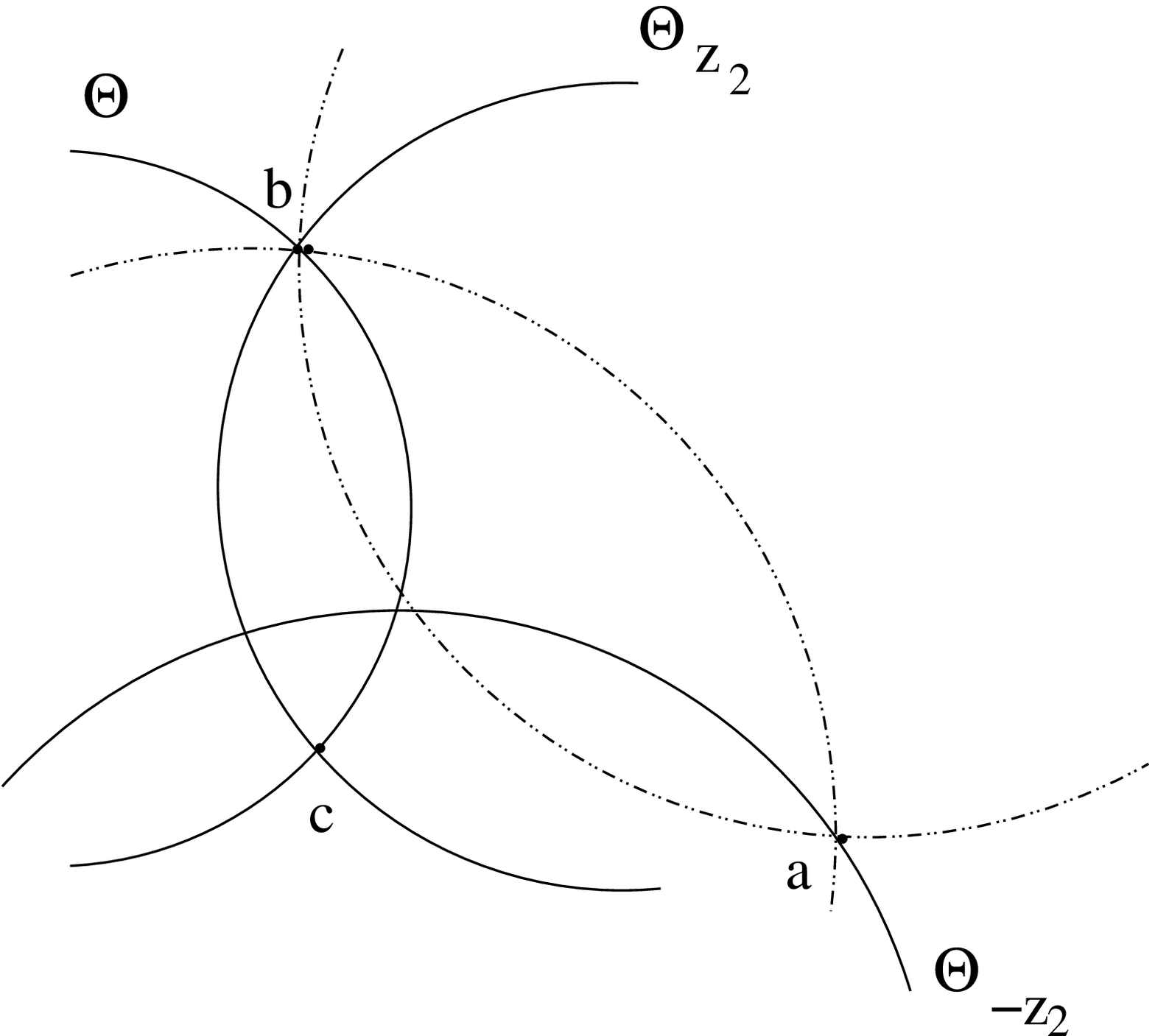}
\end{center}
\caption{}
\label{divisores}
\end{figure}

Soit  $\widetilde{F}$ la forme trilin\'eaire associ\'ee \`a la cubique $\mathcal{C}$. Donc, l'\'equation lin\'eaire
 $\widetilde{F}(a,b,\cdot)=0$ pour $a,b\in J^1, a\neq b$, d\'efinit  l'hyperplan tangent \`a la cubique en un point $p$ lisse
 sur la droite $a\cdot b$. Autrement
$\widetilde{F} (a,a,\cdot)\equiv 0$. On notera $\varphi: J^1 \times J^1 \dashrightarrow \mathbb{P}^{8*} \simeq |3\Theta|$
l'application d\'efinie par
\begin{displaymath}
        (a,b) \mapsto  \{ \widetilde{F}(a,b,\cdot)=0 \}.
\end{displaymath}
On pose $\mathcal{L}:=\mathcal{O}_{J^1}(3\Theta)$ et $V:= H^0(J^1,\mathcal{L})^*$. On a le diagramme commutatif suivant
\begin{eqnarray*}
\shorthandoff{;:!?}
\xymatrix{
	J^1\times J^1  \ar@{-->}[r]^{\varphi} \ar[rd]_{\phi} & \mathbb{P}(V^*) \\
	 &  \mathbb{P}(V\otimes V)  \ar[u]_{\gamma}
     }
\end{eqnarray*}
o\`u $\phi$ est l'application d\'efinie par le syst\`eme lin\'eaire $\mathcal{L}\boxtimes \mathcal{L}:= p_1^*(\mathcal{L})
\otimes p_2^*(\mathcal{L})$ et $\gamma$ est l'application lin\'eaire donn\'ee
par la forme trilin\'eaire $\widetilde{F}$. Donc, $\varphi$ est d\'efinie par un sous-syst\`eme lin\'eaire de
 $\mathcal{L}\boxtimes \mathcal{L}$.
\begin{Pro}  \label{phi-non-deg}
L'image de $\varphi$ dans $|3\Theta|\simeq \mathbb{P}(V^*)$ est non-d\'eg\'en\'er\'ee.
\end{Pro}
\begin{Dem}
L'action du groupe $A[3]$ sur $J^1$ induit une action sur $J^1\times J^1$ par $\eta \cdot (a,b)=( a+\eta, b+\eta)$.
Comme $\phi$ et $\gamma$ sont $\He$-\'equivariantes, alors $\varphi$ est aussi $\He$-\'equivariante.
Soit $\mathbb{P}(W)$ l'hyperplan engendr\'e par les \'el\'ements de l'image de $\varphi$, avec $W\subset
 H^0(J^1,\mathcal{L})$. On en d\'eduit
que  $W$ est un sous-espace stable par l'action du groupe de Heisenberg. Mais on sait que
$H^0(J^1,\mathcal{L})$ est une r\'epresentation irr\'eductible de $\He$ alors soit  $W={0}$, soit $W=V^*$.  Cela entra\^{\i}ne que
l'image de $\varphi$ est non-d\'eg\'en\'er\'ee.
\end{Dem}
Dans la suite pour $x \in \mathrm{Pic}(C)$, on notera $\Theta_x$ le diviseur th\^eta translat\'e $t^*_x\Theta$.
\begin{Lem}  \label{translate}
Soit $C$ une courbe de genre 2 et soient $a,b \in \mathrm{Pic}^1(C)=J^1$ distincts. Il existe
exactement deux translat\'es du diviseur th\^eta  $\Theta_x$ et $\Theta_y$, avec $\ x,y \in JC$, tels que $\Theta_x\cap
\Theta_y=\{ a,b\}$.
\end{Lem}
\begin{Dem}
Calculer le nombre de diviseurs lin\'eairement  \'equivalents \`a  $\Theta$  qui passent par $a$ et $b$ \'equivaut \`a
calculer le produit d'intersection
\begin{displaymath}
        (\Theta_a \cdot \Theta_b)=(\Theta \cdot \Theta )= 2g-2 =2.
\end{displaymath}
Soit $\alpha: C \rightarrow \mathrm{Pic}^1(C)$ le plongement canonique. Puisque l'application $C^{(2)} \rightarrow JC$ d\'efinie par
\begin{displaymath}
       (p,q)\mapsto \mathcal{O}_C(p-\iota q ) = \alpha(p)-\alpha(\iota q)
\end{displaymath}
est surjective il existe $p,q \in C$ tels que $a-b =\mathcal{O}_C(p-\iota q) $. En posant $x=\alpha(p)-a$ et
$y=\alpha(\iota p)-b$, on  v\'erifie que   $\Theta_x \cap \Theta_y=\{ a,b\}$.
\end{Dem}
\begin{Rmq}
La propri\'et\'e \'enonc\'ee dans le lemme ~\ref{translate} est sym\'etrique de celle-ci: pour tout $x,y\in JC$, $x\neq y$
il existe exactement deux translat\'es $\Theta_a$ et $\Theta_b$ avec $a,b\in \mathrm{Pic}^1(C)$ tels que  $\Theta_a \cap
\Theta_b=\{ x,y \}$.
\end{Rmq}
\\
\\
Gr\^ace au lemme \ref{translate} on peut d\'efinir une application rationnelle $\psi: J^1\times J^1 \dashrightarrow |3\Theta|$ par
\begin{displaymath}
(a,b) \mapsto \Theta_x + \Theta_y   +\Theta_{-x-y}, 
\end{displaymath}
o\`u   $\Theta_x \cap \Theta_y=\{ a,b \}$.
Comme pour $\varphi$ cette application n'est pas d\'efinie sur 
la diagonale $\Delta \subset J^1 \times J^1$. Consid\'erons les applications $\widetilde{\psi}$ et $\beta$ d\'efinies sur
$J^1 \times C^{(2)}$ par
\begin{displaymath}
\widetilde{\psi}:(a, p+q)  \mapsto  \Theta_{\alpha(p)-a} + \Theta_{\alpha(q)-a} + \Theta_{2a -\alpha(p)-\alpha(q)}
 \in |3\Theta| ,
\end{displaymath}
et
\begin{displaymath}
\beta:  (a, p+q)  \mapsto  (a, a + \omega_C -\alpha(p)-\alpha(q)) \in  J^1 \times J^1,
\end{displaymath}
o\`u $\omega_C$ d\'esigne le fibr\'e canonique.
\begin{Lem} \label{morphisme}
Le diagramme 
\begin{eqnarray*}
\shorthandoff{;:!?}
\xymatrix{
	J^1\times C^{(2)}  \ar[rd]^{\widetilde{\psi}} \ar[d]_{\beta} &  \\
	 J^1 \times J^1   \ar@{-->}[r]^{\psi} & |3\Theta|
     }
\end{eqnarray*}
commute. De plus, l'application $\beta$ est l'\'eclatement de $J^1 \times J^1$ le long de la diagonale $\Delta$.
\end{Lem}
\begin{Dem}
Soit $(a, p+q) \in J^1 \times C^{(2)}$. \`A l'aide du lemme \ref{translate} on a
\begin{eqnarray*}
\psi\circ \beta (a, p+q) &=& \psi (a, a+ \omega_{C} -\alpha(p)-\alpha(q))\\
&=& \Theta_{\alpha(p)-a} + \Theta_{\alpha(q)-a} + \Theta_{2a-\alpha(p)-\alpha(q)}\\
&=& \widetilde{\psi}(a, p+q),
\end{eqnarray*}
pour $p\neq \iota q$.
Le morphisme $\pi :C^{(2)} \rightarrow JC$ donn\'e par $\pi: p+q \mapsto \omega_C-\alpha(p)-\alpha(q)$ est l'\'eclatement
de $JC$ en $0$.On a donc le diagramme commutatif suivant
\begin{eqnarray*}
\shorthandoff{;:!?}
\xymatrix{
	& J^1\times C^{(2)}  \ar[ld]_{\mathrm{Id} \times \pi} \ar[d]^{\beta}   \\
	 J^1 \times JC   \ar[r]^{\gamma} & J^1 \times J^1       
    }
\end{eqnarray*} 
avec $\gamma$ l'isomorphisme donn\'e par $(a, \xi) \mapsto (a, a + \xi)$. On en d\'eduit que $\beta$ est l'\'eclatement de 
la $J^1 \times J^1$ le long de la diagonale. 
\end{Dem}
%\begin{Lem} \label{morphisme}
%L'application $\psi$ est un morphisme.
%\end{Lem}
%\begin{Dem}
%Il faut montrer que $\psi$ s'\'etend \`a la diagonale $\Delta \subset J^1 \times J^1$. Rappelons que le fibr\'e tangent 
%de $J^1$ est trivial. D'apr\`es la proposition 11.1.4 \cite{CAV}, la projectivisation de la compos\'ee
%\begin{displaymath}
%\mathcal{T}_C \stackrel{d\alpha}{\longrightarrow} \mathcal{T}_{J^1}\simeq J^1 \times \mathbb{C}^2
%\longrightarrow \mathbb{C}^2
%\end{displaymath} 
%est l'application $C \rightarrow \mathbb{P}^1$ donn\'ee par le syst\`eme lin\'eaire $|K_C|$. On fixe une section du 
%fibr\'e tangent $a \mapsto (a,t) \in J^1 \times \mathbb{C}^2$, $t\neq 0$, telle que le point correspondant \`a $t$ en
%$\mathbb{P}^1$
%n'est pas un point de Weierstrass. Il existe donc exactement deux points $p, \iota p \in C$ tels que $t$ est tangent \`a
%$\Theta$ en $\alpha (p)$ et $\alpha (\iota p)$. On pose $x=\alpha (p) -a$, $y=\alpha (\iota p)- a$ et on d\'efinit 
%l'application $\psi $ dans $\Delta$ par $\psi(a,a) = \Theta_x + \Theta_y +\Theta_{-x-y}$.
%\end{Dem}
\begin{Pro} \label{psi-non-deg}
L'image de $\psi$ dans $|3\Theta|$ est non-d\'eg\'en\'er\'ee.
\end{Pro}
\begin{Dem}
\`A l'aide du lemme \ref{translate}, on obtient que
\begin{displaymath}
       \psi(a+\eta, b+\eta) = \Theta_{x+\eta} + \Theta_{y+\eta} +\Theta_{-x-y-2\eta}=t_{\eta}\cdot \psi(a,b).
\end{displaymath}
L'image de $\psi$ dans  $\mathbb{P}H^0(J^1,\mathcal{L})$ engendre donc un sous-espace stable par l'action 
du groupe de Heisenberg, donc \'egal \`a $\mathbb{P}H^0(J^1,\mathcal{L})$.
\end{Dem}
\begin{Lem} \label{H_c}
On fixe $b \in J^1$. Alors $\psi^*\mathcal{O}(1)_{|_{J^1\times \{b \} }} \in |3\Theta|$.
\end{Lem}
\begin{Dem}
Soit $c\in J^1$ fix\'e,  $c\neq b, c\neq \iota(b)$ o\`u $\iota$ est l'application induite par l'involution hyperelliptique, 
et soit $H_c$ l'hyperplan dans
$|3\Theta|$ de tous les diviseurs qui passent par $c$.  On a
\begin{displaymath}
\psi^*(H_c)_{|_{J^1\times\{ b\} }}= \{ a \in J^1 \mid c \in   \Theta_x + \Theta_y   +\Theta_{-x-y}, \quad \textrm{o\`u}
 \quad  \Theta_x \cap
\Theta_y=\{ a,b \} \}.
\end{displaymath}
On suppose donc $b,c \in J^1$. D'apr\`es le lemme \ref{translate}
il existe $z_1, z_2\in JC$ tels que  $\Theta_{z_1} \cap \Theta_{z_2}= \{ b,c \}$. On prouvera l'\'egalit\'e suivante
\begin{displaymath}
\psi^*(H_c)_{|_{J^1\times\{ b\} }}= \Theta_{z_1} + \Theta_{z_2} + \Theta_{-z_1-z_2}.
\end{displaymath}
On observe que pour tout $a\in \Theta_{z_1} \cup \Theta_{z_2}$,
les diviseurs th\^eta translat\'es qui s'intersectent en $a$ et $b$ contiennent tous $c$ (c.f. Fig.\ref{diviseurs}), donc
$\Theta_{z_1} + \Theta_{z_2} \subset
\psi^*(H_c)_{|_{J^1\times\{ b\} }}$.  On montrera que le diviseur $\Theta_{-z_1-z_2}$ est contenu dans
$\psi^*(H_c)_{|_{J^1\times\{ b\} } }$. \\
\begin{figure}[!h]
\begin{center}
\includegraphics[width=4.8cm, height=5.5cm]{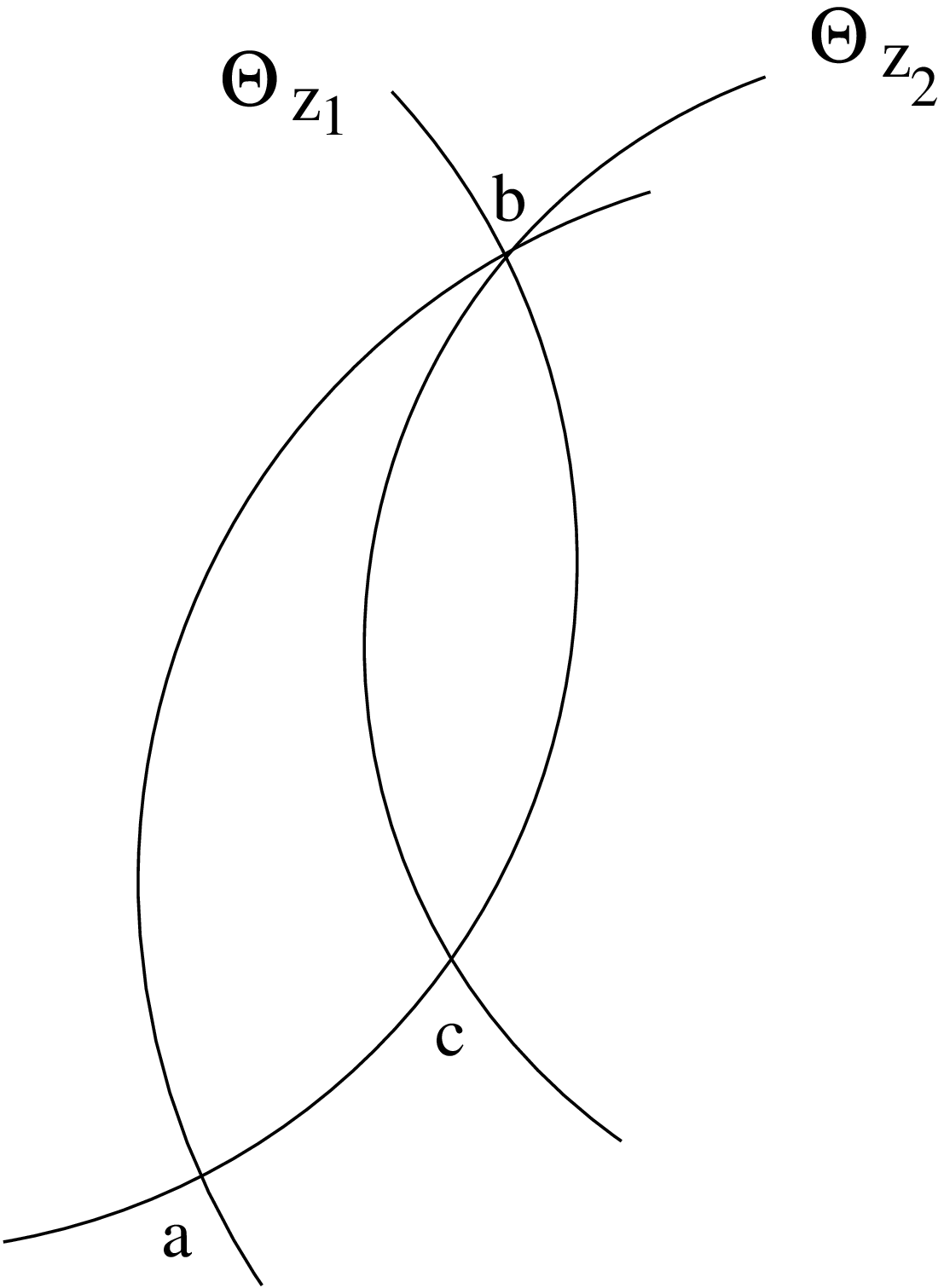}
\end{center}
\caption{}
\label{diviseurs}
\end{figure}
Apr\`es une translation on peut supposer $b,c \in \Theta\subset \mathrm{Pic}^1(C)$, i.e. , $z_1=0$ et
$z_2=\iota(c)-b=\iota(b)-c$,
tels que $\Theta \cap \Theta_{z_2}= \{ b,c \}$. En effet,
\begin{displaymath}
b \in \Theta_{z_2}  \Leftrightarrow b+z_2 =\iota(c) \in \Theta, \quad   c \in \Theta_{z_2}  \Leftrightarrow c+z_2 =\iota(b)
\in \Theta .
\end{displaymath}
Soit $a\in \Theta_{-z_2}$. On veut montrer que $a \in \psi^*(H_c)_{|_{J^1\times\{ b\} }}$. Plus pr\'ecis\'ement, 
on montrera que $c \in \Theta_{-x-y}$ o\`u $\Theta_x \cap \Theta_y = \{ a,b \}$. Il existe $p \in C$ tel que $x=\alpha(p)-a$ et 
$y=\alpha(\iota p)-b$. On a donc
\begin{eqnarray*}
a-z_2 & =& a-\iota(c) +b \\
&=& \alpha (p)-x +\alpha(\iota p) - y -\iota(c) \\
&=& \omega_{C} - \iota(c) -x-y \\
&=& c-x-y \in \Theta.
\end{eqnarray*}
On obtient ainsi $c \in \Theta_{-x-y}$ (c.f Fig.\ref{divisores}).
\\
Soit $a \in \psi^*(H_c)_{|_{J^1\times\{ b\} }}$. On veut d\'emontrer que $a \in \Theta_{z_1} + \Theta_{z_2} +
\Theta_{-z_1-z_2}$ avec $\{b, c\} =\Theta_{z_1} \cap \Theta_{z_2}$. On a 
\begin{displaymath}
c \in \Theta_x + \Theta_y +\Theta_{-x-y}, \quad \textnormal{avec} \quad \{a, b\} =\Theta_x \cap \Theta_y.
\end{displaymath}
Supposons $c \in \Theta_x$. Par le lemme \ref{translate} il existe exactement deux translat\'es du diviseur th\^eta passant
par $b$ et $c$, on peut donc supposer $x=z_1$ et on obtient $a \in  \Theta_{z_1}$. Si on suppose $c \in \Theta_y$, le
 m\^eme argument montre que $a \in  \Theta_{z_1}$. Prenons $c \in \Theta_{-x-y}$. Pour simplifier, on peut supposer 
$b,c \in \Theta$, avec $z_1=0$ et $z_2=\iota(c)-b$. On montrera que $a \in \Theta_{-z_2}$. Il existe $p \in C$ tel que
$x=\alpha(p)-a$ et $y=\alpha(\iota p)-b$. On a donc 
\begin{eqnarray*}
c-x-y & =& c - \alpha(p)+a -\alpha(\iota p) + b \\
&=& c - \omega_C +b +a \\
&=& c - \iota(b) + a \\
&=& a - z_2 \in \Theta,
\end{eqnarray*}
et donc $a \in \Theta_{-z_2}$, ce qui d\'emontre l'\'egalit\'e requise.  
\end{Dem}
\begin{The}  \label{secantes}
Les applications rationnelles $\varphi$ et $\psi$ co\" incident.
\end{The}
\begin{Dem}
Le lemme \ref{H_c} implique que $\psi$ est d\'efinie par un sous-syst\`eme lin\'eaire de $\mathcal{L}\boxtimes 
\mathcal{L}$ et les propositions \ref{phi-non-deg} et
\ref{psi-non-deg} montrent que les  sous-syst\`emes qui d\'efinissent $\varphi$ et $\psi$ sont tous les deux de dimension 
9, la dimension de $V$. On note
$\Gamma_{\varphi}, \Gamma_{\psi} $ ces sous-syst\`emes.  On note $B_s(\Gamma_{\varphi})$ et  $B_s(\Gamma_{\psi})$ les 
points base de  $\Gamma_{\varphi}$ et $
\Gamma_{\psi}$ respectivement. On montrera que la diagonale $\Delta $ dans $J^1\times J^1$ est contenue dans ces lieux 
de base.
\\
On fixe $c\in \mathbb{P}(V)$ et soit $H_c$ l'hyperplan associ\'e dans $|3\Theta|$. Alors $\varphi^{-1}(H_c)=\{ (a,b) \mid
 \widetilde{F}(a,b, c)=0
\} $ contient $\Delta$. En effet, vu que $J^1=\Sing \mathcal{C}$, alors $\widetilde{F}(x,x,\cdot)\equiv 0$ pour tout $x\in J^1$ 
et en particulier $\widetilde{F}(x,x, c)=0 \ \forall x\in J^1$. On a donc $\Delta \subset B_s(\Gamma_{\varphi})$.
\\
Dans le lemme \ref{H_c} on a montr\'e que pour $c\in J^1$ fix\'e, 
\begin{displaymath}
\psi^{-1}(H_c)_{|_{J^1\times \{ b\}}} = \Theta_{z_1}+\Theta_{z_2}+\Theta_{-z_1-z_2},
\end{displaymath}
o\`u $\Theta_{z_1} \cap \Theta_{z_2} = \{ b,c \}$. Pour tout $b\in J^1$ on a donc $b \in \psi^{-1}
(H_c)_{|_{J^1\times \{ b\}}}$ et ainsi $\Delta\subset
\psi^{-1}(H_c)$.  Cela est suffisant parce que les hyperplans $H_c$, avec $c \in J^1$, engendrent $\mathbb{P}(V^*)$.
\\
On note $\mathcal{I}$ le faisceau d'id\'eaux de $\Delta$ en $J^1\times J^1$. On a montr\'e que $\Gamma_{\varphi},
\Gamma_{\psi}  \subset H^0(J^1\times J^1, (\mathcal{L}\boxtimes\mathcal{L}) \otimes \mathcal{I})$.
On consid\`ere l'involution $i: (a,b)\mapsto (b,a)$ sur $J^1\times J^1$. On notera  $H^0(J^1\times
J^1,(\mathcal{L}\boxtimes\mathcal{L}) \otimes \mathcal{I})_+$ la partie $i$-invariante. Comme $\varphi(a,b)=\varphi(b,a)$
et $\psi(a,b)=\psi(b,a)$, on a $\Gamma_{\varphi}, \Gamma_{\psi} \subset H^0(J^1\times J^1,(\mathcal{L}\boxtimes
\mathcal{L}) \otimes \mathcal{I})_+$.On a la suite exacte
\begin{displaymath}
0 \longrightarrow H^0(J^1\times J^1,(\mathcal{L}\boxtimes\mathcal{L})\otimes \mathcal{I})\longrightarrow  H^0(J^1\times
 J^1,\mathcal{L}\boxtimes\mathcal{L})
\stackrel{rest_{|_{\Delta}}}{\longrightarrow} H^0(\Delta,\mathcal{L}\boxtimes\mathcal{L}_{|_{\Delta}} ).
\end{displaymath}
En utilisant  la formule de K\"unneth et le fait que $H^0(J^1,\mathcal{L}) \otimes H^0(J^1,\mathcal{L}) \simeq
 \bigwedge^2H^0(J^1,\mathcal{L})
\oplus S^2H^0(J^1,\mathcal{L})$ on observe que l'espace  $H^0(J^1\times J^1,(\mathcal{L}\boxtimes\mathcal{L})\otimes
 \mathcal{I})_+$ est contenu dans $S^2H^0(J^1,\mathcal{L})$. Donc la suite ci-dessus restreinte \`a la partie
$i$-invariante donne la suite exacte
\begin{displaymath}
0 \longrightarrow H^0(J^1\times J^1,(\mathcal{L}\boxtimes\mathcal{L}) \otimes \mathcal{I})_+ \longrightarrow
S^2H^0(J^1,\mathcal{L})
{\longrightarrow} H^0(J^1,\mathcal{L}^2) \longrightarrow 0 ,
\end{displaymath}
o\`u la surjectivit\'e de la multiplication est assur\'ee par la proposition 7.3.4.~\cite{CAV}. En calculant les 
dimensions de $S^2H^0(J^1,\mathcal{L})$ et
$H^0(J^1,\mathcal{L}^2)$ on trouve que la dimension de  $H^0(J^1\times J^1,(\mathcal{L}\boxtimes\mathcal{L}) \otimes
 \mathcal{I})_+$ est \'egal \`a
45-36=9. Ainsi on obtient que $\Gamma_{\varphi}=\Gamma_{\psi}= H^0(J^1\times J^1,\mathcal{L}\boxtimes\mathcal{L}
(-\Delta))_+ $, donc
 $\varphi$ et $\psi$ sont d\'efinies par le m\^eme syst\`eme lin\'eaire. Il existe donc un isomorphisme $\gamma$ tel que 
le diagramme
\begin{eqnarray*}
\shorthandoff{;:!?}
\xymatrix@!0 @R=10mm @C=1.8cm{
              & \mathbb{P}(V^*) \ar[dd]^{\gamma}_{\wr} \\
	J^1\times J^1  \ar@{-->}[ru]^{\varphi} \ar@{-->}[rd]^{\psi} &  \\
             & \mathbb{P}(V^*)
     }
\end{eqnarray*}
commute. On sait  que $V^*=H^0(J^1,\mathcal{L})$ est la seule repr\'esentation irr\'eductible du groupe de Heisenberg 
o\`u ${\bf\mu_3}$ agit par
multiplication. Puisque $\varphi$ et $\psi$ sont $\He$-\'equivariantes, elles d\'efinissent des repr\'esentations
 irr\'eductibles de
$\He$. Par le lemme de Schur l'application lin\'eaire associ\'ee \`a $\gamma$ est une homoth\'etie. Donc $\gamma=$Id,
ce qui ach\`eve la preuve.
\end{Dem}

\section{Sextiques planes} \label{sextiques_planes}

Consid\'erons une courbe elliptique $E$ et le plongement  $\varphi_3 : E \hookrightarrow
\mathbb{P}^2$ d\'efini par le fibr\'e en droites $\mathcal{O}_E(3\cdot 0)$. Soit  $\He \simeq {\bf\mu_3}
\times \mathbb{Z}/3
\times \widehat{\mathbb{Z}/3}$ le groupe de Heisenberg associ\'e \`a la polarisation $\mathcal{O}_E(3\cdot
 0)$. Comme l'application
$\varphi_3$ est $\He$-\'equivariante, l'image est une cubique plane  $\He$-invariante. On peut choisir des
coordonn\'ees
$X_0,X_1, X_2$ de $\mathbb{P}^2$ telles que $\He$ agit par
\begin{equation}  \label{eq:action}
(t,x,x^*)X_b= tx^*(b-x)X_{b-x}
\end{equation}
pour tout $  (t,x,x^*)\in \He $.
Par rapport \`a ces coordonn\'ees, l'image de $\varphi_M$ est donn\'ee par l'\'equation
\begin{displaymath}
f_{\lambda}:= X_0^3+ X_1^3+X_2^3 - 3\lambda X_0X_1X_2=0,
\end{displaymath}
connue comme la forme de Hesse de la cubique plane. Comme $\He$  agit sur les deriv\'ees partielles de
 $f_{\lambda}$, les neuf
points d'inflexion de la cubique forment une orbite du groupe $\He$.
\\
L'image de l'application
\begin{displaymath}
\mathcal{D}: [X_0, X_1, X_2] \mapsto \left[\frac{\partial f_{\lambda}}{\partial X_0} , \frac{\partial
 f_{\lambda}}{\partial X_1} ,
\frac{\partial f_{\lambda}}{\partial X_2} \right]
\end{displaymath}
est par d\'efinition la courbe duale de $E$. Dans le cas o\`u $E$ est lisse, la courbe duale est une sextique plane avec 
9 points de rebroussement comme points
singuliers, images des 9 points d'inflexion. Puisque cette application est $\He$-\'equivariante, la sextique duale est
 invariante par
$\He$ et ses points de rebroussement forment aussi une orbite du groupe de Heisenberg. Soient $Y_i$ pour $i=0,1,2$,
les coordonn\'ees homog\`enes de $\mathbb{P}^{2*}$ induites par l'application $\mathcal{D}$. \\
Une fois qu'on a calcul\'e la dimension de l'espace de sextiques planes $\He$-invariantes,  on peut donner une base
 explicite pour cet espace, \`a savoir
\\
$$ Y_0^6 +Y_1^6+  Y_2^6, \  Y_0^3Y_1^3 + Y_0^3Y_2^3 + Y_1^3Y_2^3,  \ Y_0Y_1Y_2( Y_0^3 +Y_1^3+  Y_2^3) , \  Y_0^2Y_1^2Y_2^2. $$\\
On \'ecrit donc l'\'equation d'une sextique plane g\'en\'erale  $\He$-invariante comme suit \\
$$ S := Y_0^6 +Y_1^6+  Y_2^6 + a_1(  Y_0^3Y_1^3 + Y_0^3Y_2^3 + Y_1^3Y_2^3) + $$
$$a_2 Y_0Y_1Y_2( Y_0^3 +Y_1^3+  Y_2^3) +a_3 Y_0^2Y_1^2Y_2^2=0.$$
Pour trouver les coefficients $a_i$ en termes du param\`etre $\lambda$ on utilisera les propri\'et\'es des points de
 rebroussement de la sextique.  \\
Les points d'inflexion de la cubique sont les points d'intersection de $E$ avec le Hessien de $E$. Un simple calcul 
exhibe les points
d'inflexions comme l'orbite du point  $p=(0:1:-1)$. Alors les points de rebroussement de la courbe duale sont  donn\'es 
par l'orbite du point
$\mathcal{D}(p)=(\lambda:1:1)$.  On obtient deux  \'equations lin\'eairement ind\'ependantes en \'evaluant les 
d\'eriv\'ees partielles $\frac{\partial F}{\partial Y_0}$
et $ \frac{\partial F}{\partial Y_1}$ au point singulier $(\lambda:1:1)$.  Pour donner une troisi\`eme \'equation, on
 consid\`ere la droite $\ell$
tangente \`a la sextique avec multiplicit\'e 3 en $\mathcal{D}(p)$. En utilisant le fait que cette droite est l'image 
par $\mathcal{D}$ du point
d'inflexion en $p$, on trouve que $Y_1=Y_2$ est l'\'equation qui d\'efinit $\ell$ dans $\mathbb{P}^{2*}$. Donc,  
la d\'eriv\'ee seconde du polyn\^ome
\begin{displaymath}
S_{|_{Y_1=Y_2=1}} = Y_0^6 + 2 + a_1(2Y_0^3 +1) + a_2Y_0(Y_0^3 +2) + a_3Y_0^2 =0,
\end{displaymath}
doit s'annuler en $Y_0=\lambda$. On obtient ainsi le syst\`eme d'\'equations lin\'eaires
$$\begin{array}{r}
6\lambda^5 + 6 a_1\lambda^2 +  a_2(4\lambda^3+2) + 2 a_3\lambda  =0 \\
 6 + 3 a_1(\lambda^3+1) +  a_2\lambda(\lambda^3+5) + 2 a_3\lambda^2  =0  \\
9 a_1\lambda^2 +  a_2(4\lambda^3+5) + 4 a_3\lambda  =0,
\end{array} $$
qui a comme solution les valeurs $a_1= 4\lambda^3 -2, a_2= -6\lambda^2, a_3= -3\lambda(\lambda^3 -4)$. Ceci prouve le
lemme suivant
\begin{Lem}
L'\'equation de $E^*$ dans $\mathbb{P}^{2*}$ est \\
 $$
F_{\lambda} := Y_0^6 +Y_1^6+  Y_2^6 +(4\lambda^3 -2)(  Y_0^3Y_1^3 + Y_0^3Y_2^3 + Y_1^3Y_2^3) $$
$$ -6\lambda^2Y_0Y_1Y_2( Y_0^3 +Y_1^3+  Y_2^3) -3\lambda(\lambda^3 -4) Y_0^2Y_1^2Y_2^2=0.$$
\end{Lem}
\begin{Cor} \label{duale}
Soient $E$ et $E'$ deux cubiques lisses dans $\mathbb{P}^2$ telles que leurs courbes duales ont le m\^eme ensemble
de points de rebroussement. Alors $E^* = (E')^*$. 
%La courbe duale d'une cubique lisse est uniquement caracteris\'ee par ses points de rebroussement.\\
\end{Cor}

\section{Le produit des courbes elliptiques}

On pose $A:= JC$.
Soit $\eta \in A[3], \eta\neq 0$ et soit $f: C_{\eta} \rightarrow C$ le rev\^etement \'etale 3-cyclique associ\'e \`a
$\eta$. Rappelons que
\begin{displaymath}
Gal(C_{\eta}/\mathbb{P}^1) = \langle \sigma, j \mid \sigma^3 =j^2 = 1 , \ j\sigma =\sigma^2 j \rangle 
\end{displaymath}
avec $j: C_{\eta}\rightarrow C_{\eta}$ un rel\`evement de l'involution hyperelliptique $\iota$ sur $C$. On a montr\'e
 (c.f.~\ref{prym}) que la
vari\'et\'e de Prym $P= Prym(C_{\eta}/C)$ est isomorphe au produit de courbes elliptiques $E_{\eta} \times E_{\eta}$, avec
$E_{\eta}:= C_{\eta}/ \langle j \rangle$. On \'ecrira $E$ \`a la place de $E_{\eta}$ s'il n'y a pas de confusion. 
Soit $M$ la polarisation sur $P\simeq E \times E$ de type (1,3) induite par la polarisation principale de $JC_{\eta}$.
Les automorphismes $\sigma$ et $j$ induits sur $JC_{\eta}$ se restreignent \`a des automorphismes de $P$.
\begin{Lem}
Les automorphismes $\sigma$ et $j$ sur $E\times E$ prennent la forme matricielle
$$ T = \left(
\begin{array}{cc}
0& -1\\
1 & -1
\end{array} \right), \quad
\widetilde{J}= \left(
\begin{array}{cc}
1& -1\\
0 & -1
\end{array} \right)
$$
respectivement.
\end{Lem}
\begin{Dem}
On consid\`ere l'isomorphisme $\psi : E\times E \stackrel{\sim}{\rightarrow} P$, donn\'e par $\psi (x,y) =x+\sigma y$. Puisque $\sigma \in
Gal(C_{\eta}/C)$, il v\'erifie l'\'equation $\sigma^2 + \sigma +1 =0 $ dans $ P \subset \Ker \Nm_f $. En  appliquant $\sigma$ \`a  $x+\sigma y \in P$
on a
\begin{displaymath}
\sigma(x+\sigma y) = \sigma x+ \sigma^2y = -y - \sigma y + \sigma x.
\end{displaymath}
On en d\'eduit que $ T= {0 \ -1 \choose 1 \ -1}$.\\
En appliquant $j$ \`a  $x+\sigma y \in P$ on obtient
\begin{displaymath}
jx+j\sigma y = x+ \sigma^2y = (x-y) - \sigma y,
\end{displaymath}
puisque par d\'efinition $j$ est induit l'identit\'e sur $E$ (c.f. section~\ref{decomposition}). On a donc que $\widetilde{J}=
  {1 \ -1 \choose 0 \ -1}$.
\end{Dem}
Gr\^ace au r\'esultat suivant (Lemme 5.4 ~\cite{BL1}) on conna\^{\i}t la forme de la polarisation $M$ sur $E\times E$.
\begin{Lem}  \label{polarisation}
La polarisation principale $\Theta$ sur $JC_{\eta}$ induit la polarisation $\Xi$ sur $E\times E$ laquelle est invariante par l'action de $\sigma$.
De plus,
\begin{displaymath}
\Xi = [E \times \{ 0\} + \{ 0\} \times E + \Delta ]  \ \in \ H^2(E\times E, \mathbb{Z})
\end{displaymath}
o\`u $\Delta$ d\'esigne la diagonale dans $E \times E$.
\end{Lem}
En identifiant  $E$ avec sa vari\'et\'e duale $\widehat{E}$ via
la polarisation canonique, l'isog\'enie  $\phi_{M} : E \times E \rightarrow E \times E $ est donn\'ee par une matrice de la forme (c.f. section~\ref{polar_prym})
$$ \phi_M= \left(
\begin{array}{cc}
2& \beta\\
\beta & 2
\end{array} \right) ,
$$
o\`u $\beta$ est un automorphisme  de $E$. Vu que $M$ est $\sigma$-invariant on a
$$\phi_M \circ T^{-1} = \widehat{T} \circ \phi_M,$$ avec
$\widehat{T}= {}^tT$ l'application duale vue comme automorphisme de $E \times E$. On en d\'eduit $\beta =-1$. On a donc
\begin{displaymath}
K(M)= \Ker \phi_M = \{ (x,-x) \mid 3x=0\}.
\end{displaymath}
Notons $J$ l'involution sur $E\times E$ don\'ee par
$$J = -\widetilde{J}T^2 =  \left(
\begin{array}{cc}
0& -1\\
-1 & 0
\end{array} \right) .$$
\\
\begin{Rmqs} \\
1) Les matrices $T$ et $J$ satisfont les relations $T^3=J^2=1$ et $TJ=JT^2$, donc elles engendrent le groupe sym\'etrique $\mathfrak{S}_3$.\\
2) On a $K(M)=\F (T) \subset \F (J)= \mathcal{A}$, o\`u $\mathcal{A}$ d\'esigne l'anti-diagonale dans $E\times E$.
\end{Rmqs}
\begin{Lem}
Les automorphismes $T$ et $J$ induisent l'identit\'e dans $H^0(M)$.
\end{Lem}
\begin{Dem}
Soit $\vartheta \neq 0$ une section globale du fibr\'e $\mathcal{O}_E(0)$. Alors
\begin{displaymath}
\Theta=p_1^*\vartheta \cdot p_2^*\vartheta \cdot \alpha^* \vartheta \ \in \ H^0(M),
\end{displaymath}
o\`u les $p_i: E\times E \rightarrow E$ sont les projections et $\alpha: E\times E \rightarrow E$ est donn\'ee par $\alpha(x,y)=x-y$. On a donc
\begin{displaymath}
T^* \Theta= T^*p_1^*\vartheta \cdot T^* p_2^*\vartheta \cdot T^*\alpha^* \vartheta = (-p_2)^*\vartheta \cdot \alpha^*\vartheta \cdot
(-p_1)^*\vartheta = \Theta.
\end{displaymath}
Pour tout $x \in K(M)= \F (T)$ on a
\begin{displaymath}
T^* t_x^*\Theta= t_x^*T^* \Theta  = t_x^*\Theta .
\end{displaymath}
Par le th\'eor\`eme de Stone-von Neuman les sections $ t_x^*\Theta, \ x \in K(M) $ engendrent l'espace vectoriel $H^0(M)$, donc $T$ induit
l'identit\'e  sur $H^0(M)$. De m\^eme,
 \begin{displaymath}
J^* \Theta= J^*p_1^*\vartheta \cdot J^* p_2^*\vartheta \cdot J^*\alpha^* \vartheta = (-p_2)^*\vartheta \cdot (-p_1)^*\vartheta \cdot
(-\alpha)^*\vartheta = \Theta,
\end{displaymath}
puisque le fibr\'e $\mathcal{O}_E(0)$ est sym\'etrique.
Comme $ K(M)\subset \F(J)$ $t_x$ et $J$ commutent pour tout $x\in K(M)$, alors $J^*t_x^* \Theta =
 t_x^*\Theta$. On conclut que $J$ induit l'identit\'e sur $H^0(M)$.
\end{Dem}
Consid\'erons l'application $\varphi_M: E\times E \rightarrow |M|^*\simeq\mathbb{P}^2$, associ\'ee \`a $M$. Alors
$\varphi_M$ est une application
surjective de degr\'e 6, ramifi\'ee le long d'une courbe plane de degr\'e 18 (c.f. example 10.1.5 \cite{CAV}). Le lemme
pr\'ec\'edent et le fait que
$\Deg \varphi_M = |\mathfrak{S}_3|=6$ impliquent le corollaire
\begin{Cor}
L'application $\varphi_M$ est un rev\^etement de Galois avec groupe de Galois $\mathfrak{S}_3=\langle T,J \rangle$.
\end{Cor}
Donc $\varphi_M$ se factorise par
\begin{displaymath}
E\times E \stackrel{\varphi_T}{\longrightarrow}Y  \stackrel{\varphi_J}{\longrightarrow} \mathbb{P}^2,
\end{displaymath}
o\`u  $\varphi_T$ et $\varphi_J$ sont des rev\^etements de degr\'es 3 et 2 respectivement. La surface $Y$ a comme
singularit\'es les 9 points
correspondants aux points fixes par $T$. Donc le rev\^etement double $\varphi_J: Y \rightarrow \mathbb{P}^2$ est
ramifi\'e sur une courbe
$D\subset \mathbb{P}^2$ avec 9 singularit\'es. Comme  $\varphi_T$ est \'etale en codimension 1, le diviseur de
ramification de $\varphi_M$ est 3$D$ avec $\Deg D= 6$. 
Observons que la restriction de $\varphi_M$ \`a l'anti-diagonale donne une
bijection entre $E\simeq \mathcal{A}$ et $D$, et que $\varphi_M$ envoie les 9 points de $K(M)= \Fix (T)$
dans le points singuliers de $D$.\\
\\
On consid\`ere  $E$ plong\'ee dans $\mathbb{P}^2 \simeq |3\cdot 0|^*$. Soit
$\Lambda : E \times E \rightarrow \mathbb{P}^{2*}$ l'application qui envoie $(x,y) \in E\times E$ dans 
$\ell_{x,-y}$, la droite qui relie $x$ \`a $-y$ et soit $\delta : E \rightarrow E \times E$, l'application d\'efinie par
$x \mapsto (x,-x)$.
\begin{Lem}
Il existe un isomorphisme $\lambda: \mathbb{P}^{2*} \stackrel{\sim}{\longrightarrow} \mathbb{P}^2$ tel que le
diagramme suivant commute
\begin{eqnarray}   \label{diagramme}
\shorthandoff{;:!?}
\xymatrix{
	E  \ar@{^{(}->}[r]^{\delta} \ar@{^{(}->}[d]^{\varphi_{3\cdot 0}} & E \times E \ar[d]^{\Lambda} \ar[r]^{ \ \varphi_M} &
    \mathbb{P}^2  \\
	 \mathbb{P}^2 \ar[r]^{\mathcal{D}} & \mathbb{P}^{2*} \ar[ru]_{\lambda} &  
     }
\end{eqnarray}
o\`u $\mathcal{D}$ est l'application duale.
\end{Lem}
\begin{Dem}
Il suffit de montrer qu'il existe une droite $\ell \subset \mathbb{P}^{2*}$ telle que $\Lambda^{-1}(\ell)= E \times \{ 0\}
+ \{ 0\} \times E + \Delta = \Xi$, puisque $M= \mathcal{O}_{E\times E}(\Xi)$ (c.f. lemme \ref{polarisation}). Prenons
$\ell$ la droite d'incidence de l'origine $0 \in E \subset \mathbb{P}^2$. On a
\begin{eqnarray*}
\Lambda^{-1}(\ell) & = & \{ (x,y) \in E\times E \mid 0 \in \ell_{x,-y} \}\\
&=& E \times \{ 0\} + \{ 0\} \times E + \{ (x,y)\mid (x) + (-y) + 0 \equiv 3 \cdot 0\} \\
&=&  E \times \{ 0\} + \{ 0\} \times E + \Delta 
\end{eqnarray*}
En particulier, l'image de $(x, -x)$ par $\Lambda$ est la droite tangente \`a $E$ en $x$.
\end{Dem}
\begin{Cor}
La sextique $D$ est isomorphe \`a la duale de la cubique $E$.
\end{Cor}
Pour chaque $\tilde{\eta} \in \He$ d'image non-nulle dans $A[3]$, on note $\mathbb{P}_{\tilde{\eta}}^2$, le
projectivis\'e de l'espace de points fixes par l'action de $\tilde{\eta}$. On utilisera cette notation tant\^ot 
pour l'action de $\tilde{\eta}$ sur $V$, tant\^ot pour l'action de $\tilde{\eta}$ sur $V^*$.
\begin{Lem}
Le diagramme
\begin{eqnarray*}
\shorthandoff{;:!?}
\xymatrix{
	P  \ar[rd]^{\varphi_M} \ar[d] &  \\
          P/\langle \sigma\rangle \ar@{^{(}->}[d]_{f_*} \ar[r]^{\theta_|} &  \mathbb{P}^2_{\tilde{\eta}}
        \ar@{^{(}->}[d]\\
         \su \ar[r]^{\theta} & |3\Theta|
     }
\end{eqnarray*}
commute.
\end{Lem}
\begin{Dem}
Soit $L \in P$. Par la formule de projection on a  
\begin{displaymath}
f_{*}L \otimes \eta \simeq f_{*}( L \otimes f^*\eta) \simeq f_{*}(L \otimes \mathcal{O}_{C_{\eta}}) \simeq 
f_{*}L.
\end{displaymath}
Ainsi l'image de $P/\langle \sigma \rangle $ dans $|3\Theta|$ est fix\'ee par l'action de $\eta$. Il existe
donc $\tilde{\eta} \in \He$, rel\`evement de $\eta$, tel que l'image de $P/\langle \sigma \rangle $ est 
contenue dans le plan $\mathbb{P}_{\tilde{\eta}}^2$.
\\
Soit $L \in P$ un fibr\'e en droites et $(x,y)$ son point correspondant dans $E \times E$ via l'isomorphisme 
$P \simeq E \times E$. On a 
\begin{displaymath}
\iota^*(f_{*}L)^* \simeq \iota^*(f_{*}L^{-1}) \simeq f_{*}(j^*L^{-1})
\end{displaymath}
puisque $j$ est un rel\`evement de l'involution hyperelliptique (c.f. section \ref{decomposition}). Le fibr\'e $j^*L^{-1} \in
P$ correspond au point $j(-x,-y)=(-x+y, y)\in E \times E$. Comme
\begin{displaymath}
\sigma(-x+y, y) = (-y, -x) = J(x,y),
\end{displaymath}  
on en d\'eduit que la restriction de l'involution $\mathfrak{i}$ \`a $P/ \langle \sigma \rangle \hookrightarrow
\su$ correspond \`a l'involution $J$ sur $E \times E/ \langle \sigma \rangle $ sous l'isomorphisme $P \simeq E 
\times E$.
\end{Dem}
\begin{Rmq} \label{distingue}
On observe que l'action de $\eta \in A[3]$ sur $\mathbb{P}^8$ fixe les trois plans
$\mathbb{P}^2_{\tilde{\eta}}$ correspondants aux trois rel\`evements de $\eta$ dans $\He$.
Un de ces trois plans est distingu\'e, celui qui contient le diviseur $\Theta + \Theta_{\eta} + \Theta_{-\eta}$, image de
$f_*\mathcal{O}_{C_{\eta}}$ par l'application $\theta$.
\end{Rmq}
\begin{Cor}
La sextique $D$ est le lieu de ramification de $\theta_{|_{P/\langle \sigma \rangle}}$, i.e. $D=\mathcal{B}\cap
\mathbb{P}^2_{\tilde{\eta}}$
\end{Cor}
\begin{Def}
Pour chaque $\eta \in A[3]$ non-nul, on note $D_{\eta}$ l'intersection de la sextique $\mathcal{B} \subset 
|3 \Theta|$ avec le plan $\mathbb{P}^2_{\tilde{\eta}}$ de points fixes par l'action de $\eta$ qui contient le
diviseur  $\Theta + \Theta_{\eta} + \Theta_{-\eta}$. 
\end{Def}
Par le diagramme (\ref{diagramme}) $D_{\eta}$ est la duale de $E_{\eta}$, o\`u $E_{\eta}$ est la courbe
elliptique de la d\'ecomposition $Prym (C_{\eta}/C) \simeq E_{\eta} \times E_{\eta}$. 
\\

D'un autre c\^ot\'e, on veut \'etudier les intersections  $\mathcal{C} \cap \mathbb{P}^2_{\tilde{\eta}}
\subset |3\Theta|^*$.
On choisit les coordonn\'ees sur $(\mathbb{Z}/3)^2$ de fa\c con que l'image de $\tilde{\eta}$ dans $A[3]\simeq
(\mathbb{Z}/3)^4$ soit $(00,10)$. Il suffit donc d'analyser l'intersection de $\mathcal{C}$ avec le plan de points fixes
correspondant \`a $\tilde{\eta}=(1,00,10) \in \He$. On notera ce plan $\mathbb{P}^2$, s'il n'y a pas de confusion.
Le plan $\mathbb{P}^2$ est donc d\'efini par les \'equations $X_{10}=X_{11}=X_{12}=X_{20}=X_{21}=X_{22}=0 $ (c.f.
section~\ref{coble's_cubic}).
On note $X_i := X_{0i}, \ i=0,1,2$ les coordonn\'ees homog\`enes de $\mathbb{P}^2$. Consi\'erons $F=0$ l'\'equation 
qui d\'efinit $\mathcal{C}$ (c.f.\'equation (\ref{cubique})). On trouve que
\begin{displaymath}
F_{\eta}:= \beta_0 (X_0^3 + X_1^3 + X_2^3) +3\beta_1 X_0X_1X_2=0
\end{displaymath}
est l'\'equation qui d\'efinit $\mathcal{C} \cap \mathbb{P}^2$, o\`u $\beta_0$ et $\beta_1$ sont des scalaires qui d\'ependent de $J^1$. On v\'erifie ais\'ement
que
\begin{equation} \label{partielles}
\frac{\partial F_{\eta}}{\partial X_i} = \frac{\partial F}{\partial X_{0i}}_{|_{\mathbb{P}^2}},  \qquad
\frac{\partial F}{\partial X_b}_{|_{\mathbb{P}^2}}=3 {Q_b}_{|_{\mathbb{P}^2}} = 0,
\end{equation}
pour $i=0,1,2$ et $b \notin \{ 00, 01, 02 \}$. Ceci implique que $F_{\eta} \neq 0$ et que
\begin{displaymath}
\Sing (\mathcal{C} \cap \mathbb{P}^2) = \mathbb{P}^2 \cap \Sing \mathcal{C}.
\end{displaymath}
De fa\c con g\'en\'erale,
on a que $\Sing (\mathcal{C} \cap \mathbb{P}^2_{\tilde{\eta}}) = \mathbb{P}^2_{\tilde{\eta}} \cap J^1 = \emptyset $ pour
 tout $\tilde{\eta} \in \He$, car $\eta$ agit par
translation sur $J^1$ et donc $J^1$ ne contient pas de points fixes par $\eta$. Par cons\'equent, les cubiques
$\mathcal{C} \cap \mathbb{P}^2_{\tilde{\eta}}$ sont lisses et contenues dans $\mathbb{P}^2_{\tilde{\eta}} \cap
\mathcal{C}_{sm}$. Les \'equations (\ref{partielles})  assurent  la commutativit\'e du diagramme suivant
\begin{eqnarray*}
\shorthandoff{;:!?}
\xymatrix{
	\mathbb{P}^2_{\tilde{\eta}}  \ar@{^{(}->}[d]^i \ar[r]^{\mathcal{D}_|} & \mathbb{P}^{2*}_{\tilde{{\eta}}}
     \ar@{^{(}->}[d]^i \\
	 \mathbb{P}^8 \ar@{-->}[r]^{\mathcal{D}} & \mathbb{P}^{8*}
     }
\end{eqnarray*}
On a ainsi
\begin{displaymath}
i (\mathcal{C} \cap \mathbb{P}^2_{\tilde{\eta}})^* = \overline{(i\circ \mathcal{D}_|) (\mathcal{C} \cap
 \mathbb{P}^2_{\tilde{\eta}})} =
\overline{(\mathcal{D} \circ i) (\mathcal{C} \cap \mathbb{P}^2_{\tilde{\eta}})} \subset 
\overline{\mathcal{D}(\mathcal{C}_{sm})} =\mathcal{C}^*,
\end{displaymath}
d'o\`u
\begin{displaymath}
(\mathcal{C} \cap \mathbb{P}^2_{\tilde{\eta}})^* \subset \mathcal{C}^* \cap  \mathbb{P}^{2*}_{\tilde{\eta}}.
\end{displaymath}
On a donc (c.f. section \ref{degre_duale})
\begin{equation} \label{grado}
6 = \Deg (\mathcal{C} \cap \mathbb{P}^2_{\tilde{\eta}})^*  \leq \Deg (\mathcal{C}^* \cap  \mathbb{P}^{2*}_{\tilde{\eta}})
= \Deg \mathcal{C}^* \leq 6. 
\end{equation}
Il en d\'ecoule que $\Deg \mathcal{C}^* = 6$ et donc l'inclusion ci-dessus est en fait une \'egalit\'e. 
\\
\\
Pour tout $\tilde{\eta}=(t,x,x^*) \in \He$ et pour tout $b\in (\mathbb{Z}/3)^2$ on a 
\begin{displaymath}
\tilde{\eta}\cdot Q_b = t^2x^*(2b-2x) \cdot Q_{b-x}.
\end{displaymath}
Donc, pour tout $p \in V$, $\mathcal{D}(\tilde{\eta}p)= (t^2,x,2x^*)\mathcal{D}(p)$. On en d\'eduit
que $\mathbb{P}^{2*}_{\tilde{\eta}}= \mathbb{P}(\Fix (\tilde{\eta}')) \subset
|3\Theta|$, avec $\tilde{\eta}' = (t^2,x,2x^*)$.
\\
\begin{Def}
Pour chaque $\eta \in A[3]$ non-nul, on note $\widetilde{D}_{\eta}$ l'intersection de la sextique $\mathcal{C}^*$ 
avec le plan $\mathbb{P}^2_{\tilde{\eta}} \subset |3\Theta|$ de points fixes par l'action de $\eta$ qui contient le
diviseur  $\Theta + \Theta_{\eta} + \Theta_{-\eta}$. 
\end{Def}
\begin{Rmq} On observe que les cubiques $\mathcal{C}\cap \mathbb{P}^2_{\tilde{\eta}}$ sont isomorphes pour les
$\tilde{\eta}=(t,x,x^*)$ de m\^eme image dans $A[3]$. On  a donc une copie de la courbe $\widetilde{D}_{\eta'}$
dans chacun des plans $\mathbb{P}^{2*}_{\tilde{\eta}} = \mathbb{P}(\Fix (\tilde{\eta}'))$, les plans de points fixes par
l'action de ${\eta}'=(x, 2x^*)$.
\end{Rmq}
\\
\begin{Pro}  \label{sextiques}
Les sextiques $D_{\eta}$ et  $\widetilde{D}_{\eta}$ co\"{\i}ncident.
\end{Pro}
\begin{Dem}
\'Etant donn\'e que $D_{\eta}$ et $\widetilde{D}_{\eta}$ sont des courbes duales de cubiques lisses planes, elles sont 
des sextiques avec 9
points de rebroussement correspondants aux 9 points d'inflexion. Rappelons qu'on a d\'efini l'application $\psi: J^1\times
J^1 \dashrightarrow |3\Theta|$ par
\begin{displaymath}
(a,b) \mapsto \Theta_x + \Theta_y +  \Theta_{-x-y},
\end{displaymath}
o\`u $\Theta_x$ et $\Theta_y $ sont les seuls translat\'es du diviseur th\^eta passant par $a$ et $b$. Par le th\'eor\`eme
\ref{secantes} on a
\begin{displaymath}
\im \psi =\im \varphi \subset \mathcal{C}^*.
\end{displaymath}
De plus, par la remarque du lemme ~\ref{singularite} l'image de  $\psi$ est contenue dans le lieu singulier de $ \mathcal{C}^*$. On a donc
\begin{displaymath}
\im \psi \cap  \mathbb{P}^2_{\tilde{\eta}} \subset \Sing \mathcal{C}^*\cap  \mathbb{P}^2_{\tilde{\eta}} \subseteq \Sing
(\mathcal{C}^*\cap \mathbb{P}^2_{\tilde{\eta}} ) = \Sing \widetilde{D}_{\eta}.
\end{displaymath}
Par ailleurs, on a $\im \psi\subset \Sing\mathcal{B}$ puisque les diviseurs de la forme $\Theta_x + \Theta_y +  \Theta_{-x-y}$ 
proviennent de fibr\'es vectoriels d\'ecomposables, et donc des points singuliers de $\su$.  Ainsi on obtient
\begin{displaymath}
\im \psi \cap  \mathbb{P}^2_{\tilde{\eta}} \subset \Sing \mathcal{B} \cap  \mathbb{P}^2_{\tilde{\eta}}  \subseteq \Sing
 (\mathcal{B} \cap \mathbb{P}^2_{\tilde{\eta}} ) = \Sing D_{\eta}.
\end{displaymath}
Soit $K_{\eta}:= \langle \eta \rangle^{\bot} /  \langle \eta \rangle \simeq  (\mathbb{Z}/3) \times
(\widehat{\mathbb{Z}/3})$. Soient $X_0, X_1, X_2$ des coordonn\'ees de $\mathbb{P}^2_{\tilde{\eta}}$, elles
sont donc fix\'ees par 
l'action de $\tilde{\eta}=(t,\eta) \in \He$. Observons qu'un \'el\'ement $(s,\lambda) \in \He$, pr\'eserve
 $\mathbb{P}^2_{\tilde{\eta}}$ si et seulement si
\begin{displaymath}
(t, \eta)((s, \lambda)\cdot X_i) =((s, \lambda)\cdot X_i) = ((s, \lambda)(t, \eta)) \cdot X_i,
\end{displaymath}
pour $i=0,1,2$, i.e. si et seulement si $(t, \eta)$ et $(s, \lambda)$ commutent , ce qui se traduit par $\langle 
\lambda , \eta \rangle =0$. Par
cons\'equent le groupe  ${\bf\mu_3} \times K_{\eta} \simeq \He$
agit sur  les coordonn\'ees de $\mathbb{P}^2_{\tilde{\eta}}$ comme dans (\ref{eq:action}). 
Observons que par construction $D_{\eta}$ et $\widetilde{D}_{\eta}$ sont $K_{\eta}-$invariantes (c.f. section
\ref{sextiques_planes}).  \\
Soit $x\in A[3]$, tel que $\langle x , \eta \rangle =0$. Le diviseur  $\Theta_x +\Theta_{\eta+x} +\Theta_{2\eta +x} \in
|3\Theta|$  est un point fixe pour l'action de $\eta$, on a donc
\begin{displaymath}
\Theta_x +\Theta_{\eta+x} +\Theta_{2\eta +x} \in \im \psi \cap \mathbb{P}^2_{\tilde{\eta}},
\end{displaymath}
et son orbite  par l'action de
$K_{\eta}$ est aussi contenue dans $\im \psi \cap \mathbb{P}^2_{\tilde{\eta}}$. Puisque $|K_{\eta}|=9$,  les points de
cette orbite sont les points singuliers de $D_{\eta}$ et $\widetilde{D}_{\eta}$.  Par  le corollaire ~\ref{duale} cela
suffit  pour conclure que les sextiques $D_{\eta}$ et $\widetilde{D}_{\eta}$ sont identiques.
\end{Dem}

\section[Sextiques invariantes]{Sextiques dans $\mathbb{P}^8$ invariantes par le groupe de Heisenberg}

Soit $V= H^0(J^1,\mathcal{L})$, un espace vectoriel de dimension 9 avec des coordonn\'ees $Z_b$ o\`u $b\in (\mathbb{Z}/3)^2$.
Consid\'erons l'action de $\He$
dans l'espace des sextiques $S^6V$. Comme cette action se factorise par le quotient ab\'elien A[3], on consid\'erera
l'espace des sextiques $A[3]$-invariantes $(S^6V)^{A[3]}$, au lieu de l'espace des sextiques $\He$-invariantes.  
\begin{Lem}
L'espace vectoriel  $(S^6V)^{A[3]}$ est de dimension 43.
\end{Lem}
\begin{Dem}
Par la th\'eorie des repr\'esentations (\S 2.2 ~\cite{F-H}) on a
\begin{equation} \label{trace}
\sum_{a \in A[3]} Tr (a \mid S^6V)=  \dim S^6_0 V \cdot |A[3] |.
\end{equation}
Par ailleurs, par le lemme 2 (\S 5, Ch. V ~\cite{Bou}) on a
\begin{displaymath}
\sum_{n=0}^{\infty} Tr (s \mid S^nV) T^n = \frac{1}{\Det (1-sT)},
\end{displaymath}
o\`u $s$ est un endomorphisme de  $V$. Dans notre cas l'action d'un \'el\'ement de $\tilde{a} \in \He$ sur $V$ a comme
valeurs propres les \'el\'ements de ${\bf\mu_3}$.
En fait, chaque racine cubique de l'unit\'e appara\^{\i}t avec multiplicit\'e 3. On a donc
\begin{eqnarray*}
\sum_{n=0}^{\infty} Tr (\tilde{a} \mid S^nV) T^n & = & \frac{1}{\Det (1-\tilde{a}T)}\\
& =& \frac{1}{(1-T^3)^3} \\
&=&  1+ 3T^3 + 6T^6 + \cdots.
\end{eqnarray*}
Pour $n=6$ on obtient que $Tr(a \mid S^6V)=6$ pour tout $a \in A[3]$ non nul. Donc de l'\'equation (\ref{trace})
il en r\'esulte
\begin{eqnarray*}
\sum_{a\in A[3] } Tr(a \mid S^6V) &=& 6\cdot 80 + \dim S^6V\\
&=& 480+3003\\
&=& \dim S_0^6V \cdot | A[3] |
\end{eqnarray*}
et on obtient que $\dim S_0^6V= 3483/81 =43$.
\end{Dem}

Pour chaque $\eta \in A[3]$, $\eta\neq 0$ on note $V_{\eta}$ le sous-espace de $V$ (de dimension 3) des points fixes 
par l'action de $\tilde{\eta}\in \He$, o\`u $\tilde{\eta}$ est le rel\`evement de $\eta$ tel que 
$\mathbb{P}^2_{\tilde{\eta}} \subset |3\Theta|$ est le plan distingu\'e (c.f. remarque de la page \pageref{distingue}). 
On observe que
$V_{\eta}=V_{-\eta}$. On consid\`ere le groupe $K_{\eta}= \langle \eta \rangle^{\bot} / \langle \eta \rangle $.
Puisque l'action de $\He \simeq {\bf\mu_3} \times K_{\eta}$ sur $S^6V_{\eta}$ se factorise par le quotient $K_{\eta}$ on
consid\'erera l'espace des sextiques $K_{\eta}$-invariantes au lieu de $(S^6V_{\eta})^{\He}$. Soit
\begin{displaymath}
res_{\eta}: (S^6V)^{A[3]}  \rightarrow  (S^6V_{\eta})^{K_{\eta}}
\end{displaymath}
l'application restriction. On d\'efinit l'application
\begin{displaymath}
\nu:= \sum_{\eta} res_{\eta}:  (S^6V)^{A[3]}  \rightarrow \bigoplus (S^6V_{\eta})^{K_{\eta}}
\end{displaymath}
o\`u la somme parcourt les $\eta \in A[3]$, $\eta \neq 0$ , modulo $\eta\sim -\eta$.  Le but est de d\'ecrire le noyau de cette application pour
savoir sous quelles conditions les restrictions aux sous-espaces des points fixes caract\'erisent compl\`etement une sextique invariante par le
groupe de Heisenberg. \\

On va d'abord donner une base de l'espace $ (S^6V)^{A[3]}$.  On note $K:= \{ (x,00) \mid x \in (\mathbb{Z}/3)^2 \}$ et
$\widehat{K}:= \{ (00,x^*) \mid x^* \in (\widehat{\mathbb{Z}/3})^2  \}$. Un mon\^ome de degr\'e 6, $Z_{\mu_1}
\cdots Z_{\mu_6}$, est invariant par l'action de $\widehat{K}$ si et seulement si $\sum \mu_i=0$. On peut donc construire des polyn\^omes
$A[3]$-invariants en faisant agir $K$ sur chacun des mon\^omes $\widehat{K}$-invariants et en prenant la somme sur tous les \'el\'ements de $K$.
Avec cette m\'ethode on obtient une base de  $ (S^6V)^{A[3]}$ comme suit
\begin{eqnarray*}
\sum_{\sigma \in (\mathbb{Z}/3)^2} Z_{\sigma}^6 & &\\
\sum_{\sigma} Z_{\sigma + \mu}^3Z_{\sigma}^3 & & \textrm{avec} \ \mu \neq 0, \  \textrm{modulo} \ \mu \sim -\mu   \\
\sum_{\sigma} Z_{\sigma}^3Z_{\sigma +\mu_1}Z_{\sigma + \mu_2} Z_{\sigma + \mu_3} && \textrm{avec } \ \sum \mu_i=0 \\
\frac{1}{3}\sum_{\sigma} Z_{\sigma}^2Z_{\sigma +\mu}^2Z_{\sigma + 2\mu}^2 & & \textrm{avec} \ \mu \neq 0  \\
\sum_{\sigma} Z_{\sigma}^2Z_{\sigma +\mu_1}Z_{\sigma + \mu_2} Z_{\sigma + \mu_1 +\mu_2}^2 & &\textrm{avec } \ \mu_1, \mu_2 \neq 0, \langle
\mu_1 \rangle \neq   \langle   \mu_2 \rangle\\
\sum_{\sigma} Z_{\sigma}^2Z_{\sigma +\mu_1}Z_{\sigma + \mu_1^{-1}} Z_{\sigma + \mu_2} Z_{\sigma +\mu_2^{-1}} & &\textrm{avec } \ \mu_1, \mu_2
\neq 0, \ \langle  \mu_1 \rangle \neq   \langle   \mu_2 \rangle\\
\frac{1}{3}\sum_{\sigma} \left( \prod_{i} Z_{\sigma + \mu_i} \right) &  & \textrm{avec } \ \sum_{i=1}^6 \mu_i=0, \mu_i\neq \mu_j \  \textrm{pour} \ i\neq j .
\end{eqnarray*}
On v\'erifie que ces \'el\'ements forment une base pour  $ (S^6V)^{A[3]}$ (cette base est \'ecrite explicitement dans
 l'annexe). On a utilis\'e le programme Maple pour calculer 
le rang et le noyau de l'application $\nu$ (c.f. Annexe). On a trouv\'e que $\nu$ est de rang 39 et qu'une base pour 
$\Ker \nu$ est donn\'ee par
\begin{eqnarray*}
w_1:=\sum_\sigma Z_{\sigma +20}^3Z_{\sigma }Z_{\sigma + 01} Z_{\sigma + 02}- \sum_\sigma Z_{\sigma +10}^3Z_{\sigma }Z_{\sigma + 01} Z_{\sigma + 02}
&& \\
w_2:=\sum_\sigma Z_{\sigma +02}^3Z_{\sigma }Z_{\sigma + 10} Z_{\sigma + 20}- \sum_\sigma Z_{\sigma +01}^3Z_{\sigma }Z_{\sigma + 10} Z_{\sigma + 20}
&& \\
w_3:=\sum_\sigma Z_{\sigma +02}^3Z_{\sigma }Z_{\sigma + 11} Z_{\sigma + 22}- \sum_\sigma Z_{\sigma +01}^3Z_{\sigma }Z_{\sigma + 11} Z_{\sigma + 22}
&& \\
w_4:=\sum_\sigma Z_{\sigma +20}^3Z_{\sigma }Z_{\sigma + 12} Z_{\sigma + 21}- \sum_\sigma Z_{\sigma + 10}^3Z_{\sigma }Z_{\sigma + 12} Z_{\sigma + 21}.
\end{eqnarray*}
L'involution $\iota^* : J^1 \rightarrow J^1 $ s'\'etend \`a une involution $\iota$ de $V$ agissant par
\begin{displaymath}
\iota: [Z_{00}, Z_{01}, \ldots , Z_{22}] \mapsto [Z_{00}, Z_{-(01)}, \ldots , Z_{-(22)}] = [Z_{00}, Z_{02}, \ldots , Z_{11}].
\end{displaymath}
Donc $\iota$ agit aussi sur l'espace des sextiques $(S^6V)^{A[3]}$. Soit  $(S^6V)^{A[3]}= W_+ \oplus W_-$ la d\'ecomposition en les
parties $\iota$-invariante et $\iota$-anti-invariante respectivement.\\
\\
\begin{Rmq}  \label{remarque}
On observe que $\iota (w_k) =-w_k$ pour $k=1,2,3$. Donc $\Ker \nu \subset W_-$, d'o\`u $W_+ \cap \Ker \nu =\{ 0\}$. \\
\end{Rmq}
\\
Pour tout $\eta \in A[3]$ et $E \in \mathcal{B}$ on a
\begin{displaymath}
\mathfrak{i}(E \otimes \eta) = \iota^*(E \otimes \eta)^* \simeq \iota^*E^* \otimes \iota^*\eta^{-1} \simeq E \otimes \eta, 
\end{displaymath}
d'o\`u $\mathcal{B}$ est $A[3]$-invariante. D'un autre cot\'e, puisque l'application duale $\mathcal{D}$ est 
$\He$-\'equivariante on a que la sextique $\mathcal{C}^*$ est aussi $A[3]$-invariante dans $\mathbb{P}(V) \simeq
\mathbb{P}^8$. Soient  $\mathsf{R}$ et $\mathsf{S}$  les polyn\^omes homog\`enes dans $(S^6V)^{A[3]}$ qui d\'efinissent $\mathcal{B}$ et $\mathcal{C}^*$ respectivement.
\begin{Pro} \label{poly}
Le polyn\^ome $\mathsf{R}$ est  $\iota$-invariant et $\mathsf{S} \in W_+ \cup W_-$.
\end{Pro}
\begin{Dem}
La cubique de Coble est d\'efinie par le polyn\^ome
\begin{displaymath}
F= \sum_{b\in K} X_bQ_b=0.
\end{displaymath}
En appliquant l'involution $\iota$ on obtient
\begin{displaymath}
\iota(\sum_b X_bQ_b) = \sum_b \iota(X_b)\iota(Q_b) = \sum_b X_{-b}Q_{-b}
\end{displaymath}
et donc $F$ est $\iota$-invariante. \'Etant donn\'e que l'application duale $\mathcal{D}$ est d\'efinie par
\begin{displaymath}
\mathcal{D}: x \mapsto [Q_{00} (x) : \cdots  : Q_{22}(x)]
\end{displaymath}
et $\iota Q_b (x) =Q_{-b} (x) = Q_b(\iota x)$, alors $\mathcal{D}$ est $\iota$-\'equivariante, d'o\`u $\mathsf{S} \in W_+ \cup W_-$. \\
\\
Soit $\Theta_E = \mathrm{Supp } \ \theta(E)$ avec $E\in \su$. On a montr\'e (c.f. proposition \ref{theta}) que
\begin{displaymath}
\iota^*(\Theta_E) = \Theta_{E^*}.
\end{displaymath}
Si $E\in \mathcal{B}= \{ E \in \su \mid \iota^* E \sim_{s} E^*\}$, alors  $E^*\in \mathcal{B}$ et donc $\mathcal{B}$ est $\iota$-invariante.
Par cons\'equent,  $\mathsf{R} \in W_+ \cup W_-$. L'action de l'involution $\iota$ sur $V$ induit une d\'ecomposition $V=V_+ \oplus V_-$,
en les sous-espaces associ\'es aux valeurs propres $\{ \pm 1\}$. Si on suppose $\mathsf{R} \in W_-$, alors la sextique $\mathcal{B}$ contient l'espace
$\mathbb{P}(V_+) \simeq \mathbb{P}^4$. En effet, si $p\in V_+$ alors
\begin{displaymath}
\mathsf{R}(p) = \mathsf{R}(\iota p) = \iota \mathsf{R}(p) = -\mathsf{R}(p),
\end{displaymath}
et $\mathsf{R}(p)=0$. On analysera l'involution
$\iota^*$ sur $\mathcal{B} \subset \su$ pour donner une estimation de la dimension du lieu de points fixes.  \\
\\
Soit $E \in  \mathcal{B}$ un fibr\'e vectoriel  $\iota$-invariant.  \\
{\it Cas 1.}  Supposons $E$ stable. On a donc $\iota^*E \simeq E^* \simeq E$, i.e. il existe un isomorphisme $u: E
\stackrel{\sim}{\rightarrow} E^*$.
Comme $E$ est stable, l'homomorphisme  $u^{-1} \circ {}^tu : E \rightarrow E$ est de la forme $\lambda \cdot$Id, avec $\lambda \neq 0 $. On a
$$ \Det u = \Det {}^tu = \lambda ^3 \Det u,$$
donc $\lambda^3 =1$; et puisque ${}^tu = \lambda\cdot\mathrm{Id} \circ u$, on a $\lambda^2=1$. Donc ${}^tu =u$.  On a ainsi que
$u$  induit une forme quadratique sur $E$, autrement dit, $E \in H^1(C, \mathcal{O}SO_3)$. Par la suite exacte (~\ref{suite}) on obtient
que la dimension du lieu des fibr\'es $E\in \mathcal{B}$ qui
sont stables et $\iota$-invariants est au plus la dimension de  $\mathcal{U}_C(2,0)$, l'espace de modules des fibr\'es vectoriels semi-stables de rang 2
et de degr\'e  0. Il en d\'ecoule
\begin{displaymath}
\dim \F (\iota^*) \cap  \mathcal{B}_s \leq 3.
\end{displaymath}
{\it Cas 2.}  Supposons $E\sim_{s} L_1 \oplus L_2 \oplus L_3$, avec $L_i $ des fibr\'es en droites. On a
\begin{displaymath}
\iota^*(L_1 \oplus L_2 \oplus L_3 )  \simeq  L_1^{-1} \oplus L_2^{-1} \oplus L_3^{-1}  \sim_{s} L_1 \oplus L_2 \oplus L_3.
\end{displaymath}
Supposons $L_1 \simeq L_2^{-1}$. Comme $\Det E \simeq L_1 \otimes L_2 \otimes L_3  \simeq \fs$ on a $L_3\simeq \fs$ et donc
$E \sim_{s} L_1 \oplus L_1^{-1} \oplus \fs$. Par contre, si pour tout $i=1,2,3$, $L_i\simeq L_i^{-1}$ alors $L_i \in JC[2]$, ce qui ne
rajoute qu'un nombre fini de points fixes. Donc dans ce cas la $\dim \F (\iota^*) \leq \dim$ Pic$(C)$ =2. \\
{\it Cas 3.} Supposons $E\sim_{s} L \oplus F$ avec $F\in \mathcal{U}_C(2,L^{-1})$ stable. On a
 \begin{displaymath}
\iota^*(L \oplus F) \simeq L^{-1} \oplus F^* \sim_{s} L \oplus F
\end{displaymath}
car  $E\in \F (\iota^*)$, ce qui implique $L \simeq L^{-1}$ et donc $L \in JC[2]$, et $F^* \simeq F$. En combinant avec
 l'estimation faite dans
le cas 2 on  obtient
\begin{displaymath}
\dim \F (\iota^*) \cap  \mathcal{B}_{ss} \leq \dim (\mathcal{U}_C(2,0) \cup \textrm{Pic}(C))=3.
\end{displaymath}
En conclusion, $\mathcal{B}$ ne contient pas l'espace $\mathbb{P}(V_+)\simeq \mathbb{P}^4$ et donc $\mathsf{R} \in W_+$.
\end{Dem}
\begin{Lem}  \label{orthogonaux}
Soient  $\eta, \lambda \in A[3] \smallsetminus \{ 0\}, \langle \eta \rangle \neq \langle \lambda \rangle$. Si $\langle 
\eta , \lambda \rangle =0$,
alors $\mathbb{P}^2_{\tilde{\eta}} \cap \mathbb{P}^2_{\tilde{\lambda}} \neq \emptyset $, avec $\tilde{\eta}$ et
$\tilde{\lambda}$ rel\`evements de $\eta$ et $\lambda$ respectivement.
\end{Lem}
\begin{Dem}
Supposons  $\eta, \lambda \in A[3] \smallsetminus \{ 0\}$ orthogonaux. Donc le groupe $K_{\eta, \lambda}$
engendr\'e par $\eta$ et $\lambda$, est un
sous-groupe maximal isotrope. Donc le groupe $K_{\tilde{\eta}, \tilde{\lambda}}$ engendr\'e par $\tilde{\eta}$ et
$\tilde{\lambda}$  est un rel\`evement de $K_{\eta, \lambda}$ dans $\He$, car $\tilde{\eta}$ et $\tilde{\lambda}$
commutent. La repr\'esentation du groupe de Heisenberg 
induit une repr\'esentation de  $K_{\tilde{\eta}, \tilde{\lambda}}$ dans
$V= H^0(J^1, \mathcal{L})$. Ainsi l'intersection de  $\mathbb{P}^2_{\tilde{\eta}}$ et  $\mathbb{P}^2_{\tilde{\lambda}}$ 
est donn\'ee par le projectivis\'e de l'espace
vectoriel $V^{K_{\tilde{\eta}, \tilde{\lambda}}}$, l'espace de sections invariantes par $\tilde{\eta}$ et
$\tilde{\lambda}$. Par la proposition 3 ~\cite{Mum1} on a
\begin{displaymath}
\dim V^{K_{\tilde{\eta}, \tilde{\lambda}}}= 1,
\end{displaymath}
ce qui montre le lemme.
\end{Dem}
\begin{Lem} \label{chaine}
Pour tout $\eta, \lambda \in A[3] \smallsetminus \{ 0\}$, il existe une suite de sous-espaces
\begin{displaymath}
V_{\eta}= V_{\beta_0}, V_{\beta_1}, \ldots V_{\beta_n} = V_{\lambda}
\end{displaymath}
avec $\beta_i \in A[3] \smallsetminus \{ 0\} $, telle que $V_{\beta_i} \cap V_{\beta_{i+1}} \neq \{ 0 \}$ pour tout $ 0
 \leq i \leq n-1$.
\end{Lem}
\begin{Dem}
Par le lemme ~\ref{orthogonaux} il suffit de montrer que  pour tout $\eta, \lambda \in A[3] \smallsetminus \{ 0\}$, il
 existe une suite $\eta= \beta_0, \beta_1,
 \ldots \beta_n = \lambda$ dans  $A[3] \smallsetminus \{ 0\}$ telle que $\beta_i \bot  \beta_{i+1}$ pour tout $ 0 \leq 
i \leq n-1$.  \\
Pour tout $y^*$ il existe $x\neq 0$ tel que $y^*(x)=0$, donc pour tout $(y,y^*) \in A[3]$ il
existe $(x,00) \in K$, avec $x\neq 0$  tel que
\begin{displaymath}
\langle (x,00), (y,y^*) \rangle = y^*(x)=0.
\end{displaymath}
Comme le sous-groupe $K$ est isotrope, cela suffit pour d\'emontrer la proposition.
\end{Dem}
\begin{Pro}  \label{scalaire}
Soient  $\tilde{\eta}$, $\tilde{\lambda} \in \He$ tels que leurs images dans $A[3]$ v\'erifient les conditions du lemme
\ref{orthogonaux}. Soit $E\in \su$ tel que $\theta(E)\in \mathbb{P}^2_{\tilde{\eta}} \cap \mathbb{P}^2_{\tilde{\lambda}} \subset |3\Theta|$. Alors $E$ n'appartient pas \`a  $\mathcal{B}$.
\end{Pro}
\begin{Dem}
On consid\`ere l'application \'etale de degr\'e 3
\begin{displaymath}
f:  C_{\eta} \longrightarrow C
\end{displaymath}
associ\'ee \`a $\eta \in A[3]\setminus \{ 0\}$ et soit $P= Prym(C_{\eta}/C)$, la vari\'et\'e de Prym munie de la polarisation $M$, induite par la polarisation
principale de $JC_{\eta}$. Soit $\sigma$ l'automorphisme de $C_{\eta}$ qui
engendre le groupe de Galois $Gal(C_{\eta}/C)$. \\
Pour tout fibr\'e $L \in P$ le fibr\'e
de rang 3 $f_* L$  est contenu dans $\su$ et est invariant par l'action de $\eta$. En effet,  par la formule de projection on a
\begin{displaymath}
f_* (L) \otimes \eta \simeq f_*( L \otimes f^* \eta) \simeq  f_*( L \otimes \mathcal{O}_{C_{\eta}}) \simeq f_* L.
\end{displaymath}
On a le diagramme commutatif suivant
\begin{eqnarray*}
\shorthandoff{;:!?}
\xymatrix{
	P  \ar[rd]^{\varphi_{M}} \ar[d] &  \\
          P/\langle \sigma\rangle \ar@{^{(}->}[d]_{f_*} \ar[r]^{\theta_|} &  \mathbb{P}^2_{\tilde{\eta}} \ar@{^{(}->}[d]\\
         \su \ar[r]^{\theta} & |3\Theta|
     }
\end{eqnarray*}
Ce diagramme montre que si   $\theta(E)\in \mathbb{P}^2_{\tilde{\eta}}$ alors $E$ est de la forme $f_*(L)$ avec $L $ un fibr\'e en droites dans la vari\'et\'e
de Prym $P$.  Supposons $E \in \su$ tel que $\theta(E)\in \mathbb{P}^2_{\tilde{\eta}} \cap \mathbb{P}^2_{\tilde{\lambda}} $. Donc $E\simeq f_*(L)$ avec $L\in P$
et $E$ est invariant par l'action de $\lambda$. En utilisant la formule de projection on obtient
\begin{displaymath}
E \simeq f_* (L) \simeq f_*(L) \otimes \lambda \simeq f_*( L \otimes f^* \lambda).
\end{displaymath}
Ceci implique que, soit  $L \otimes f^* \lambda \simeq \sigma^*L$, soit   $L \otimes f^* \lambda \simeq \sigma^{2*}L$. Par ailleurs,  on a
d\'emontr\'e  pr\'ec\'edemment que $P \simeq Z \times Z$, o\`u $Z$ est une courbe elliptique et que l'automorphisme $\sigma$ de $JC_{\eta}$
restreints \`a $P$ prend la forme
\begin{displaymath}
\sigma : (x,y) \mapsto (-y, x-y),    \qquad \forall x,y \in Z.
\end{displaymath}
On a que
\begin{displaymath}
\Gamma := \mathcal{A} + \sigma(\mathcal{A}) + \sigma^2(\mathcal{A}) \subset Z \times Z,
\end{displaymath}
avec $\mathcal{A}$ l'anti-diagonale dans $Z \times Z$,
est le lieu de ramification de l'application $\varphi_M : P \rightarrow \mathbb{P}^2_{\tilde{\eta}}$.  Soit $(x,y) \in Z\times Z$ le point correspondant au
fibr\'e $L$. Supposons $\sigma^* L \otimes L^{-1} \simeq f^* \lambda$. Donc le fibr\'e $f^*\lambda$ correspond au point
\begin{displaymath}
(\sigma - 1)(x,y) = (-y, x-y) -(x,y) = (-y-x, -2y).
\end{displaymath}
Observons que par construction  $\sigma^* f^*\lambda \simeq f^*\lambda$.\\
Supposons $E\in \mathcal{B}$, i.e. $\iota^*E \simeq E^*$. Alors $L$ appartient au lieu de ramification de $\varphi_M $, autrement dit,
$(x,y)\in \Gamma$.  \\
Supposons $(x,y) \in \mathcal{A}$, alors $y=-x$ et  $f^*\lambda$ correspond au point $(0, 2x)$.  On a
\begin{displaymath}
\sigma (0,2x) = (-2x, -2x)  = (0,2x),
\end{displaymath}
et donc $2x=0$,  ce qui se traduit par $f^*\lambda \simeq \mathcal{O}_{C_{\eta}}$.  Mais ceci est une contradiction avec l'hypoth\`ese $\langle \eta
\rangle \neq \langle \lambda \rangle $ puisque $ \Ker f^*= \langle \eta\rangle$.  \\
Si $(x,y)\in \sigma^i(\mathcal{A})$, avec $i=1,2$ alors $\sigma^{-i}(x,y)\in \mathcal{A}$.
En proc\'edant comme avant on obtient $(\sigma^{-i})^*f^*\lambda\simeq \mathcal{O}_{C_{\eta}}$ et on arrive \`a la m\^eme contradiction.
\\
De mani\`ere analogue, on d\'emontre que si $L \otimes f^* \lambda \simeq \sigma^{2*}L$ alors $f^*L \simeq \mathcal{O}_{C_{\eta}}$ et donc on
contredit l'hypoth\`ese $\langle \eta \rangle \neq \langle \lambda \rangle $.
\end{Dem}
\begin{The}
Les polyn\^omes $\mathsf{S}$ et $\mathsf{R}$ d\'efinissent la m\^eme hypersurface, i.e. $\mathcal{C}^* = \mathcal{B}$.
\end{The}
\begin{Dem}
Par la proposition ~\ref{sextiques}  on a
\begin{displaymath}
\mathcal{C}^* \cap \mathbb{P}^2_{\tilde{\eta}} = \mathcal{B}  \cap \mathbb{P}^2_{\tilde{\eta}},
\end{displaymath}
pour tout $\tilde{\eta} \in \He$ avec projection non-nulle dans $A[3]$, d'o\`u
\begin{displaymath}
\mathsf{S}_{|_{V_{\eta}}}  = \epsilon(\eta) \mathsf{R}_{|_{V_{\eta}}},
\end{displaymath}
avec  $\epsilon(\eta)$ un scalaire non-nul. Pour tout $\eta, \lambda \in A[3]\setminus \{ 0 \}$ tels que $V_{\eta} \cap
V_{\lambda} \neq \{ 0\}$ on a
\begin{displaymath}
\mathsf{S}_{|_{V_{\eta}\cap V_{\lambda}}}  = \epsilon(\eta)\mathsf{R}_{|_{V_{\eta} \cap V_{\lambda}}} = \epsilon(\lambda) \mathsf{R}_{|_{V_{\eta} \cap
V_{\lambda}}}.
\end{displaymath}
Par la proposition ~\ref{scalaire}  $\mathsf{R}_{|_{V_{\eta} \cap V_{\lambda}}} \neq 0$,
donc $\epsilon(\eta) =\epsilon(\lambda)$. Gr\^ace au lemme ~\ref{chaine},  on d\'eduit que le scalaire $\epsilon$ ne d\'epend pas de $\eta$.
On obtient donc l'\'egalit\'e
\begin{displaymath}
\mathsf{S}_{|_{V_{\eta}}}  - \epsilon \mathsf{R}_{|_{V_{\eta}}} =0 \quad \forall \eta \in A[3] \setminus \{  0\},
\end{displaymath}
i.e. $\mathsf{S} - \epsilon \mathsf{R}$ appartient au noyau de l'application  $\nu: (S^6V)^{A[3]}  \rightarrow \bigoplus (S^6V_{\eta})^{K_{\eta}}$. Donc, par la
remarque de la page~\pageref{remarque},  $\mathsf{S} - \epsilon \mathsf{R} \in W_-$. Si $\mathsf{S} \in W_ -$ alors $\epsilon \mathsf{R} 
= \mathsf{S} - (\mathsf{S} - \epsilon \mathsf{R}) \in
W_-$, ce qui contredit la proposition ~\ref{poly}. Donc $\mathsf{S} \in W_+$ et $\mathsf{S} - \epsilon \mathsf{R}$ est $\iota$-invariante.  On a donc
\begin{displaymath}
\mathsf{S}F - \epsilon \mathsf{R} \in W_+ \cap \Ker \nu = \{ 0 \}.
\end{displaymath}
Ce qui ach\`eve la d\'emonstration.
\end{Dem}

\chapter[Vari\'et\'es de Prym]{Vari\'et\'es de Prym associ\'ees aux rev\^etements $n$-cycliques d'une courbe hyperelliptique}

Soit $H$ une courbe projective, hyperelliptique et non singuli\`ere de genre $g$ et soit $JH \simeq \mathrm{Pic}^0(H)$
la jacobienne de $H$. On fixe un \'el\'ement $\eta \in JH$ d'ordre $n$. On note $f : C \rightarrow H$ le rev\^etement
\'etale $n$-cyclique associ\'e \`a $\eta$. Soit $JC$ la jacobienne de $C$.
\\
Soit $\Nm : JC \rightarrow JH$, $\sum n_ip_i \mapsto \sum n_if(p_i)$, l'application Norme de $f$.
On appelle vari\'et\'e de Prym $P=Prym (C/H)$ associ\'ee au rev\^etement $f$ la
composante neutre de $\Ker\Nm$. Si $\sigma : C \rightarrow C$ est l'automorphisme
d'ordre $n$ qui engendre le groupe de Galois $Gal(C/H)$, on peut \'ecrire $P= \im (1-\sigma).$
La vari\'et\'e $P$ est une sous-vari\'et\'e ab\'elienne de $JC$ de dimension $(n-1)(g-1)$ avec polarisation $\Xi$
induite par la polarisation principale de $JC$.\\
Le but de ces lignes est de montrer que lorsque $n$ n'est pas divisible par 4 la  vari\'et\'e de Prym $P$ est
isomorphe \`a un produit de jacobiennes.

\section{D\'ecomposition de la vari\'et\'e de Prym} \label{decomposition}
On pose $Y=f^*(JH)$. C'est une sous-vari\'et\'e ab\'elienne suppl\'ementaire de $P$ dans $JC$ de dimension $g$.
On peut aussi d\'ecrire $Y$ comme l'image de l'endomorphisme Norme $1+\sigma+ \cdots + \sigma^{n-1}$.
\\
\begin{Pro}
La polarisation $\Xi$ sur $P$ est du type  $(\underbrace{1,...,1}_{(n-2)(g-1)},\underbrace{n,...,n}_{g-1})$.
\end{Pro}
\begin{Dem}
Puisque $f:C \rightarrow H$ est \'etale l'application $f^*: JH \rightarrow JC$ n'est pas injective (~\cite{CAV}
11.4.3.). En fait, $\Ker f^*= \langle \eta \rangle$, un sous-groupe de $JH[n]$ d'ordre $n$. On a donc le diagramme
commutatif suivant
$$
\shorthandoff{;:!?}
\xymatrix @!0 @R=15mm @C=1.5cm{
       JH \ar[rr]^{f^*} \ar[rd]_h & & JC \\
       & JH/\langle \eta \rangle \simeq Y \ar[ru]_{i_Y} &
     }
$$
o\`u $h$ est une isog\'enie de degr\'e $n$. Soit $\Theta$ un diviseur th\^eta dans $JC$ et soit $M=\mathcal{O}_{Y}
(i_Y^* \Theta)$. On a $h^*M \simeq\mathcal{O}_{JH}(n\Theta_H) $ car $(f^*)^*\Theta \equiv n\Theta_H$
(\cite{CAV} 12.3.1.). D'apr\`es ~\cite{Mum3} (Lemme 2, pag. 232), $K(M) \simeq \langle \eta \rangle ^{\perp}/ \langle \eta
\rangle$, o\`u $\langle \eta \rangle ^{\perp}$ est l'orthogonale de $\langle \eta \rangle$ par rapport \`a la forme
de Weil $e^{h^*M}: K(h^*M) \times K(h^*M) \rightarrow \mathbb{C}^*$. Puisque $K(h^*M)\simeq (\mathbb{Z}/n)^{2g}$ on
obtient donc $K(M) \simeq (\mathbb{Z}/n)^{2(g-1)}$. Ainsi $M$ est une polarisation du type
$(1,n, \ldots,n)$ et par ~\cite{CAV} (Cor. 12.1.5.), $\Xi = i_P^* \Theta$ est du type $(1, \ldots,1,n, \ldots,n)$.
\end{Dem}
On note $i : H \rightarrow H $ l'involution hyperelliptique. Par construction on a $C = Spec ({\bf
\mathcal{A}})$, o\`u ${\bf \mathcal{A}}:=\fs \bigoplus \eta \bigoplus \cdots \bigoplus \eta^{n-1}$ est muni d'une
structure de $\fs$-alg\`ebre donn\'ee par un isomorphisme $\tau :\fs \stackrel{\sim}{\rightarrow} \eta^{n}.$ On
consid\`ere le changement de base
\begin{eqnarray}
\shorthandoff{;:!?}
\xymatrix @!0 @R=14mm @C=1.4cm{
       Spec(i^* {\bf \mathcal{A}})= & i^* C  \ar[r]^j \ar[d]_f & C \ar[d]^f &=Spec({\bf \mathcal{A}})\\
       & H \ar[r]^i & H &
     }
\end{eqnarray}
Comme l'involution $i$ agit sur $JH$ par $(-1)_{JH}$, on peut choisir un isomorphisme $\varphi: i^*\eta \stackrel{\sim}
{\rightarrow} \eta^{-1}$ de fa\c con que $\varphi^{\otimes n} = Id$ via $\tau$. On obtient ainsi un isomorphisme de
$\fs$-alg\`ebres ${\bf \mathcal{A}} \rightarrow i^*{\bf \mathcal{A}}= i^*\fs \bigoplus i^*\eta \bigoplus \cdots
\bigoplus i^*\eta^ {n-1} $. De cette fa\c con on peut identifier $i^*C$ \`a $C$ et $j$ \`a un automorphisme
v\'erifiant $j^2=1_{C}$. Observons que ce rel\`evement \`a $C$ de l'involution hyperelliptique n'est pas canonique
car il d\'epend du choix de l'isomorphisme $\varphi$.
Dans ~\cite{BL1}(Prop. 2.1.) on montre la proposition suivante

\begin{Pro}
Le rev\^etement $C \rightarrow \mathbb{P}^1$ est galoisien avec groupe de Galois $Gal(C/\mathbb{P}^1)= D_n =
\langle j,\sigma \mid j^2=\sigma^n =1, \ j\sigma j=\sigma^{-1}\rangle$.
\end{Pro}

Le groupe di\'edral $ D_n=\langle j,\sigma \rangle$ contient les involutions $j_{\nu}=j\sigma^{\nu}$ pour
$\nu=0,\ldots, n-1$. Soient $f_{\nu} :C \rightarrow  C_{\nu}:=C/\langle j_{\nu} \rangle$ les rev\^etements
doubles ramifi\'es associ\'es \`a ces involutions. Soit $g_{\nu}$ le genre de $C_{\nu}$.
On a donc le diagramme suivant

\begin{eqnarray}
\shorthandoff{;:!?}
\xymatrix{
        & C \ar[ld]_{f_0} \ar[d]^f \ar[rd]^{f_1} & \\
	C_0 \ar[rd] & H \ar[d]^{\pi} & C_1 \ar[ld]\\
        &\mathbb{P}^1 &
     }
\end{eqnarray}
Soit $W=\{x_1, \ldots, x_{2g+2} \}\subset \mathbb{P}^1$ l'ensemble des points de Weierstrass pour le rev\^etement
$\pi:H \rightarrow \mathbb{P}^1;$  on pose $S=\{ x\in W \mid \ (\pi \circ f)^{-1}(x) \textrm{ ne contient pas}$
\\ de point fixe par $j \}$ et $T=W \setminus S$.
\begin{Pro}
\hspace{1cm}\\
a) Pour $n$ impair les courbes $C_{\nu}$ sont de genre $g_{\nu}=\frac{1}{2}(n-1)(g-1)$ pour $\nu=0,\dots,n-1$.\\
b) Pour $n$ pair $ g_{\nu} = \frac{n}{2}(g-1)+1-\frac{|T|}{2}$ pour $\nu=0,\dots,n-1$.
\end{Pro}
\begin{Dem}
\hspace{1cm}\\
\emph{a)} Il suffit de montrer la proposition pour $\nu=0$. On observe que les images par $f$ des points fixes par l'
involution $j=j_0$ sont des points de ramification du rev\^etement $H \rightarrow \mathbb{P}^1$ puisque le
diagramme (1) est commutatif. Soit $q \in H$ un point fixe par $i$ et $p \in C$ un rel\`evement de $q$. Comme $f$ est
un rev\^etement non-ramifi\'e, $p$ est un point fixe par une involution $j_m$ de $C$ avec $0 \leq m \leq n$. Puisque $n$
est impair, il existe un unique $k$ modulo $n$ qui v\'erifie l'\'equation $2k \equiv m$ mod $n$. Donc, $\sigma^k p \in
f^{-1}(q)$ est le seul point fixe par $j$ dans la fibre $f^{-1}(q)$. En effet, comme $j\sigma=\sigma^{-1}j$ on a
\begin{eqnarray*}
j\sigma^k p = j\sigma^k j_m p = \sigma^{m-k} p= \sigma^{2k-k}p= \sigma^k p.
\end{eqnarray*}
En conclusion $S= \emptyset $ et on a autant de points fixes par  $j_{\nu}$ que de points de ramification du rev\^etement
$H \rightarrow \mathbb{P}^1$. Par la formule
de Hurwitz on a
\begin{eqnarray*}
g(C)-1=2(g_{\nu}-1) +g+1,
\end{eqnarray*}
et on trouve que $g_{\nu} = \frac{1}{2}(g(C)-g)= \frac{1}{2}(n(g-1)+1-g)= \frac{1}{2}(n-1)(g-1)$.\\
\emph{b)} Dans le cas $n$ pair, sur certaines fibres au-dessus des points de
ramification du rev\^etement $\pi$, l'involution $j$ n'a pas de point fixe.
On observe que si $p\in \F (j)$ alors $\sigma^{\frac{n}{2}}p \in\F (j)$
et qu'ils sont les seuls points fix\'es sur la fibre $f^{-1}(f(p))$. Donc, par la formule de Hurwitz
\begin{equation*}
g(C)-1 = 2g_{\nu}-2 + |T|
\end{equation*}
et on obtient $g_{\nu} =\frac{n}{2}(g-1) + 1- \frac{|T|}{2}$.
\end{Dem}

Les automorphismes $\sigma$ et $j$ induisent des automorphismes dans $JC$, not\'es aussi $\sigma$ et $j$. Les
applications $f_{\nu}^*: JC_{\nu} \rightarrow JC$  sont injectives car les rev\^etements doubles sont ramifi\'es.
On peut donc consid\'erer les jacobiennes $JC_{\nu}$ comme des sous-vari\'et\'es ab\'eliennes de $JC$.
Pour tout point $c \in C$ on obtient le diagramme commutatif
\begin{eqnarray}
\shorthandoff{;:!?}
\xymatrix{
	C \ar@{^{(}->}[r]^{\alpha_c} \ar[d]_{f_\nu} & JC \ar[d]^{\Nmi}  \\
	C_\nu \ar@{^{(}->}[r]^{\alpha_{f_\nu (c)}} &  JC_\nu
     }
\end{eqnarray}
pour $\nu=0,\ldots,n-1$, o\`u $\alpha_c$ est l'application de Abel-Jacobi. On peut \'ecrire $JC_{\nu}= \im (1+
j_{\nu})$ dans $JC$. De plus, comme $\sigma(1+j_{\nu})=(1+j_{\nu -2})\sigma$ l'automorphisme
$\sigma$ se restreint \`a des isomorphismes
\begin{eqnarray*}
\sigma: JC_{\nu} \rightarrow JC_{\nu-2} \qquad \textrm{pour} \quad \nu \in \mathbb{Z}/n\mathbb{Z}.
\end{eqnarray*}
On en d\'eduit que, lorsque $n$ est pair, il y a deux classes d'isomorphismes de jacobiennes qu'on note $JC_0$ et
$JC_1$; lorsque $n$ est impair toutes les jacobiennes $JC_{\nu}$ sont isomorphes.

\begin{Pro}
Les sous-vari\'et\'es $JC_{\nu}$ sont contenues dans la vari\'et\'e de Prym $P$.
\end{Pro}
\begin{Dem}
On a $P=\Ker (1+ \sigma +\cdots + \sigma^{n-1})^0$. D'autre part
\begin{eqnarray*}
(1+\sigma +\cdots +\sigma^{n-1})(1+j_{\nu})=  (1+\sigma +\cdots + \sigma^{n-1}+j + j_1+ \cdots +
j_{n-1})
\end{eqnarray*}
se factorise par l'application Norme du rev\^etement $C \rightarrow \mathbb{P}^1$ et donc est l'application nulle.
En cons\'equence $JC_{\nu} =\im (1+j_{\nu}) \subset \Ker (1+ \sigma +\cdots +\sigma^{n-1})$, et puisque $JC_{\nu}$ est
connexe on a $JC_{\nu} \subset P$.
\end{Dem}
Soient $p_0$ et $p_1$ les projections de $JC_0 \times JC_1$ sur les deux facteurs correspondants.

\begin{The} \label{prym}
L'application $\psi = p_0+p_1 : JC_0 \times JC_1 \rightarrow P$ est un isomorphisme de vari\'et\'es
ab\'eliennes pour $n \geq 2$ non divisible par 4.
\end{The}
\begin{Dem}
Par la prop. 2.3 \emph{a)} dim$P=(n-1)(g-1)=$dim$(JC_0 \times JC_1)$ pour $n$ impair. Lorsque $n$ est pair on
obtient de la proposition 2.3 que dim$JC_0=g_0=\frac{n}{2}(g-1)+1-t$ et dim$JC_1 = \frac{n}{2}(g-1)+1-s$
o\`u $t+s=g+1$. On a donc dim$(JC_0 \times JC_1)=(n-1)(g-1)$. Comme $\psi$ est bien un morphisme de vari\'et\'es
ab\'eliennes il suffit de montrer que $\psi$ est injective.\\
Soit $(x,y) \in JC_0 \times JC_1$ tel que $\psi(x,y)=x+y =0$. Alors $x=(-y) \in JC_0 \cap JC_1$.
Le lemme suivant montre que n\'ecessairement $x=0$ et donc aussi $y=0$, ce qui termine la preuve. Voir~\cite{Mum3}
pag. 346 pour le cas n=2.
\end{Dem}

\begin{Lem}
$ JC_0 \cap JC_1=\{ 0\}$
\end{Lem}
\begin{Dem}
\hspace{1cm}\\
\emph{a) Cas n impair.} Soit $F \in JC_0 \cap JC_1 \subset \Ker(1-j)\cap\Ker(1-j_1)\subset \F(j,\sigma)$; on a donc
$F \in \F (\sigma) \cap P \subset \F (\sigma) \cap \Ker (1+\sigma+\cdots +\sigma^{n-1}) \subset JC[n]$. Par
ailleurs, comme $\sigma^*F \simeq F$ et $f$ \'etant \'etale cela
implique que $F=f^*L$ pour un fibr\'e en droites $L \in JH$. Si de plus $j^*F \simeq F$, alors
\begin{eqnarray*}
j^*f^* L \simeq f^* i^* L \simeq f^* L.
\end{eqnarray*}
Mais dans $JH$ l'involution $i$ agit par $-1_{JH}$, i.e. $i^*L \simeq L^{-1}$. On a donc $f^*L^{-1} \simeq
(f^*L)^{-1} \simeq f^*L$ et $F=f^*L$ est un point de 2-torsion. On conclut que
$F \in \F (j,\sigma )\cap P \subset JC[2] \cap JC[n] = \{0\}$.
\\
\emph{b) Cas n=2m, m impair.} Dans ce cas on utilise le r\'esultat suivant qui est un cas particulier du Lemme
de descente d\^u \`a Kempf (~\cite{D-N}, Th\'eor\`eme 2.3):
\begin{Lem}
Soit $X$ une vari\'et\'e alg\'ebrique int\`egre sur laquelle op\`ere un groupe fini
$G$. Soit $F$ un $G$-fibr\'e vectoriel sur $X$. Alors $F$ descend \`a $X/G$ si et seulement si pour tout point
$x \in X$, le stabilisateur de $x$ dans $G$ agit trivialement sur $F_x$.
\end{Lem}
Soient $X_i:= C_i / \langle \sigma^m \rangle$ et $ X:= C/ \langle \sigma^m \rangle$, o\`u $\sigma^m$ est une
involution qui commute avec $j$ et $j_1$. On consid\`ere la tour de courbes suivante
\begin{eqnarray}
\shorthandoff{;:!?}
\def\commutatif{\ar@{}[r]|{\circlearrowleft}}
\xymatrix{
	&C \ar[ld]_{q_0} \ar[d]^q \ar[rd]^{q_1}  & \\
	C_0 \ar[d]_{r_0} \commutatif & X \ar[ld]_{f_0} \ar[d]^f \ar[rd]^{f_1} \commutatif &C_1 \ar[d]^{r_1} \\
        X_0 \ar[rd] & H \ar[d]^{\pi} &X_1 \ar[ld] \\
        & \mathbb{P}^1 &
     }
\end{eqnarray}
\\
Soit $F \in JC_0 \cap JC_1$. On sait qu'il existe des fibr\'es en droites $M_i \in JC_i, i=0,1$ tels
que $q_i^*M_i\simeq F$. Puisque $j$ et $\sigma^m$ commutent $M_0$ est invariante par $\sigma^m$. En effet,
\begin{eqnarray*}
q_0^*\sigma^{m*}M_0 \simeq \sigma^{m*}q_0^*M_0 \simeq \sigma^{m*}F \simeq F \simeq q_0^*M_0.
\end{eqnarray*}
Comme $q_0$ est ramifi\'e, $q_0^*$ est injective et donc $\sigma^{m*}M_0 \simeq M_0$. Observons que les points
de ramification du rev\^etement $r_0: C_0 \rightarrow X_0$, i.e. les points fixes par $\sigma^m$, peuvent \^etre
relev\'es aux points fixes par $j\sigma^m=j_m$ dans $C$. En effet, soit $p \in \F (\sigma^m)$ dans $C_0$ et soit
$\tilde{p} \in C$ tel que $q_0(\tilde{p})=p$. On a
\begin{eqnarray*}
q_0(\sigma^m \tilde{p})= \sigma^m q_0(\tilde{p})= \sigma^m p =p,
\end{eqnarray*}
donc $\sigma^m \tilde{p} \in q_0^{-1}(p)= \{ \tilde{p}, j\tilde{p} \}$. Comme $q: C \rightarrow X $ est non-ramifi\'e,
$\sigma^m\tilde{p} \neq \tilde{p}$. Ainsi $\sigma^m\tilde{p} = j\tilde{p}$ et on a $j\sigma^m \tilde{p}=j_m\tilde{p}=
\tilde{p}$.\\
L'action de $\sigma^m$ sur les fibres de $M_0$ au-dessus des points de
ramification de $r_0$ est la m\^eme que celle de $j_m$ sur les fibres de $F$ au-dessus des points fixes par
$j_m$ dans $C$ puisque $q_0^*M_0\simeq F$.\\
Soit $x\in \F(j_m) \subset C$, donc $x=\sigma^{\frac{m+1}{2}}y$ ou bien $x=\sigma^{\frac{3m+1}{2}}y$ o\`u $y \in
\F(j_1)$. On observe que $\langle j_m \rangle =Stab(x) $ est un sous-groupe conjugu\'e de
$\langle j_1 \rangle$ qui par hypoth\`ese agit trivialement sur $F_y$, donc $j_m$ agit aussi
trivialement sur $F_x$. On en d\'eduit que $\sigma^m$ agit trivialement sur $M_{0,q_0(x)}$ et par le lemme de
descente il existe un fibr\'e $N_0 \in JX_0$ tel que $r_0^*N_0 \simeq
M_0$. De fa\c con analogue on montre l'existence d'un fibr\'e $N_1 \in JX_1$  tel que $r_1^*N_1 \simeq M_1$.\\
Comme $q^*f_0^*N_0 \simeq q^*f_1^*N_1 \simeq F$ on a
\begin{eqnarray*}
\beta:=f_0^*N_0 \otimes (f_1^*N_1)^{-1} \in \Ker q^*.
\end{eqnarray*}
Puisque $q$ est un rev\^etement double non-ramifi\'e, $\beta^2 \simeq \mathcal{O}_X$, i.e., $\beta \in JX[2]$ mais
aussi $\beta \in JX_0 \times JX_1 \simeq  Prym(X/H) $, ce dernier isomorphisme ayant \'et\'e
d\'emontr\'e dans le cas \emph{a)}.
Comme
\begin{eqnarray*}
q^*\beta \simeq \mathcal{O}_C \simeq \sigma^*\mathcal{O}_C \simeq \sigma^*q^*\beta
\simeq q^*\sigma^*\beta
\end{eqnarray*}
on a $\sigma^*\beta \in \Ker q^*$. En fait $\sigma^*\beta \simeq \beta$ et puisque $Prym(X/H) \subset \Ker (1+\sigma +\cdots
+ \sigma^{m-1})$, on a $\beta \in JX[m]\cap JX[2]=\{0\}$. Ainsi, $f_0^*N_0 \simeq f_1^*N_1 \in JX_0 \cap JX_1 =
\{ 0 \}$  par \emph{a)} et donc $F\simeq q^* \mathcal{O}_X \simeq \mathcal{O}_C$.
\end{Dem}

\section{La polarisation} \label{polar_prym}

On consid\`ere la polarisation $\Xi$ dans $JC_0 \times JC_1$. On a $\Xi \equiv \psi^*\Theta$, o\`u $\Theta$ est la
polarisation principale dans $JC$. On pose $\phi =\phi_{\Xi} : JC_0 \times JC_1 \longrightarrow \widehat{JC_0}
\times \widehat{JC_1}$. Donc $\phi$ est de la forme
\begin{displaymath}
\phi = \left( \begin{array}{cc} \alpha & \hat{\beta}\\ \beta & \delta \end{array} \right)
\end{displaymath}
L'application $\alpha: JC_0 \rightarrow \widehat{JC_0}$ est la restriction \`a $JC_0$ de la polarisation principale
de $JC$. Or, l'inclusion $f_0^*: JC_0 \rightarrow JC$ est le pullback d'un rev\^etement double ramifi\'e, et
par ~\cite{CAV} (12.3.1.) on obtient $(f_0^*)^*\Theta \equiv 2\Theta_0$, o\`u $\Theta_0$ est la polarisation
principale dans $JC_0$. Ainsi $\alpha= \phi_{2\Theta_0}=2\phi_{\Theta_0}$. De fa\c con analogue, on obtient
$\delta = 2\phi_{\Theta_1}: JC_1 \rightarrow \widehat{JC_1}$, o\`u $\Theta_1$ est la polarisation principale dans
$JC_1$.\\ L'application $\beta$ est l'application qui fait commuter le diagramme
\begin{eqnarray}
\shorthandoff{;:!?}
\xymatrix{
	JC \ar[r]^{\phi_\Theta}  & \widehat{JC} \ar[d]^{\widehat{f_1^*}} \\
	JC_0 \ar@{^{(}->}[u]^{f_0^*} \ar[r]^{\beta} &  \widehat{JC_1}
     }
\end{eqnarray}
donc $\beta = \widehat{f_1^*} \circ\phi_\Theta \circ f_0^*$ et $\hat{\beta} = \widehat{f_0^*} \circ\phi_\Theta
\circ f_1^*$. Pour expliciter $\widehat{f_i^*}$ on consid\`ere le diagramme commutatif
\begin{eqnarray}
\shorthandoff{;:!?}
\xymatrix{
	JC \ar[r]^{\phi_\Theta}  & \widehat{JC} \\
	JC_i \ar@{^{(}->}[u]^{f_i^*} \ar[r]^{\phi_{\Theta_i}} &  \widehat{JC_i} \ar[u]_{\widehat{\Nmi}}
     }
\end{eqnarray}
pour $i=0,1$, obtenu en appliquant le foncteur $Pic^0$ au diagramme (3). Ensuite, en dualisant (6) on a
\begin{eqnarray}
\shorthandoff{;:!?}
\xymatrix{
	\widehat{JC} \ar[r]^{\phi_\Theta^{-1}} \ar[d]_{\widehat{f_i^*}} & JC \ar[d]^{\Nmi}\\
	\widehat{JC_i}  \ar[r]^{\phi_{\Theta_i}^{-1}} &  JC_i
     }
\end{eqnarray}
pour $i=0,1$. On a donc $\widehat{f_i^*} = \phi_{\Theta_i} \circ \Nmi \circ \phi_{\Theta}^{-1}$ et en
utilisant le fait que $\Nmi = 1+j_i$ on obtient
\begin{displaymath}
\beta = \phi_{\Theta_1} \circ (1+j_1) \circ f^*_0 \textrm{ \ \ et \ \ } \hat{\beta} = \phi_{\Theta_0}
\circ (1+j) \circ f^*_1.
\end{displaymath}

\backmatter
%% Created by Maple 8.00 (IBM INTEL LINUX)
%% Source Worksheet: rango3.mws
%% Generated: Mon Jun 23 17:24:02 2003

%\documentclass{article}
%\usepackage{maple2e}
 %\def\emptyline{\vspace{12pt}}
%\DefineParaStyle{Maple Output}
%\DefineParaStyle{Warning}
%\DefineCharStyle{2D Math}
%\DefineCharStyle{2D Output}
%\begin{document}
%\pagestyle{empty}
\chapter{Annexe}
\addcontentsline{chapter}{toc}{Annexe}

\begin{maplegroup}

Variables:
\emptyline
\end{maplegroup}
\begin{maplegroup}
\begin{mapleinput}
\mapleinline{active}{1d}{Z:= array(0..2,0..2);}{%
}
\end{mapleinput}

\mapleresult
\begin{maplelatex}
\mapleinline{inert}{2d}{Z := array(0 .. 2,0 .. 2,[]);}{%
\[
Z := \mathrm{array}(0 .. 2, \,0 .. 2, \,[])
\]
}
\end{maplelatex}
\end{maplegroup}
\begin{maplegroup}
\emptyline
Tableau de polyn\^omes invariants par le groupe de Heisenberg:
\end{maplegroup}
\emptyline
\begin{maplegroup}
\begin{mapleinput}
\mapleinline{active}{1d}{T:= array(1..43);
                       }{%
}
\end{mapleinput}
\begin{mapleinput}
\end{mapleinput}
\mapleresult
\begin{maplelatex}
\mapleinline{inert}{2d}{T := array(1 .. 43,[]);}{%
\[
T := \mathrm{array}(1 .. 43, \,[])
\]
}
\end{maplelatex}
\end{maplegroup}
\emptyline
\begin{maplegroup}

Proc\'edure,  l'action du sous-groupe isotrope  sur les  indices de la
base:
\end{maplegroup}
\begin{maplegroup}
\begin{mapleinput}
\mapleinline{active}{1d}{q:=proc(n,t) RETURN( (n+t) mod 3 ) end;}{%
}
\end{mapleinput}

\mapleresult
\begin{maplelatex}
\mapleinline{inert}{2d}{q := proc (n, t) RETURN(`mod`(n+t,3)) end proc;}{%
\[
q := \textbf{proc} (n, \,t)\,\mathrm{RETURN}((n + t)\,\mathrm{mod
}\,3)\,\textbf{end proc} 
\]
}
\end{maplelatex}
\end{maplegroup}
\emptyline
\begin{maplegroup}

L'action du sous-groupe isotrope sur les mon\^omes:
\end{maplegroup}
\begin{maplegroup}
\begin{mapleinput}
\mapleinline{active}{1d}{Act:=proc(r,s,P)   RETURN( subs( \{seq( seq( Z[i,j]=Z[q(r,i),q(s,j)]
,i=0..2),j=0..2 )\} ,P )  ) end;}{%
}
\end{mapleinput}

\mapleresult
\begin{maplelatex}
\mapleinline{inert}{2d}{Act := proc (r, s, P) 
RETURN(subs(\{seq(seq(Z[i,j] =
Z[q(r,i),q(s,j)],i = 0 .. 2),j = 0 .. 2)\},P)) end proc;}{%
\maplemultiline{
\mathit{Act} := \textbf{proc} (r, \,s, \,P) \\
\mapleIndent{1} \mathrm{RETURN}(\mathrm{subs}(\{\mathrm{seq}(
\mathrm{seq}({Z_{i, \,j}}={Z_{\mathrm{q}(r, \,i), \,\mathrm{q}(s
, \,j)}}, \,i=0 .. 2), \,j=0 .. 2)\}, \,P)) \\
\textbf{end proc}  }
}
\end{maplelatex}

\end{maplegroup}
\begin{maplegroup}
\emptyline
Proc\'edure pour engendrer  un  polyn\^ome invariant:
\end{maplegroup}
\emptyline
\begin{maplegroup}
\begin{mapleinput}
\mapleinline{active}{1d}{Gen:=proc(P) RETURN( sum( 'sum('Act(l,k,P)','l'=0..2)',k=0..2))
end;}{%
}
\end{mapleinput}

\mapleresult
\begin{maplelatex}
\mapleinline{inert}{2d}{Gen := proc (P) RETURN(sum('sum('Act(l,k,P)',('l') = 0 .. 2)',k = 0
.. 2)) end proc;}{%
\maplemultiline{
\mathit{Gen} := \textbf{proc} (P)\, \\
\mapleIndent{1}\mathrm{RETURN}(\mathrm{sum}(
\hbox{'}\mathrm{sum}(\hbox{'}\mathrm{Act}(l, \,k, \,P)\hbox{'}, 
\,\hbox{'}l\hbox{'}=0 .. 2)\hbox{'}, \,k=0 .. 2))\,\textbf{end proc} }
}
\end{maplelatex}
\end{maplegroup}
\begin{maplegroup}
\emptyline
Remplissage du tableau:
\emptyline
\end{maplegroup}
\begin{maplegroup}
\begin{mapleinput}
\mapleinline{active}{1d}{T[1]:=Gen(Z[0,0]^6);}{%
}
\end{mapleinput}

\mapleresult
\begin{maplelatex}
\mapleinline{inert}{2d}{T[1] :=
Z[0,0]^6+Z[1,0]^6+Z[2,0]^6+Z[0,1]^6+Z[1,1]^6+Z[2,1]^6+Z[0,2]^6+Z[1,2]^
6+Z[2,2]^6;}{%
\maplemultiline{
{T_{1}} := {Z_{0, \,0}}^{6} + {Z_{1, \,0}}^{6} + {Z_{2, \,0}}^{6}
 + {Z_{0, \,1}}^{6} + {Z_{1, \,1}}^{6} + {Z_{2, \,1}}^{6} + {Z_{0
, \,2}}^{6} \\ 
+ {Z_{1, \,2}}^{6} + {Z_{2, \,2}}^{6} }
}
\end{maplelatex}
\end{maplegroup}
\emptyline
\begin{maplegroup}
\begin{mapleinput}
\mapleinline{active}{1d}{T[2]:=Gen(Z[0,0]^3*Z[0,1]^3);}{%
}
\end{mapleinput}

\mapleresult
\begin{maplelatex}
\mapleinline{inert}{2d}{T[2] :=
Z[0,0]^3*Z[0,1]^3+Z[1,0]^3*Z[1,1]^3+Z[2,0]^3*Z[2,1]^3+Z[0,1]^3*Z[0,2]^
3+Z[1,1]^3*Z[1,2]^3+Z[2,1]^3*Z[2,2]^3+Z[0,2]^3*Z[0,0]^3+Z[1,2]^3*Z[1,0
]^3+Z[2,2]^3*Z[2,0]^3;}{%
\maplemultiline{
{T_{2}} := {Z_{0, \,0}}^{3}\,{Z_{0, \,1}}^{3} + {Z_{1, \,0}}^{3}
\,{Z_{1, \,1}}^{3} + {Z_{2, \,0}}^{3}\,{Z_{2, \,1}}^{3} + {Z_{0, 
\,1}}^{3}\,{Z_{0, \,2}}^{3}\\ + {Z_{1, \,1}}^{3}\,{Z_{1, \,2}}^{3}
 + {Z_{2, \,1}}^{3}\,{Z_{2, \,2}}^{3} 
 + {Z_{0, \,2}}^{3}\,{Z_{0, \,0}}^{3} + {Z_{1, \,2}}^{3}\,
{Z_{1, \,0}}^{3} + {Z_{2, \,2}}^{3}\,{Z_{2, \,0}}^{3} }
}
\end{maplelatex}
\end{maplegroup}
\emptyline
\begin{maplegroup}
\begin{mapleinput}
\mapleinline{active}{1d}{T[3]:=Gen(Z[0,0]^3*Z[1,0]^3);}{%
}
\end{mapleinput}

\mapleresult
\begin{maplelatex}
\mapleinline{inert}{2d}{T[3] :=
Z[0,0]^3*Z[1,0]^3+Z[1,0]^3*Z[2,0]^3+Z[2,0]^3*Z[0,0]^3+Z[0,1]^3*Z[1,1]^
3+Z[1,1]^3*Z[2,1]^3+Z[2,1]^3*Z[0,1]^3+Z[0,2]^3*Z[1,2]^3+Z[1,2]^3*Z[2,2
]^3+Z[2,2]^3*Z[0,2]^3;}{%
\maplemultiline{
{T_{3}} := {Z_{0, \,0}}^{3}\,{Z_{1, \,0}}^{3} + {Z_{1, \,0}}^{3}
\,{Z_{2, \,0}}^{3} + {Z_{2, \,0}}^{3}\,{Z_{0, \,0}}^{3} + {Z_{0, 
\,1}}^{3}\,{Z_{1, \,1}}^{3} \\ + {Z_{1, \,1}}^{3}\,{Z_{2, \,1}}^{3}
 + {Z_{2, \,1}}^{3}\,{Z_{0, \,1}}^{3} 
\mbox{} + {Z_{0, \,2}}^{3}\,{Z_{1, \,2}}^{3} + {Z_{1, \,2}}^{3}\,
{Z_{2, \,2}}^{3} + {Z_{2, \,2}}^{3}\,{Z_{0, \,2}}^{3} }
}
\end{maplelatex}
\end{maplegroup}
\emptyline
\begin{maplegroup}
\begin{mapleinput}
\mapleinline{active}{1d}{T[4]:=Gen(Z[0,0]^3*Z[1,1]^3); }{%
}
\end{mapleinput}
\mapleresult
\begin{maplelatex}
\mapleinline{inert}{2d}{T[4] :=
Z[0,0]^3*Z[1,1]^3+Z[1,0]^3*Z[2,1]^3+Z[2,0]^3*Z[0,1]^3+Z[0,1]^3*Z[1,2]^
3+Z[1,1]^3*Z[2,2]^3+Z[2,1]^3*Z[0,2]^3+Z[0,2]^3*Z[1,0]^3+Z[1,2]^3*Z[2,0
]^3+Z[2,2]^3*Z[0,0]^3;}{%
\maplemultiline{
{T_{4}} := {Z_{0, \,0}}^{3}\,{Z_{1, \,1}}^{3} + {Z_{1, \,0}}^{3}
\,{Z_{2, \,1}}^{3} + {Z_{2, \,0}}^{3}\,{Z_{0, \,1}}^{3} + {Z_{0, 
\,1}}^{3}\,{Z_{1, \,2}}^{3}\\ + {Z_{1, \,1}}^{3}\,{Z_{2, \,2}}^{3}
 + {Z_{2, \,1}}^{3}\,{Z_{0, \,2}}^{3} 
\mbox{} + {Z_{0, \,2}}^{3}\,{Z_{1, \,0}}^{3} + {Z_{1, \,2}}^{3}\,
{Z_{2, \,0}}^{3} + {Z_{2, \,2}}^{3}\,{Z_{0, \,0}}^{3} }
}
\end{maplelatex}
\end{maplegroup}
\emptyline
\begin{maplegroup}
\begin{mapleinput}
\mapleinline{active}{1d}{T[5]:=Gen(Z[0,0]^3*Z[1,2]^3); }{%
}
\end{mapleinput}
\mapleresult
\begin{maplelatex}
\mapleinline{inert}{2d}{T[5] :=
Z[0,0]^3*Z[1,2]^3+Z[1,0]^3*Z[2,2]^3+Z[2,0]^3*Z[0,2]^3+Z[0,1]^3*Z[1,0]^
3+Z[1,1]^3*Z[2,0]^3+Z[2,1]^3*Z[0,0]^3+Z[0,2]^3*Z[1,1]^3+Z[1,2]^3*Z[2,1
]^3+Z[2,2]^3*Z[0,1]^3;}{%
\maplemultiline{
{T_{5}} := {Z_{0, \,0}}^{3}\,{Z_{1, \,2}}^{3} + {Z_{1, \,0}}^{3}
\,{Z_{2, \,2}}^{3} + {Z_{2, \,0}}^{3}\,{Z_{0, \,2}}^{3} + {Z_{0, 
\,1}}^{3}\,{Z_{1, \,0}}^{3} \\+ {Z_{1, \,1}}^{3}\,{Z_{2, \,0}}^{3}
 + {Z_{2, \,1}}^{3}\,{Z_{0, \,0}}^{3} 
\mbox{} + {Z_{0, \,2}}^{3}\,{Z_{1, \,1}}^{3} + {Z_{1, \,2}}^{3}\,
{Z_{2, \,1}}^{3} + {Z_{2, \,2}}^{3}\,{Z_{0, \,1}}^{3} }
}
\end{maplelatex}
\end{maplegroup}

\emptyline
\begin{maplegroup}
\begin{mapleinput}
\mapleinline{active}{1d}{T[6]:=Gen(Z[0,0]^4*Z[0,1]*Z[0,2]);}{%
}
\end{mapleinput}

\mapleresult
\begin{maplelatex}
\mapleinline{inert}{2d}{T[6] :=
Z[0,0]^4*Z[0,1]*Z[0,2]+Z[1,0]^4*Z[1,1]*Z[1,2]+Z[2,0]^4*Z[2,1]*Z[2,2]+Z
[0,1]^4*Z[0,2]*Z[0,0]+Z[1,1]^4*Z[1,2]*Z[1,0]+Z[2,1]^4*Z[2,2]*Z[2,0]+Z[
0,2]^4*Z[0,0]*Z[0,1]+Z[1,2]^4*Z[1,0]*Z[1,1]+Z[2,2]^4*Z[2,0]*Z[2,1];}{%
\maplemultiline{
{T_{6}} := {Z_{0, \,0}}^{4}\,{Z_{0, \,1}}\,{Z_{0, \,2}} + {Z_{1, 
\,0}}^{4}\,{Z_{1, \,1}}\,{Z_{1, \,2}} + {Z_{2, \,0}}^{4}\,{Z_{2, 
\,1}}\,{Z_{2, \,2}}\\\mbox{}  + {Z_{0, \,1}}^{4}\,{Z_{0, \,2}}\,{Z_{0, \,0}
} + {Z_{1, \,1}}^{4}\,{Z_{1, \,2}}\,{Z_{1, \,0}}
 + {Z_{2, \,1}}^{4}\,{Z_{2, \,2}}\,{Z_{2, \,0}}\\ \mbox{} + {Z_{0, 
\,2}}^{4}\,{Z_{0, \,0}}\,{Z_{0, \,1}} + {Z_{1, \,2}}^{4}\,{Z_{1, 
\,0}}\,{Z_{1, \,1}} + {Z_{2, \,2}}^{4}\,{Z_{2, \,0}}\,{Z_{2, \,1}
} }
}
\end{maplelatex}
\end{maplegroup}

\emptyline
\begin{maplegroup}
\begin{mapleinput}
\mapleinline{active}{1d}{T[7]:=Gen(Z[0,0]*Z[0,1]*Z[0,2]*Z[1,0]^3);}{%
}
\end{mapleinput}
\mapleresult
\begin{maplelatex}
\mapleinline{inert}{2d}{T[7] :=
Z[0,0]*Z[0,1]*Z[0,2]*Z[1,0]^3+Z[1,0]*Z[1,1]*Z[1,2]*Z[2,0]^3+Z[2,0]*Z[2
,1]*Z[2,2]*Z[0,0]^3+Z[0,1]*Z[0,2]*Z[0,0]*Z[1,1]^3+Z[1,1]*Z[1,2]*Z[1,0]
*Z[2,1]^3+Z[2,1]*Z[2,2]*Z[2,0]*Z[0,1]^3+Z[0,2]*Z[0,0]*Z[0,1]*Z[1,2]^3+
Z[1,2]*Z[1,0]*Z[1,1]*Z[2,2]^3+Z[2,2]*Z[2,0]*Z[2,1]*Z[0,2]^3;}{%
\maplemultiline{
{T_{7}} := {Z_{0, \,0}}\,{Z_{0, \,1}}\,{Z_{0, \,2}}\,{Z_{1, \,0}}
^{3} + {Z_{1, \,0}}\,{Z_{1, \,1}}\,{Z_{1, \,2}}\,{Z_{2, \,0}}^{3}
 + {Z_{2, \,0}}\,{Z_{2, \,1}}\,{Z_{2, \,2}}\,{Z_{0, \,0}}^{3}\\ + {
Z_{0, \,1}}\,{Z_{0, \,2}}\,{Z_{0, \,0}}\,{Z_{1, \,1}}^{3} 
\mbox{} + {Z_{1, \,1}}\,{Z_{1, \,2}}\,{Z_{1, \,0}}\,{Z_{2, \,1}}
^{3} + {Z_{2, \,1}}\,{Z_{2, \,2}}\,{Z_{2, \,0}}\,{Z_{0, \,1}}^{3}\\
 + {Z_{0, \,2}}\,{Z_{0, \,0}}\,{Z_{0, \,1}}\,{Z_{1, \,2}}^{3} 
\mbox{} + {Z_{1, \,2}}\,{Z_{1, \,0}}\,{Z_{1, \,1}}\,{Z_{2, \,2}}
^{3} + {Z_{2, \,2}}\,{Z_{2, \,0}}\,{Z_{2, \,1}}\,{Z_{0, \,2}}^{3}
 }
}
\end{maplelatex}
\end{maplegroup}

\emptyline
\begin{maplegroup}
\begin{mapleinput}
\mapleinline{active}{1d}{T[8]:=Gen(Z[0,0]*Z[0,1]*Z[0,2]*Z[2,0]^3);}{%
}
\end{mapleinput}
\mapleresult
\begin{maplelatex}
\mapleinline{inert}{2d}{T[8] :=
Z[0,0]*Z[0,1]*Z[0,2]*Z[2,0]^3+Z[1,0]*Z[1,1]*Z[1,2]*Z[0,0]^3+Z[2,0]*Z[2
,1]*Z[2,2]*Z[1,0]^3+Z[0,1]*Z[0,2]*Z[0,0]*Z[2,1]^3+Z[1,1]*Z[1,2]*Z[1,0]
*Z[0,1]^3+Z[2,1]*Z[2,2]*Z[2,0]*Z[1,1]^3+Z[0,2]*Z[0,0]*Z[0,1]*Z[2,2]^3+
Z[1,2]*Z[1,0]*Z[1,1]*Z[0,2]^3+Z[2,2]*Z[2,0]*Z[2,1]*Z[1,2]^3;}{%
\maplemultiline{
{T_{8}} := {Z_{0, \,0}}\,{Z_{0, \,1}}\,{Z_{0, \,2}}\,{Z_{2, \,0}}
^{3} + {Z_{1, \,0}}\,{Z_{1, \,1}}\,{Z_{1, \,2}}\,{Z_{0, \,0}}^{3}
 + {Z_{2, \,0}}\,{Z_{2, \,1}}\,{Z_{2, \,2}}\,{Z_{1, \,0}}^{3}\\ + {
Z_{0, \,1}}\,{Z_{0, \,2}}\,{Z_{0, \,0}}\,{Z_{2, \,1}}^{3} 
 + {Z_{1, \,1}}\,{Z_{1, \,2}}\,{Z_{1, \,0}}\,{Z_{0, \,1}}
^{3} + {Z_{2, \,1}}\,{Z_{2, \,2}}\,{Z_{2, \,0}}\,{Z_{1, \,1}}^{3}\\
 + {Z_{0, \,2}}\,{Z_{0, \,0}}\,{Z_{0, \,1}}\,{Z_{2, \,2}}^{3} 
 + {Z_{1, \,2}}\,{Z_{1, \,0}}\,{Z_{1, \,1}}\,{Z_{0, \,2}}
^{3} + {Z_{2, \,2}}\,{Z_{2, \,0}}\,{Z_{2, \,1}}\,{Z_{1, \,2}}^{3}
 }
}
\end{maplelatex}
\end{maplegroup}
\emptyline
\begin{maplegroup}
\begin{mapleinput}
\mapleinline{active}{1d}{T[9]:=Gen(Z[0,0]^4*Z[1,0]*Z[2,0]);}{%
}
\end{mapleinput}
\mapleresult
\begin{maplelatex}
\mapleinline{inert}{2d}{T[9] :=
Z[0,0]^4*Z[1,0]*Z[2,0]+Z[1,0]^4*Z[2,0]*Z[0,0]+Z[2,0]^4*Z[0,0]*Z[1,0]+Z
[0,1]^4*Z[1,1]*Z[2,1]+Z[1,1]^4*Z[2,1]*Z[0,1]+Z[2,1]^4*Z[0,1]*Z[1,1]+Z[
0,2]^4*Z[1,2]*Z[2,2]+Z[1,2]^4*Z[2,2]*Z[0,2]+Z[2,2]^4*Z[0,2]*Z[1,2];}{%
\maplemultiline{
{T_{9}} := {Z_{0, \,0}}^{4}\,{Z_{1, \,0}}\,{Z_{2, \,0}} + {Z_{1, 
\,0}}^{4}\,{Z_{2, \,0}}\,{Z_{0, \,0}} + {Z_{2, \,0}}^{4}\,{Z_{0, 
\,0}}\,{Z_{1, \,0}}\\ + {Z_{0, \,1}}^{4}\,{Z_{1, \,1}}\,{Z_{2, \,1}
} + {Z_{1, \,1}}^{4}\,{Z_{2, \,1}}\,{Z_{0, \,1}} 
 + {Z_{2, \,1}}^{4}\,{Z_{0, \,1}}\,{Z_{1, \,1}}\\ + {Z_{0, 
\,2}}^{4}\,{Z_{1, \,2}}\,{Z_{2, \,2}} + {Z_{1, \,2}}^{4}\,{Z_{2, 
\,2}}\,{Z_{0, \,2}} + {Z_{2, \,2}}^{4}\,{Z_{0, \,2}}\,{Z_{1, \,2}
} }
}
\end{maplelatex}
\end{maplegroup}
\emptyline
\begin{maplegroup}
\begin{mapleinput}
\mapleinline{active}{1d}{T[10]:=Gen(Z[0,0]*Z[1,0]*Z[2,0]*Z[0,1]^3);}{%
}
\end{mapleinput}
\mapleresult
\begin{maplelatex}
\mapleinline{inert}{2d}{T[10] :=
Z[0,0]*Z[1,0]*Z[2,0]*Z[0,1]^3+Z[1,0]*Z[2,0]*Z[0,0]*Z[1,1]^3+Z[2,0]*Z[0
,0]*Z[1,0]*Z[2,1]^3+Z[0,1]*Z[1,1]*Z[2,1]*Z[0,2]^3+Z[1,1]*Z[2,1]*Z[0,1]
*Z[1,2]^3+Z[2,1]*Z[0,1]*Z[1,1]*Z[2,2]^3+Z[0,2]*Z[1,2]*Z[2,2]*Z[0,0]^3+
Z[1,2]*Z[2,2]*Z[0,2]*Z[1,0]^3+Z[2,2]*Z[0,2]*Z[1,2]*Z[2,0]^3;}{%
\maplemultiline{
{T_{10}} := {Z_{0, \,0}}\,{Z_{1, \,0}}\,{Z_{2, \,0}}\,{Z_{0, \,1}
}^{3} + {Z_{1, \,0}}\,{Z_{2, \,0}}\,{Z_{0, \,0}}\,{Z_{1, \,1}}^{3
} + {Z_{2, \,0}}\,{Z_{0, \,0}}\,{Z_{1, \,0}}\,{Z_{2, \,1}}^{3}
 \\
\mbox{} + {Z_{0, \,1}}\,{Z_{1, \,1}}\,{Z_{2, \,1}}\,{Z_{0, \,2}}
^{3} + {Z_{1, \,1}}\,{Z_{2, \,1}}\,{Z_{0, \,1}}\,{Z_{1, \,2}}^{3}
 + {Z_{2, \,1}}\,{Z_{0, \,1}}\,{Z_{1, \,1}}\,{Z_{2, \,2}}^{3} \\
\mbox{} + {Z_{0, \,2}}\,{Z_{1, \,2}}\,{Z_{2, \,2}}\,{Z_{0, \,0}}
^{3} + {Z_{1, \,2}}\,{Z_{2, \,2}}\,{Z_{0, \,2}}\,{Z_{1, \,0}}^{3}
 + {Z_{2, \,2}}\,{Z_{0, \,2}}\,{Z_{1, \,2}}\,{Z_{2, \,0}}^{3} }
}
\end{maplelatex}
\end{maplegroup}
\emptyline
\begin{maplegroup}
\begin{mapleinput}
\mapleinline{active}{1d}{T[11]:=Gen(Z[0,0]*Z[1,0]*Z[2,0]*Z[0,2]^3);}{%
}
\end{mapleinput}
\mapleresult
\begin{maplelatex}
\mapleinline{inert}{2d}{T[11] :=
Z[0,0]*Z[1,0]*Z[2,0]*Z[0,2]^3+Z[1,0]*Z[2,0]*Z[0,0]*Z[1,2]^3+Z[2,0]*Z[0
,0]*Z[1,0]*Z[2,2]^3+Z[0,1]*Z[1,1]*Z[2,1]*Z[0,0]^3+Z[1,1]*Z[2,1]*Z[0,1]
*Z[1,0]^3+Z[2,1]*Z[0,1]*Z[1,1]*Z[2,0]^3+Z[0,2]*Z[1,2]*Z[2,2]*Z[0,1]^3+
Z[1,2]*Z[2,2]*Z[0,2]*Z[1,1]^3+Z[2,2]*Z[0,2]*Z[1,2]*Z[2,1]^3;}{%
\maplemultiline{
{T_{11}} := {Z_{0, \,0}}\,{Z_{1, \,0}}\,{Z_{2, \,0}}\,{Z_{0, \,2}
}^{3} + {Z_{1, \,0}}\,{Z_{2, \,0}}\,{Z_{0, \,0}}\,{Z_{1, \,2}}^{3
} + {Z_{2, \,0}}\,{Z_{0, \,0}}\,{Z_{1, \,0}}\,{Z_{2, \,2}}^{3}
 \\
\mbox{} + {Z_{0, \,1}}\,{Z_{1, \,1}}\,{Z_{2, \,1}}\,{Z_{0, \,0}}
^{3} + {Z_{1, \,1}}\,{Z_{2, \,1}}\,{Z_{0, \,1}}\,{Z_{1, \,0}}^{3}
 + {Z_{2, \,1}}\,{Z_{0, \,1}}\,{Z_{1, \,1}}\,{Z_{2, \,0}}^{3} \\
\mbox{} + {Z_{0, \,2}}\,{Z_{1, \,2}}\,{Z_{2, \,2}}\,{Z_{0, \,1}}
^{3} + {Z_{1, \,2}}\,{Z_{2, \,2}}\,{Z_{0, \,2}}\,{Z_{1, \,1}}^{3}
 + {Z_{2, \,2}}\,{Z_{0, \,2}}\,{Z_{1, \,2}}\,{Z_{2, \,1}}^{3} }
}
\end{maplelatex}
\end{maplegroup}
\emptyline
\begin{maplegroup}
\begin{mapleinput}
\mapleinline{active}{1d}{T[12]:=Gen(Z[0,0]^4*Z[1,1]*Z[2,2]);}{%
}
\end{mapleinput}
\mapleresult
\begin{maplelatex}
\mapleinline{inert}{2d}{T[12] :=
Z[0,0]^4*Z[1,1]*Z[2,2]+Z[1,0]^4*Z[2,1]*Z[0,2]+Z[2,0]^4*Z[0,1]*Z[1,2]+Z
[0,1]^4*Z[1,2]*Z[2,0]+Z[1,1]^4*Z[2,2]*Z[0,0]+Z[2,1]^4*Z[0,2]*Z[1,0]+Z[
0,2]^4*Z[1,0]*Z[2,1]+Z[1,2]^4*Z[2,0]*Z[0,1]+Z[2,2]^4*Z[0,0]*Z[1,1];}{%
\maplemultiline{
{T_{12}} := {Z_{0, \,0}}^{4}\,{Z_{1, \,1}}\,{Z_{2, \,2}} + {Z_{1
, \,0}}^{4}\,{Z_{2, \,1}}\,{Z_{0, \,2}} + {Z_{2, \,0}}^{4}\,{Z_{0
, \,1}}\,{Z_{1, \,2}}\\ + {Z_{0, \,1}}^{4}\,{Z_{1, \,2}}\,{Z_{2, \,
0}} + {Z_{1, \,1}}^{4}\,{Z_{2, \,2}}\,{Z_{0, \,0}} 
+ {Z_{2, \,1}}^{4}\,{Z_{0, \,2}}\,{Z_{1, \,0}}\\ + {Z_{0, 
\,2}}^{4}\,{Z_{1, \,0}}\,{Z_{2, \,1}} + {Z_{1, \,2}}^{4}\,{Z_{2, 
\,0}}\,{Z_{0, \,1}} + {Z_{2, \,2}}^{4}\,{Z_{0, \,0}}\,{Z_{1, \,1}
} }
}
\end{maplelatex}
\end{maplegroup}
\emptyline
\begin{maplegroup}
\begin{mapleinput}
\mapleinline{active}{1d}{T[13]:=Gen(Z[0,0]*Z[1,1]*Z[2,2]*Z[0,1]^3);}{%
}
\end{mapleinput}
\mapleresult
\begin{maplelatex}
\mapleinline{inert}{2d}{T[13] :=
Z[0,0]*Z[1,1]*Z[2,2]*Z[0,1]^3+Z[1,0]*Z[2,1]*Z[0,2]*Z[1,1]^3+Z[2,0]*Z[0
,1]*Z[1,2]*Z[2,1]^3+Z[0,1]*Z[1,2]*Z[2,0]*Z[0,2]^3+Z[1,1]*Z[2,2]*Z[0,0]
*Z[1,2]^3+Z[2,1]*Z[0,2]*Z[1,0]*Z[2,2]^3+Z[0,2]*Z[1,0]*Z[2,1]*Z[0,0]^3+
Z[1,2]*Z[2,0]*Z[0,1]*Z[1,0]^3+Z[2,2]*Z[0,0]*Z[1,1]*Z[2,0]^3;}{%
\maplemultiline{
{T_{13}} := {Z_{0, \,0}}\,{Z_{1, \,1}}\,{Z_{2, \,2}}\,{Z_{0, \,1}
}^{3} + {Z_{1, \,0}}\,{Z_{2, \,1}}\,{Z_{0, \,2}}\,{Z_{1, \,1}}^{3
} + {Z_{2, \,0}}\,{Z_{0, \,1}}\,{Z_{1, \,2}}\,{Z_{2, \,1}}^{3}
 \\
\mbox{} + {Z_{0, \,1}}\,{Z_{1, \,2}}\,{Z_{2, \,0}}\,{Z_{0, \,2}}
^{3} + {Z_{1, \,1}}\,{Z_{2, \,2}}\,{Z_{0, \,0}}\,{Z_{1, \,2}}^{3}
 + {Z_{2, \,1}}\,{Z_{0, \,2}}\,{Z_{1, \,0}}\,{Z_{2, \,2}}^{3} \\
\mbox{} + {Z_{0, \,2}}\,{Z_{1, \,0}}\,{Z_{2, \,1}}\,{Z_{0, \,0}}
^{3} + {Z_{1, \,2}}\,{Z_{2, \,0}}\,{Z_{0, \,1}}\,{Z_{1, \,0}}^{3}
 + {Z_{2, \,2}}\,{Z_{0, \,0}}\,{Z_{1, \,1}}\,{Z_{2, \,0}}^{3} }
}
\end{maplelatex}
\end{maplegroup}
\emptyline
\begin{maplegroup}
\begin{mapleinput}
\mapleinline{active}{1d}{T[14]:=Gen(Z[0,0]*Z[1,1]*Z[2,2]*Z[0,2]^3);}{%
}
\end{mapleinput}
\mapleresult
\begin{maplelatex}
\mapleinline{inert}{2d}{T[14] :=
Z[0,0]*Z[1,1]*Z[2,2]*Z[0,2]^3+Z[1,0]*Z[2,1]*Z[0,2]*Z[1,2]^3+Z[2,0]*Z[0
,1]*Z[1,2]*Z[2,2]^3+Z[0,1]*Z[1,2]*Z[2,0]*Z[0,0]^3+Z[1,1]*Z[2,2]*Z[0,0]
*Z[1,0]^3+Z[2,1]*Z[0,2]*Z[1,0]*Z[2,0]^3+Z[0,2]*Z[1,0]*Z[2,1]*Z[0,1]^3+
Z[1,2]*Z[2,0]*Z[0,1]*Z[1,1]^3+Z[2,2]*Z[0,0]*Z[1,1]*Z[2,1]^3;}{%
\maplemultiline{
{T_{14}} := {Z_{0, \,0}}\,{Z_{1, \,1}}\,{Z_{2, \,2}}\,{Z_{0, \,2}
}^{3} + {Z_{1, \,0}}\,{Z_{2, \,1}}\,{Z_{0, \,2}}\,{Z_{1, \,2}}^{3
} + {Z_{2, \,0}}\,{Z_{0, \,1}}\,{Z_{1, \,2}}\,{Z_{2, \,2}}^{3}
 \\
\mbox{} + {Z_{0, \,1}}\,{Z_{1, \,2}}\,{Z_{2, \,0}}\,{Z_{0, \,0}}
^{3} + {Z_{1, \,1}}\,{Z_{2, \,2}}\,{Z_{0, \,0}}\,{Z_{1, \,0}}^{3}
 + {Z_{2, \,1}}\,{Z_{0, \,2}}\,{Z_{1, \,0}}\,{Z_{2, \,0}}^{3} \\
\mbox{} + {Z_{0, \,2}}\,{Z_{1, \,0}}\,{Z_{2, \,1}}\,{Z_{0, \,1}}
^{3} + {Z_{1, \,2}}\,{Z_{2, \,0}}\,{Z_{0, \,1}}\,{Z_{1, \,1}}^{3}
 + {Z_{2, \,2}}\,{Z_{0, \,0}}\,{Z_{1, \,1}}\,{Z_{2, \,1}}^{3} }
}
\end{maplelatex}
\end{maplegroup}
\emptyline
\begin{maplegroup}
\begin{mapleinput}
\mapleinline{active}{1d}{T[15]:=Gen(Z[0,0]^4*Z[1,2]*Z[2,1]);}{%
}
\end{mapleinput}
\mapleresult
\begin{maplelatex}
\mapleinline{inert}{2d}{T[15] :=
Z[0,0]^4*Z[1,2]*Z[2,1]+Z[1,0]^4*Z[2,2]*Z[0,1]+Z[2,0]^4*Z[0,2]*Z[1,1]+Z
[0,1]^4*Z[1,0]*Z[2,2]+Z[1,1]^4*Z[2,0]*Z[0,2]+Z[2,1]^4*Z[0,0]*Z[1,2]+Z[
0,2]^4*Z[1,1]*Z[2,0]+Z[1,2]^4*Z[2,1]*Z[0,0]+Z[2,2]^4*Z[0,1]*Z[1,0];}{%
\maplemultiline{
{T_{15}} := {Z_{0, \,0}}^{4}\,{Z_{1, \,2}}\,{Z_{2, \,1}} + {Z_{1
, \,0}}^{4}\,{Z_{2, \,2}}\,{Z_{0, \,1}} + {Z_{2, \,0}}^{4}\,{Z_{0
, \,2}}\,{Z_{1, \,1}}\\ + {Z_{0, \,1}}^{4}\,{Z_{1, \,0}}\,{Z_{2, \,
2}} + {Z_{1, \,1}}^{4}\,{Z_{2, \,0}}\,{Z_{0, \,2}} 
 + {Z_{2, \,1}}^{4}\,{Z_{0, \,0}}\,{Z_{1, \,2}}\\ + {Z_{0, 
\,2}}^{4}\,{Z_{1, \,1}}\,{Z_{2, \,0}} + {Z_{1, \,2}}^{4}\,{Z_{2, 
\,1}}\,{Z_{0, \,0}} + {Z_{2, \,2}}^{4}\,{Z_{0, \,1}}\,{Z_{1, \,0}
} }
}
\end{maplelatex}
\end{maplegroup}
\emptyline
\begin{maplegroup}
\begin{mapleinput}
\mapleinline{active}{1d}{T[16]:=Gen(Z[0,0]*Z[1,2]*Z[2,1]*Z[1,0]^3);}{%
}
\end{mapleinput}
\mapleresult
\begin{maplelatex}
\mapleinline{inert}{2d}{T[16] :=
Z[0,0]*Z[1,2]*Z[2,1]*Z[1,0]^3+Z[1,0]*Z[2,2]*Z[0,1]*Z[2,0]^3+Z[2,0]*Z[0
,2]*Z[1,1]*Z[0,0]^3+Z[0,1]*Z[1,0]*Z[2,2]*Z[1,1]^3+Z[1,1]*Z[2,0]*Z[0,2]
*Z[2,1]^3+Z[2,1]*Z[0,0]*Z[1,2]*Z[0,1]^3+Z[0,2]*Z[1,1]*Z[2,0]*Z[1,2]^3+
Z[1,2]*Z[2,1]*Z[0,0]*Z[2,2]^3+Z[2,2]*Z[0,1]*Z[1,0]*Z[0,2]^3;}{%
\maplemultiline{
{T_{16}} := {Z_{0, \,0}}\,{Z_{1, \,2}}\,{Z_{2, \,1}}\,{Z_{1, \,0}
}^{3} + {Z_{1, \,0}}\,{Z_{2, \,2}}\,{Z_{0, \,1}}\,{Z_{2, \,0}}^{3
} + {Z_{2, \,0}}\,{Z_{0, \,2}}\,{Z_{1, \,1}}\,{Z_{0, \,0}}^{3}
 \\
\mbox{} + {Z_{0, \,1}}\,{Z_{1, \,0}}\,{Z_{2, \,2}}\,{Z_{1, \,1}}
^{3} + {Z_{1, \,1}}\,{Z_{2, \,0}}\,{Z_{0, \,2}}\,{Z_{2, \,1}}^{3}
 + {Z_{2, \,1}}\,{Z_{0, \,0}}\,{Z_{1, \,2}}\,{Z_{0, \,1}}^{3} \\
\mbox{} + {Z_{0, \,2}}\,{Z_{1, \,1}}\,{Z_{2, \,0}}\,{Z_{1, \,2}}
^{3} + {Z_{1, \,2}}\,{Z_{2, \,1}}\,{Z_{0, \,0}}\,{Z_{2, \,2}}^{3}
 + {Z_{2, \,2}}\,{Z_{0, \,1}}\,{Z_{1, \,0}}\,{Z_{0, \,2}}^{3} }
}
\end{maplelatex}
\end{maplegroup}
\emptyline
\begin{maplegroup}
\begin{mapleinput}
\mapleinline{active}{1d}{T[17]:=Gen(Z[0,0]*Z[1,2]*Z[2,1]*Z[2,0]^3);}{%
}
\end{mapleinput}
\mapleresult
\begin{maplelatex}
\mapleinline{inert}{2d}{T[17] :=
Z[0,0]*Z[1,2]*Z[2,1]*Z[2,0]^3+Z[1,0]*Z[2,2]*Z[0,1]*Z[0,0]^3+Z[2,0]*Z[0
,2]*Z[1,1]*Z[1,0]^3+Z[0,1]*Z[1,0]*Z[2,2]*Z[2,1]^3+Z[1,1]*Z[2,0]*Z[0,2]
*Z[0,1]^3+Z[2,1]*Z[0,0]*Z[1,2]*Z[1,1]^3+Z[0,2]*Z[1,1]*Z[2,0]*Z[2,2]^3+
Z[1,2]*Z[2,1]*Z[0,0]*Z[0,2]^3+Z[2,2]*Z[0,1]*Z[1,0]*Z[1,2]^3;}{%
\maplemultiline{
{T_{17}} := {Z_{0, \,0}}\,{Z_{1, \,2}}\,{Z_{2, \,1}}\,{Z_{2, \,0}
}^{3} + {Z_{1, \,0}}\,{Z_{2, \,2}}\,{Z_{0, \,1}}\,{Z_{0, \,0}}^{3
} + {Z_{2, \,0}}\,{Z_{0, \,2}}\,{Z_{1, \,1}}\,{Z_{1, \,0}}^{3}
 \\
\mbox{} + {Z_{0, \,1}}\,{Z_{1, \,0}}\,{Z_{2, \,2}}\,{Z_{2, \,1}}
^{3} + {Z_{1, \,1}}\,{Z_{2, \,0}}\,{Z_{0, \,2}}\,{Z_{0, \,1}}^{3}
 + {Z_{2, \,1}}\,{Z_{0, \,0}}\,{Z_{1, \,2}}\,{Z_{1, \,1}}^{3} \\
\mbox{} + {Z_{0, \,2}}\,{Z_{1, \,1}}\,{Z_{2, \,0}}\,{Z_{2, \,2}}
^{3} + {Z_{1, \,2}}\,{Z_{2, \,1}}\,{Z_{0, \,0}}\,{Z_{0, \,2}}^{3}
 + {Z_{2, \,2}}\,{Z_{0, \,1}}\,{Z_{1, \,0}}\,{Z_{1, \,2}}^{3} }
}
\end{maplelatex}
\end{maplegroup}
\emptyline
\begin{maplegroup}
\begin{mapleinput}
\mapleinline{active}{1d}{T[18]:=1/3*Gen(Z[0,0]^2*Z[0,1]^2*Z[0,2]^2);}{%
}
\end{mapleinput}
\mapleresult
\begin{maplelatex}
\mapleinline{inert}{2d}{T[18] :=
Z[0,0]^2*Z[0,1]^2*Z[0,2]^2+Z[1,0]^2*Z[1,1]^2*Z[1,2]^2+Z[2,0]^2*Z[2,1]^
2*Z[2,2]^2;}{%
\[
{T_{18}} := {Z_{0, \,0}}^{2}\,{Z_{0, \,1}}^{2}\,{Z_{0, \,2}}^{2}
 + {Z_{1, \,0}}^{2}\,{Z_{1, \,1}}^{2}\,{Z_{1, \,2}}^{2} + {Z_{2, 
\,0}}^{2}\,{Z_{2, \,1}}^{2}\,{Z_{2, \,2}}^{2}
\]
}
\end{maplelatex}
\end{maplegroup}
\emptyline
\begin{maplegroup}
\begin{mapleinput}
\mapleinline{active}{1d}{T[19]:=1/3*Gen(Z[0,0]^2*Z[1,0]^2*Z[2,0]^2);}{%
}
\end{mapleinput}
\mapleresult
\begin{maplelatex}
\mapleinline{inert}{2d}{T[19] :=
Z[0,0]^2*Z[1,0]^2*Z[2,0]^2+Z[0,1]^2*Z[1,1]^2*Z[2,1]^2+Z[0,2]^2*Z[1,2]^
2*Z[2,2]^2;}{%
\[
{T_{19}} := {Z_{0, \,0}}^{2}\,{Z_{1, \,0}}^{2}\,{Z_{2, \,0}}^{2}
 + {Z_{0, \,1}}^{2}\,{Z_{1, \,1}}^{2}\,{Z_{2, \,1}}^{2} + {Z_{0, 
\,2}}^{2}\,{Z_{1, \,2}}^{2}\,{Z_{2, \,2}}^{2}
\]
}
\end{maplelatex}
\end{maplegroup}
\emptyline
\begin{maplegroup}
\begin{mapleinput}
\mapleinline{active}{1d}{T[20]:=1/3*Gen(Z[0,0]^2*Z[1,1]^2*Z[2,2]^2);}{%
}
\end{mapleinput}
\mapleresult
\begin{maplelatex}
\mapleinline{inert}{2d}{T[20] :=
Z[0,0]^2*Z[1,1]^2*Z[2,2]^2+Z[1,0]^2*Z[2,1]^2*Z[0,2]^2+Z[2,0]^2*Z[0,1]^
2*Z[1,2]^2;}{%
\[
{T_{20}} := {Z_{0, \,0}}^{2}\,{Z_{1, \,1}}^{2}\,{Z_{2, \,2}}^{2}
 + {Z_{1, \,0}}^{2}\,{Z_{2, \,1}}^{2}\,{Z_{0, \,2}}^{2} + {Z_{2, 
\,0}}^{2}\,{Z_{0, \,1}}^{2}\,{Z_{1, \,2}}^{2}
\]
}
\end{maplelatex}
\end{maplegroup}
\emptyline
\begin{maplegroup}
\begin{mapleinput}
\mapleinline{active}{1d}{T[21]:=1/3*Gen(Z[0,0]^2*Z[2,1]^2*Z[1,2]^2);}{%
}
\end{mapleinput}
\mapleresult
\begin{maplelatex}
\mapleinline{inert}{2d}{T[21] :=
Z[0,0]^2*Z[2,1]^2*Z[1,2]^2+Z[1,0]^2*Z[0,1]^2*Z[2,2]^2+Z[2,0]^2*Z[1,1]^
2*Z[0,2]^2;}{%
\[
{T_{21}} := {Z_{0, \,0}}^{2}\,{Z_{2, \,1}}^{2}\,{Z_{1, \,2}}^{2}
 + {Z_{1, \,0}}^{2}\,{Z_{0, \,1}}^{2}\,{Z_{2, \,2}}^{2} + {Z_{2, 
\,0}}^{2}\,{Z_{1, \,1}}^{2}\,{Z_{0, \,2}}^{2}
\]
}
\end{maplelatex}
\end{maplegroup}
\emptyline
\begin{maplegroup}
\begin{mapleinput}
\mapleinline{active}{1d}{T[22]:=1/3*Gen(Z[0,0]*Z[0,1]*Z[0,2]*Z[1,0]*Z[1,1]*Z[1,2]);}{%
}
\end{mapleinput}
\mapleresult
\begin{maplelatex}
\mapleinline{inert}{2d}{T[22] :=
Z[0,0]*Z[0,1]*Z[0,2]*Z[1,0]*Z[1,1]*Z[1,2]+Z[1,0]*Z[1,1]*Z[1,2]*Z[2,0]*
Z[2,1]*Z[2,2]+Z[2,0]*Z[2,1]*Z[2,2]*Z[0,0]*Z[0,1]*Z[0,2];}{%
\maplemultiline{
{T_{22}} := {Z_{0, \,0}}\,{Z_{0, \,1}}\,{Z_{0, \,2}}\,{Z_{1, \,0}
}\,{Z_{1, \,1}}\,{Z_{1, \,2}} + {Z_{1, \,0}}\,{Z_{1, \,1}}\,{Z_{1
, \,2}}\,{Z_{2, \,0}}\,{Z_{2, \,1}}\,{Z_{2, \,2}} \\
\mbox{} + {Z_{2, \,0}}\,{Z_{2, \,1}}\,{Z_{2, \,2}}\,{Z_{0, \,0}}
\,{Z_{0, \,1}}\,{Z_{0, \,2}} }
}
\end{maplelatex}
\end{maplegroup}
\emptyline
\begin{maplegroup}
\begin{mapleinput}
\mapleinline{active}{1d}{T[23]:=1/3*Gen(Z[0,0]*Z[1,0]*Z[2,0]*Z[0,1]*Z[1,1]*Z[2,1]);}{%
}
\end{mapleinput}
\mapleresult
\begin{maplelatex}
\mapleinline{inert}{2d}{T[23] :=
Z[0,0]*Z[1,0]*Z[2,0]*Z[0,1]*Z[1,1]*Z[2,1]+Z[0,1]*Z[1,1]*Z[2,1]*Z[0,2]*
Z[1,2]*Z[2,2]+Z[0,2]*Z[1,2]*Z[2,2]*Z[0,0]*Z[1,0]*Z[2,0];}{%
\maplemultiline{
{T_{23}} := {Z_{0, \,0}}\,{Z_{1, \,0}}\,{Z_{2, \,0}}\,{Z_{0, \,1}
}\,{Z_{1, \,1}}\,{Z_{2, \,1}} + {Z_{0, \,1}}\,{Z_{1, \,1}}\,{Z_{2
, \,1}}\,{Z_{0, \,2}}\,{Z_{1, \,2}}\,{Z_{2, \,2}} \\
\mbox{} + {Z_{0, \,2}}\,{Z_{1, \,2}}\,{Z_{2, \,2}}\,{Z_{0, \,0}}
\,{Z_{1, \,0}}\,{Z_{2, \,0}} }
}
\end{maplelatex}
\end{maplegroup}
\emptyline
\begin{maplegroup}
\begin{mapleinput}
\mapleinline{active}{1d}{T[24]:=1/3*Gen(Z[0,0]*Z[1,1]*Z[2,2]*Z[0,1]*Z[1,2]*Z[2,0]);}{%
}
\end{mapleinput}

\mapleresult
\begin{maplelatex}
\mapleinline{inert}{2d}{T[24] :=
Z[0,0]*Z[1,1]*Z[2,2]*Z[0,1]*Z[1,2]*Z[2,0]+Z[1,0]*Z[2,1]*Z[0,2]*Z[1,1]*
Z[2,2]*Z[0,0]+Z[2,0]*Z[0,1]*Z[1,2]*Z[2,1]*Z[0,2]*Z[1,0];}{%
\maplemultiline{
{T_{24}} := {Z_{0, \,0}}\,{Z_{1, \,1}}\,{Z_{2, \,2}}\,{Z_{0, \,1}
}\,{Z_{1, \,2}}\,{Z_{2, \,0}} + {Z_{1, \,0}}\,{Z_{2, \,1}}\,{Z_{0
, \,2}}\,{Z_{1, \,1}}\,{Z_{2, \,2}}\,{Z_{0, \,0}} \\
\mbox{} + {Z_{2, \,0}}\,{Z_{0, \,1}}\,{Z_{1, \,2}}\,{Z_{2, \,1}}
\,{Z_{0, \,2}}\,{Z_{1, \,0}} }
}
\end{maplelatex}
\end{maplegroup}
\begin{maplegroup}
\begin{mapleinput}
\mapleinline{active}{1d}{T[25]:=1/3*Gen(Z[0,0]*Z[1,2]*Z[2,1]*Z[0,1]*Z[1,0]*Z[2,2]);}{%
}
\end{mapleinput}
\mapleresult
\begin{maplelatex}
\mapleinline{inert}{2d}{T[25] :=
Z[0,0]*Z[1,2]*Z[2,1]*Z[0,1]*Z[1,0]*Z[2,2]+Z[1,0]*Z[2,2]*Z[0,1]*Z[1,1]*
Z[2,0]*Z[0,2]+Z[2,0]*Z[0,2]*Z[1,1]*Z[2,1]*Z[0,0]*Z[1,2];}{%
\maplemultiline{
{T_{25}} := {Z_{0, \,0}}\,{Z_{1, \,2}}\,{Z_{2, \,1}}\,{Z_{0, \,1}
}\,{Z_{1, \,0}}\,{Z_{2, \,2}} + {Z_{1, \,0}}\,{Z_{2, \,2}}\,{Z_{0
, \,1}}\,{Z_{1, \,1}}\,{Z_{2, \,0}}\,{Z_{0, \,2}} \\
\mbox{} + {Z_{2, \,0}}\,{Z_{0, \,2}}\,{Z_{1, \,1}}\,{Z_{2, \,1}}
\,{Z_{0, \,0}}\,{Z_{1, \,2}} }
}
\end{maplelatex}
\end{maplegroup}
\emptyline
\begin{maplegroup}
\begin{mapleinput}
\mapleinline{active}{1d}{T[26]:=Gen(Z[0,0]^2*Z[0,1]*Z[1,1]*Z[1,2]^2);}{%
}
\end{mapleinput}
\mapleresult
\begin{maplelatex}
\mapleinline{inert}{2d}{T[26] :=
Z[0,0]^2*Z[0,1]*Z[1,1]*Z[1,2]^2+Z[1,0]^2*Z[1,1]*Z[2,1]*Z[2,2]^2+Z[2,0]
^2*Z[2,1]*Z[0,1]*Z[0,2]^2+Z[0,1]^2*Z[0,2]*Z[1,2]*Z[1,0]^2+Z[1,1]^2*Z[1
,2]*Z[2,2]*Z[2,0]^2+Z[2,1]^2*Z[2,2]*Z[0,2]*Z[0,0]^2+Z[0,2]^2*Z[0,0]*Z[
1,0]*Z[1,1]^2+Z[1,2]^2*Z[1,0]*Z[2,0]*Z[2,1]^2+Z[2,2]^2*Z[2,0]*Z[0,0]*Z
[0,1]^2;}{%
\maplemultiline{
{T_{26}} := {Z_{0, \,0}}^{2}\,{Z_{0, \,1}}\,{Z_{1, \,1}}\,{Z_{1, 
\,2}}^{2} + {Z_{1, \,0}}^{2}\,{Z_{1, \,1}}\,{Z_{2, \,1}}\,{Z_{2, 
\,2}}^{2} + {Z_{2, \,0}}^{2}\,{Z_{2, \,1}}\,{Z_{0, \,1}}\,{Z_{0, 
\,2}}^{2} \\
\mbox{} + {Z_{0, \,1}}^{2}\,{Z_{0, \,2}}\,{Z_{1, \,2}}\,{Z_{1, \,
0}}^{2} + {Z_{1, \,1}}^{2}\,{Z_{1, \,2}}\,{Z_{2, \,2}}\,{Z_{2, \,
0}}^{2} + {Z_{2, \,1}}^{2}\,{Z_{2, \,2}}\,{Z_{0, \,2}}\,{Z_{0, \,
0}}^{2} \\
\mbox{} + {Z_{0, \,2}}^{2}\,{Z_{0, \,0}}\,{Z_{1, \,0}}\,{Z_{1, \,
1}}^{2} + {Z_{1, \,2}}^{2}\,{Z_{1, \,0}}\,{Z_{2, \,0}}\,{Z_{2, \,
1}}^{2} + {Z_{2, \,2}}^{2}\,{Z_{2, \,0}}\,{Z_{0, \,0}}\,{Z_{0, \,
1}}^{2} }
}
\end{maplelatex}
\end{maplegroup}
\emptyline
\begin{maplegroup}
\begin{mapleinput}
\mapleinline{active}{1d}{T[27]:=Gen(Z[0,0]^2*Z[0,2]*Z[1,2]*Z[1,1]^2);}{%
}
\end{mapleinput}
\mapleresult
\begin{maplelatex}
\mapleinline{inert}{2d}{T[27] :=
Z[0,0]^2*Z[0,2]*Z[1,2]*Z[1,1]^2+Z[1,0]^2*Z[1,2]*Z[2,2]*Z[2,1]^2+Z[2,0]
^2*Z[2,2]*Z[0,2]*Z[0,1]^2+Z[0,1]^2*Z[0,0]*Z[1,0]*Z[1,2]^2+Z[1,1]^2*Z[1
,0]*Z[2,0]*Z[2,2]^2+Z[2,1]^2*Z[2,0]*Z[0,0]*Z[0,2]^2+Z[0,2]^2*Z[0,1]*Z[
1,1]*Z[1,0]^2+Z[1,2]^2*Z[1,1]*Z[2,1]*Z[2,0]^2+Z[2,2]^2*Z[2,1]*Z[0,1]*Z
[0,0]^2;}{%
\maplemultiline{
{T_{27}} := {Z_{0, \,0}}^{2}\,{Z_{0, \,2}}\,{Z_{1, \,2}}\,{Z_{1, 
\,1}}^{2} + {Z_{1, \,0}}^{2}\,{Z_{1, \,2}}\,{Z_{2, \,2}}\,{Z_{2, 
\,1}}^{2} + {Z_{2, \,0}}^{2}\,{Z_{2, \,2}}\,{Z_{0, \,2}}\,{Z_{0, 
\,1}}^{2} \\
\mbox{} + {Z_{0, \,1}}^{2}\,{Z_{0, \,0}}\,{Z_{1, \,0}}\,{Z_{1, \,
2}}^{2} + {Z_{1, \,1}}^{2}\,{Z_{1, \,0}}\,{Z_{2, \,0}}\,{Z_{2, \,
2}}^{2} + {Z_{2, \,1}}^{2}\,{Z_{2, \,0}}\,{Z_{0, \,0}}\,{Z_{0, \,
2}}^{2} \\
\mbox{} + {Z_{0, \,2}}^{2}\,{Z_{0, \,1}}\,{Z_{1, \,1}}\,{Z_{1, \,
0}}^{2} + {Z_{1, \,2}}^{2}\,{Z_{1, \,1}}\,{Z_{2, \,1}}\,{Z_{2, \,
0}}^{2} + {Z_{2, \,2}}^{2}\,{Z_{2, \,1}}\,{Z_{0, \,1}}\,{Z_{0, \,
0}}^{2} }
}
\end{maplelatex}
\end{maplegroup}
\emptyline
\begin{maplegroup}
\begin{mapleinput}
\mapleinline{active}{1d}{T[28]:=Gen(Z[0,0]^2*Z[1,1]*Z[2,1]*Z[0,2]^2);}{%
}
\end{mapleinput}

\mapleresult
\begin{maplelatex}
\mapleinline{inert}{2d}{T[28] :=
Z[0,0]^2*Z[1,1]*Z[2,1]*Z[0,2]^2+Z[1,0]^2*Z[2,1]*Z[0,1]*Z[1,2]^2+Z[2,0]
^2*Z[0,1]*Z[1,1]*Z[2,2]^2+Z[0,1]^2*Z[1,2]*Z[2,2]*Z[0,0]^2+Z[1,1]^2*Z[2
,2]*Z[0,2]*Z[1,0]^2+Z[2,1]^2*Z[0,2]*Z[1,2]*Z[2,0]^2+Z[0,2]^2*Z[1,0]*Z[
2,0]*Z[0,1]^2+Z[1,2]^2*Z[2,0]*Z[0,0]*Z[1,1]^2+Z[2,2]^2*Z[0,0]*Z[1,0]*Z
[2,1]^2;}{%
\maplemultiline{
{T_{28}} := {Z_{0, \,0}}^{2}\,{Z_{1, \,1}}\,{Z_{2, \,1}}\,{Z_{0, 
\,2}}^{2} + {Z_{1, \,0}}^{2}\,{Z_{2, \,1}}\,{Z_{0, \,1}}\,{Z_{1, 
\,2}}^{2} + {Z_{2, \,0}}^{2}\,{Z_{0, \,1}}\,{Z_{1, \,1}}\,{Z_{2, 
\,2}}^{2} \\
\mbox{} + {Z_{0, \,1}}^{2}\,{Z_{1, \,2}}\,{Z_{2, \,2}}\,{Z_{0, \,
0}}^{2} + {Z_{1, \,1}}^{2}\,{Z_{2, \,2}}\,{Z_{0, \,2}}\,{Z_{1, \,
0}}^{2} + {Z_{2, \,1}}^{2}\,{Z_{0, \,2}}\,{Z_{1, \,2}}\,{Z_{2, \,
0}}^{2} \\
\mbox{} + {Z_{0, \,2}}^{2}\,{Z_{1, \,0}}\,{Z_{2, \,0}}\,{Z_{0, \,
1}}^{2} + {Z_{1, \,2}}^{2}\,{Z_{2, \,0}}\,{Z_{0, \,0}}\,{Z_{1, \,
1}}^{2} + {Z_{2, \,2}}^{2}\,{Z_{0, \,0}}\,{Z_{1, \,0}}\,{Z_{2, \,
1}}^{2} }
}
\end{maplelatex}
\end{maplegroup}
\emptyline
\begin{maplegroup}
\begin{mapleinput}
\mapleinline{active}{1d}{T[29]:=Gen(Z[0,0]^2*Z[1,0]*Z[1,1]*Z[2,1]^2);}{%
}
\end{mapleinput}

\mapleresult
\begin{maplelatex}
\mapleinline{inert}{2d}{T[29] :=
Z[0,0]^2*Z[1,0]*Z[1,1]*Z[2,1]^2+Z[1,0]^2*Z[2,0]*Z[2,1]*Z[0,1]^2+Z[2,0]
^2*Z[0,0]*Z[0,1]*Z[1,1]^2+Z[0,1]^2*Z[1,1]*Z[1,2]*Z[2,2]^2+Z[1,1]^2*Z[2
,1]*Z[2,2]*Z[0,2]^2+Z[2,1]^2*Z[0,1]*Z[0,2]*Z[1,2]^2+Z[0,2]^2*Z[1,2]*Z[
1,0]*Z[2,0]^2+Z[1,2]^2*Z[2,2]*Z[2,0]*Z[0,0]^2+Z[2,2]^2*Z[0,2]*Z[0,0]*Z
[1,0]^2;}{%
\maplemultiline{
{T_{29}} := {Z_{0, \,0}}^{2}\,{Z_{1, \,0}}\,{Z_{1, \,1}}\,{Z_{2, 
\,1}}^{2} + {Z_{1, \,0}}^{2}\,{Z_{2, \,0}}\,{Z_{2, \,1}}\,{Z_{0, 
\,1}}^{2} + {Z_{2, \,0}}^{2}\,{Z_{0, \,0}}\,{Z_{0, \,1}}\,{Z_{1, 
\,1}}^{2} \\
\mbox{} + {Z_{0, \,1}}^{2}\,{Z_{1, \,1}}\,{Z_{1, \,2}}\,{Z_{2, \,
2}}^{2} + {Z_{1, \,1}}^{2}\,{Z_{2, \,1}}\,{Z_{2, \,2}}\,{Z_{0, \,
2}}^{2} + {Z_{2, \,1}}^{2}\,{Z_{0, \,1}}\,{Z_{0, \,2}}\,{Z_{1, \,
2}}^{2} \\
\mbox{} + {Z_{0, \,2}}^{2}\,{Z_{1, \,2}}\,{Z_{1, \,0}}\,{Z_{2, \,
0}}^{2} + {Z_{1, \,2}}^{2}\,{Z_{2, \,2}}\,{Z_{2, \,0}}\,{Z_{0, \,
0}}^{2} + {Z_{2, \,2}}^{2}\,{Z_{0, \,2}}\,{Z_{0, \,0}}\,{Z_{1, \,
0}}^{2} }
}
\end{maplelatex}
\end{maplegroup}
\emptyline
\begin{maplegroup}
\begin{mapleinput}
\mapleinline{active}{1d}{T[30]:=Gen(Z[0,0]^2*Z[1,0]*Z[1,2]*Z[2,2]^2);}{%
}
\end{mapleinput}

\mapleresult
\begin{maplelatex}
\mapleinline{inert}{2d}{T[30] :=
Z[0,0]^2*Z[1,0]*Z[1,2]*Z[2,2]^2+Z[1,0]^2*Z[2,0]*Z[2,2]*Z[0,2]^2+Z[2,0]
^2*Z[0,0]*Z[0,2]*Z[1,2]^2+Z[0,1]^2*Z[1,1]*Z[1,0]*Z[2,0]^2+Z[1,1]^2*Z[2
,1]*Z[2,0]*Z[0,0]^2+Z[2,1]^2*Z[0,1]*Z[0,0]*Z[1,0]^2+Z[0,2]^2*Z[1,2]*Z[
1,1]*Z[2,1]^2+Z[1,2]^2*Z[2,2]*Z[2,1]*Z[0,1]^2+Z[2,2]^2*Z[0,2]*Z[0,1]*Z
[1,1]^2;}{%
\maplemultiline{
{T_{30}} := {Z_{0, \,0}}^{2}\,{Z_{1, \,0}}\,{Z_{1, \,2}}\,{Z_{2, 
\,2}}^{2} + {Z_{1, \,0}}^{2}\,{Z_{2, \,0}}\,{Z_{2, \,2}}\,{Z_{0, 
\,2}}^{2} + {Z_{2, \,0}}^{2}\,{Z_{0, \,0}}\,{Z_{0, \,2}}\,{Z_{1, 
\,2}}^{2} \\
\mbox{} + {Z_{0, \,1}}^{2}\,{Z_{1, \,1}}\,{Z_{1, \,0}}\,{Z_{2, \,
0}}^{2} + {Z_{1, \,1}}^{2}\,{Z_{2, \,1}}\,{Z_{2, \,0}}\,{Z_{0, \,
0}}^{2} + {Z_{2, \,1}}^{2}\,{Z_{0, \,1}}\,{Z_{0, \,0}}\,{Z_{1, \,
0}}^{2} \\
\mbox{} + {Z_{0, \,2}}^{2}\,{Z_{1, \,2}}\,{Z_{1, \,1}}\,{Z_{2, \,
1}}^{2} + {Z_{1, \,2}}^{2}\,{Z_{2, \,2}}\,{Z_{2, \,1}}\,{Z_{0, \,
1}}^{2} + {Z_{2, \,2}}^{2}\,{Z_{0, \,2}}\,{Z_{0, \,1}}\,{Z_{1, \,
1}}^{2} }
}
\end{maplelatex}
\end{maplegroup}
\emptyline
\begin{maplegroup}
\begin{mapleinput}
\mapleinline{active}{1d}{T[31]:=Gen(Z[0,0]^2*Z[1,1]*Z[1,2]*Z[2,0]^2);}{%
}
\end{mapleinput}

\mapleresult
\begin{maplelatex}
\mapleinline{inert}{2d}{T[31] :=
Z[0,0]^2*Z[1,1]*Z[1,2]*Z[2,0]^2+Z[1,0]^2*Z[2,1]*Z[2,2]*Z[0,0]^2+Z[2,0]
^2*Z[0,1]*Z[0,2]*Z[1,0]^2+Z[0,1]^2*Z[1,2]*Z[1,0]*Z[2,1]^2+Z[1,1]^2*Z[2
,2]*Z[2,0]*Z[0,1]^2+Z[2,1]^2*Z[0,2]*Z[0,0]*Z[1,1]^2+Z[0,2]^2*Z[1,0]*Z[
1,1]*Z[2,2]^2+Z[1,2]^2*Z[2,0]*Z[2,1]*Z[0,2]^2+Z[2,2]^2*Z[0,0]*Z[0,1]*Z
[1,2]^2;}{%
\maplemultiline{
{T_{31}} := {Z_{0, \,0}}^{2}\,{Z_{1, \,1}}\,{Z_{1, \,2}}\,{Z_{2, 
\,0}}^{2} + {Z_{1, \,0}}^{2}\,{Z_{2, \,1}}\,{Z_{2, \,2}}\,{Z_{0, 
\,0}}^{2} + {Z_{2, \,0}}^{2}\,{Z_{0, \,1}}\,{Z_{0, \,2}}\,{Z_{1, 
\,0}}^{2} \\
\mbox{} + {Z_{0, \,1}}^{2}\,{Z_{1, \,2}}\,{Z_{1, \,0}}\,{Z_{2, \,
1}}^{2} + {Z_{1, \,1}}^{2}\,{Z_{2, \,2}}\,{Z_{2, \,0}}\,{Z_{0, \,
1}}^{2} + {Z_{2, \,1}}^{2}\,{Z_{0, \,2}}\,{Z_{0, \,0}}\,{Z_{1, \,
1}}^{2} \\
\mbox{} + {Z_{0, \,2}}^{2}\,{Z_{1, \,0}}\,{Z_{1, \,1}}\,{Z_{2, \,
2}}^{2} + {Z_{1, \,2}}^{2}\,{Z_{2, \,0}}\,{Z_{2, \,1}}\,{Z_{0, \,
2}}^{2} + {Z_{2, \,2}}^{2}\,{Z_{0, \,0}}\,{Z_{0, \,1}}\,{Z_{1, \,
2}}^{2} }
}
\end{maplelatex}
\end{maplegroup}
\emptyline
\begin{maplegroup}
\begin{mapleinput}
\mapleinline{active}{1d}{T[32]:=Gen(Z[0,0]^2*Z[0,1]*Z[1,2]*Z[1,0]^2);}{%
}
\end{mapleinput}

\mapleresult
\begin{maplelatex}
\mapleinline{inert}{2d}{T[32] :=
Z[0,0]^2*Z[0,1]*Z[1,2]*Z[1,0]^2+Z[1,0]^2*Z[1,1]*Z[2,2]*Z[2,0]^2+Z[2,0]
^2*Z[2,1]*Z[0,2]*Z[0,0]^2+Z[0,1]^2*Z[0,2]*Z[1,0]*Z[1,1]^2+Z[1,1]^2*Z[1
,2]*Z[2,0]*Z[2,1]^2+Z[2,1]^2*Z[2,2]*Z[0,0]*Z[0,1]^2+Z[0,2]^2*Z[0,0]*Z[
1,1]*Z[1,2]^2+Z[1,2]^2*Z[1,0]*Z[2,1]*Z[2,2]^2+Z[2,2]^2*Z[2,0]*Z[0,1]*Z
[0,2]^2;}{%
\maplemultiline{
{T_{32}} := {Z_{0, \,0}}^{2}\,{Z_{0, \,1}}\,{Z_{1, \,2}}\,{Z_{1, 
\,0}}^{2} + {Z_{1, \,0}}^{2}\,{Z_{1, \,1}}\,{Z_{2, \,2}}\,{Z_{2, 
\,0}}^{2} + {Z_{2, \,0}}^{2}\,{Z_{2, \,1}}\,{Z_{0, \,2}}\,{Z_{0, 
\,0}}^{2} \\
\mbox{} + {Z_{0, \,1}}^{2}\,{Z_{0, \,2}}\,{Z_{1, \,0}}\,{Z_{1, \,
1}}^{2} + {Z_{1, \,1}}^{2}\,{Z_{1, \,2}}\,{Z_{2, \,0}}\,{Z_{2, \,
1}}^{2} + {Z_{2, \,1}}^{2}\,{Z_{2, \,2}}\,{Z_{0, \,0}}\,{Z_{0, \,
1}}^{2} \\
\mbox{} + {Z_{0, \,2}}^{2}\,{Z_{0, \,0}}\,{Z_{1, \,1}}\,{Z_{1, \,
2}}^{2} + {Z_{1, \,2}}^{2}\,{Z_{1, \,0}}\,{Z_{2, \,1}}\,{Z_{2, \,
2}}^{2} + {Z_{2, \,2}}^{2}\,{Z_{2, \,0}}\,{Z_{0, \,1}}\,{Z_{0, \,
2}}^{2} }
}
\end{maplelatex}
\end{maplegroup}
\emptyline
\begin{maplegroup}
\begin{mapleinput}
\mapleinline{active}{1d}{T[33]:=Gen(Z[0,0]^2*Z[0,2]*Z[1,1]*Z[1,0]^2);}{%
}
\end{mapleinput}

\mapleresult
\begin{maplelatex}
\mapleinline{inert}{2d}{T[33] :=
Z[0,0]^2*Z[0,2]*Z[1,1]*Z[1,0]^2+Z[1,0]^2*Z[1,2]*Z[2,1]*Z[2,0]^2+Z[2,0]
^2*Z[2,2]*Z[0,1]*Z[0,0]^2+Z[0,1]^2*Z[0,0]*Z[1,2]*Z[1,1]^2+Z[1,1]^2*Z[1
,0]*Z[2,2]*Z[2,1]^2+Z[2,1]^2*Z[2,0]*Z[0,2]*Z[0,1]^2+Z[0,2]^2*Z[0,1]*Z[
1,0]*Z[1,2]^2+Z[1,2]^2*Z[1,1]*Z[2,0]*Z[2,2]^2+Z[2,2]^2*Z[2,1]*Z[0,0]*Z
[0,2]^2;}{%
\maplemultiline{
{T_{33}} := {Z_{0, \,0}}^{2}\,{Z_{0, \,2}}\,{Z_{1, \,1}}\,{Z_{1, 
\,0}}^{2} + {Z_{1, \,0}}^{2}\,{Z_{1, \,2}}\,{Z_{2, \,1}}\,{Z_{2, 
\,0}}^{2} + {Z_{2, \,0}}^{2}\,{Z_{2, \,2}}\,{Z_{0, \,1}}\,{Z_{0, 
\,0}}^{2} \\
\mbox{} + {Z_{0, \,1}}^{2}\,{Z_{0, \,0}}\,{Z_{1, \,2}}\,{Z_{1, \,
1}}^{2} + {Z_{1, \,1}}^{2}\,{Z_{1, \,0}}\,{Z_{2, \,2}}\,{Z_{2, \,
1}}^{2} + {Z_{2, \,1}}^{2}\,{Z_{2, \,0}}\,{Z_{0, \,2}}\,{Z_{0, \,
1}}^{2} \\
\mbox{} + {Z_{0, \,2}}^{2}\,{Z_{0, \,1}}\,{Z_{1, \,0}}\,{Z_{1, \,
2}}^{2} + {Z_{1, \,2}}^{2}\,{Z_{1, \,1}}\,{Z_{2, \,0}}\,{Z_{2, \,
2}}^{2} + {Z_{2, \,2}}^{2}\,{Z_{2, \,1}}\,{Z_{0, \,0}}\,{Z_{0, \,
2}}^{2} }
}
\end{maplelatex}
\end{maplegroup}
\emptyline
\begin{maplegroup}
\begin{mapleinput}
\mapleinline{active}{1d}{T[34]:=Gen(Z[0,0]^2*Z[1,0]*Z[2,1]*Z[0,1]^2);}{%
}
\end{mapleinput}

\mapleresult
\begin{maplelatex}
\mapleinline{inert}{2d}{T[34] :=
Z[0,0]^2*Z[1,0]*Z[2,1]*Z[0,1]^2+Z[1,0]^2*Z[2,0]*Z[0,1]*Z[1,1]^2+Z[2,0]
^2*Z[0,0]*Z[1,1]*Z[2,1]^2+Z[0,1]^2*Z[1,1]*Z[2,2]*Z[0,2]^2+Z[1,1]^2*Z[2
,1]*Z[0,2]*Z[1,2]^2+Z[2,1]^2*Z[0,1]*Z[1,2]*Z[2,2]^2+Z[0,2]^2*Z[1,2]*Z[
2,0]*Z[0,0]^2+Z[1,2]^2*Z[2,2]*Z[0,0]*Z[1,0]^2+Z[2,2]^2*Z[0,2]*Z[1,0]*Z
[2,0]^2;}{%
\maplemultiline{
{T_{34}} := {Z_{0, \,0}}^{2}\,{Z_{1, \,0}}\,{Z_{2, \,1}}\,{Z_{0, 
\,1}}^{2} + {Z_{1, \,0}}^{2}\,{Z_{2, \,0}}\,{Z_{0, \,1}}\,{Z_{1, 
\,1}}^{2} + {Z_{2, \,0}}^{2}\,{Z_{0, \,0}}\,{Z_{1, \,1}}\,{Z_{2, 
\,1}}^{2} \\
\mbox{} + {Z_{0, \,1}}^{2}\,{Z_{1, \,1}}\,{Z_{2, \,2}}\,{Z_{0, \,
2}}^{2} + {Z_{1, \,1}}^{2}\,{Z_{2, \,1}}\,{Z_{0, \,2}}\,{Z_{1, \,
2}}^{2} + {Z_{2, \,1}}^{2}\,{Z_{0, \,1}}\,{Z_{1, \,2}}\,{Z_{2, \,
2}}^{2} \\
\mbox{} + {Z_{0, \,2}}^{2}\,{Z_{1, \,2}}\,{Z_{2, \,0}}\,{Z_{0, \,
0}}^{2} + {Z_{1, \,2}}^{2}\,{Z_{2, \,2}}\,{Z_{0, \,0}}\,{Z_{1, \,
0}}^{2} + {Z_{2, \,2}}^{2}\,{Z_{0, \,2}}\,{Z_{1, \,0}}\,{Z_{2, \,
0}}^{2} }
}
\end{maplelatex}
\end{maplegroup}
\emptyline
\begin{maplegroup}
\begin{mapleinput}
\mapleinline{active}{1d}{T[35]:=Gen(Z[0,0]^2*Z[1,0]*Z[2,2]*Z[0,2]^2);}{%
}
\end{mapleinput}

\mapleresult
\begin{maplelatex}
\mapleinline{inert}{2d}{T[35] :=
Z[0,0]^2*Z[1,0]*Z[2,2]*Z[0,2]^2+Z[1,0]^2*Z[2,0]*Z[0,2]*Z[1,2]^2+Z[2,0]
^2*Z[0,0]*Z[1,2]*Z[2,2]^2+Z[0,1]^2*Z[1,1]*Z[2,0]*Z[0,0]^2+Z[1,1]^2*Z[2
,1]*Z[0,0]*Z[1,0]^2+Z[2,1]^2*Z[0,1]*Z[1,0]*Z[2,0]^2+Z[0,2]^2*Z[1,2]*Z[
2,1]*Z[0,1]^2+Z[1,2]^2*Z[2,2]*Z[0,1]*Z[1,1]^2+Z[2,2]^2*Z[0,2]*Z[1,1]*Z
[2,1]^2;}{%
\maplemultiline{
{T_{35}} := {Z_{0, \,0}}^{2}\,{Z_{1, \,0}}\,{Z_{2, \,2}}\,{Z_{0, 
\,2}}^{2} + {Z_{1, \,0}}^{2}\,{Z_{2, \,0}}\,{Z_{0, \,2}}\,{Z_{1, 
\,2}}^{2} + {Z_{2, \,0}}^{2}\,{Z_{0, \,0}}\,{Z_{1, \,2}}\,{Z_{2, 
\,2}}^{2} \\
\mbox{} + {Z_{0, \,1}}^{2}\,{Z_{1, \,1}}\,{Z_{2, \,0}}\,{Z_{0, \,
0}}^{2} + {Z_{1, \,1}}^{2}\,{Z_{2, \,1}}\,{Z_{0, \,0}}\,{Z_{1, \,
0}}^{2} + {Z_{2, \,1}}^{2}\,{Z_{0, \,1}}\,{Z_{1, \,0}}\,{Z_{2, \,
0}}^{2} \\
\mbox{} + {Z_{0, \,2}}^{2}\,{Z_{1, \,2}}\,{Z_{2, \,1}}\,{Z_{0, \,
1}}^{2} + {Z_{1, \,2}}^{2}\,{Z_{2, \,2}}\,{Z_{0, \,1}}\,{Z_{1, \,
1}}^{2} + {Z_{2, \,2}}^{2}\,{Z_{0, \,2}}\,{Z_{1, \,1}}\,{Z_{2, \,
1}}^{2} }
}
\end{maplelatex}
\end{maplegroup}
\emptyline
\begin{maplegroup}
\begin{mapleinput}
\mapleinline{active}{1d}{T[36]:=Gen(Z[0,0]^2*Z[1,0]*Z[0,1]*Z[1,1]^2);}{%
}
\end{mapleinput}

\mapleresult
\begin{maplelatex}
\mapleinline{inert}{2d}{T[36] :=
Z[0,0]^2*Z[1,0]*Z[0,1]*Z[1,1]^2+Z[1,0]^2*Z[2,0]*Z[1,1]*Z[2,1]^2+Z[2,0]
^2*Z[0,0]*Z[2,1]*Z[0,1]^2+Z[0,1]^2*Z[1,1]*Z[0,2]*Z[1,2]^2+Z[1,1]^2*Z[2
,1]*Z[1,2]*Z[2,2]^2+Z[2,1]^2*Z[0,1]*Z[2,2]*Z[0,2]^2+Z[0,2]^2*Z[1,2]*Z[
0,0]*Z[1,0]^2+Z[1,2]^2*Z[2,2]*Z[1,0]*Z[2,0]^2+Z[2,2]^2*Z[0,2]*Z[2,0]*Z
[0,0]^2;}{%
\maplemultiline{
{T_{36}} := {Z_{0, \,0}}^{2}\,{Z_{1, \,0}}\,{Z_{0, \,1}}\,{Z_{1, 
\,1}}^{2} + {Z_{1, \,0}}^{2}\,{Z_{2, \,0}}\,{Z_{1, \,1}}\,{Z_{2, 
\,1}}^{2} + {Z_{2, \,0}}^{2}\,{Z_{0, \,0}}\,{Z_{2, \,1}}\,{Z_{0, 
\,1}}^{2} \\
\mbox{} + {Z_{0, \,1}}^{2}\,{Z_{1, \,1}}\,{Z_{0, \,2}}\,{Z_{1, \,
2}}^{2} + {Z_{1, \,1}}^{2}\,{Z_{2, \,1}}\,{Z_{1, \,2}}\,{Z_{2, \,
2}}^{2} + {Z_{2, \,1}}^{2}\,{Z_{0, \,1}}\,{Z_{2, \,2}}\,{Z_{0, \,
2}}^{2} \\
\mbox{} + {Z_{0, \,2}}^{2}\,{Z_{1, \,2}}\,{Z_{0, \,0}}\,{Z_{1, \,
0}}^{2} + {Z_{1, \,2}}^{2}\,{Z_{2, \,2}}\,{Z_{1, \,0}}\,{Z_{2, \,
0}}^{2} + {Z_{2, \,2}}^{2}\,{Z_{0, \,2}}\,{Z_{2, \,0}}\,{Z_{0, \,
0}}^{2} }
}
\end{maplelatex}
\end{maplegroup}
\emptyline
\begin{maplegroup}
\begin{mapleinput}
\mapleinline{active}{1d}{T[37]:=Gen(Z[0,0]^2*Z[0,1]*Z[2,0]*Z[2,1]^2);}{%
}
\end{mapleinput}

\mapleresult
\begin{maplelatex}
\mapleinline{inert}{2d}{T[37] :=
Z[0,0]^2*Z[0,1]*Z[2,0]*Z[2,1]^2+Z[1,0]^2*Z[1,1]*Z[0,0]*Z[0,1]^2+Z[2,0]
^2*Z[2,1]*Z[1,0]*Z[1,1]^2+Z[0,1]^2*Z[0,2]*Z[2,1]*Z[2,2]^2+Z[1,1]^2*Z[1
,2]*Z[0,1]*Z[0,2]^2+Z[2,1]^2*Z[2,2]*Z[1,1]*Z[1,2]^2+Z[0,2]^2*Z[0,0]*Z[
2,2]*Z[2,0]^2+Z[1,2]^2*Z[1,0]*Z[0,2]*Z[0,0]^2+Z[2,2]^2*Z[2,0]*Z[1,2]*Z
[1,0]^2;}{%
\maplemultiline{
{T_{37}} := {Z_{0, \,0}}^{2}\,{Z_{0, \,1}}\,{Z_{2, \,0}}\,{Z_{2, 
\,1}}^{2} + {Z_{1, \,0}}^{2}\,{Z_{1, \,1}}\,{Z_{0, \,0}}\,{Z_{0, 
\,1}}^{2} + {Z_{2, \,0}}^{2}\,{Z_{2, \,1}}\,{Z_{1, \,0}}\,{Z_{1, 
\,1}}^{2} \\
\mbox{} + {Z_{0, \,1}}^{2}\,{Z_{0, \,2}}\,{Z_{2, \,1}}\,{Z_{2, \,
2}}^{2} + {Z_{1, \,1}}^{2}\,{Z_{1, \,2}}\,{Z_{0, \,1}}\,{Z_{0, \,
2}}^{2} + {Z_{2, \,1}}^{2}\,{Z_{2, \,2}}\,{Z_{1, \,1}}\,{Z_{1, \,
2}}^{2} \\
\mbox{} + {Z_{0, \,2}}^{2}\,{Z_{0, \,0}}\,{Z_{2, \,2}}\,{Z_{2, \,
0}}^{2} + {Z_{1, \,2}}^{2}\,{Z_{1, \,0}}\,{Z_{0, \,2}}\,{Z_{0, \,
0}}^{2} + {Z_{2, \,2}}^{2}\,{Z_{2, \,0}}\,{Z_{1, \,2}}\,{Z_{1, \,
0}}^{2} }
}
\end{maplelatex}
\end{maplegroup}
\emptyline
\begin{maplegroup}
\begin{mapleinput}
\mapleinline{active}{1d}{T[38]:=Gen(Z[0,0]^2*Z[0,1]*Z[0,2]*Z[1,0]*Z[2,0]);}{%
}
\end{mapleinput}

\mapleresult
\begin{maplelatex}
\mapleinline{inert}{2d}{T[38] :=
Z[0,0]^2*Z[0,1]*Z[0,2]*Z[1,0]*Z[2,0]+Z[1,0]^2*Z[1,1]*Z[1,2]*Z[2,0]*Z[0
,0]+Z[2,0]^2*Z[2,1]*Z[2,2]*Z[0,0]*Z[1,0]+Z[0,1]^2*Z[0,2]*Z[0,0]*Z[1,1]
*Z[2,1]+Z[1,1]^2*Z[1,2]*Z[1,0]*Z[2,1]*Z[0,1]+Z[2,1]^2*Z[2,2]*Z[2,0]*Z[
0,1]*Z[1,1]+Z[0,2]^2*Z[0,0]*Z[0,1]*Z[1,2]*Z[2,2]+Z[1,2]^2*Z[1,0]*Z[1,1
]*Z[2,2]*Z[0,2]+Z[2,2]^2*Z[2,0]*Z[2,1]*Z[0,2]*Z[1,2];}{%
\maplemultiline{
{T_{38}} := {Z_{0, \,0}}^{2}\,{Z_{0, \,1}}\,{Z_{0, \,2}}\,{Z_{1, 
\,0}}\,{Z_{2, \,0}} + {Z_{1, \,0}}^{2}\,{Z_{1, \,1}}\,{Z_{1, \,2}
}\,{Z_{2, \,0}}\,{Z_{0, \,0}} + {Z_{2, \,0}}^{2}\,{Z_{2, \,1}}\,{
Z_{2, \,2}}\,{Z_{0, \,0}}\,{Z_{1, \,0}} \\
\mbox{} + {Z_{0, \,1}}^{2}\,{Z_{0, \,2}}\,{Z_{0, \,0}}\,{Z_{1, \,
1}}\,{Z_{2, \,1}} + {Z_{1, \,1}}^{2}\,{Z_{1, \,2}}\,{Z_{1, \,0}}
\,{Z_{2, \,1}}\,{Z_{0, \,1}} + {Z_{2, \,1}}^{2}\,{Z_{2, \,2}}\,{Z
_{2, \,0}}\,{Z_{0, \,1}}\,{Z_{1, \,1}} \\
\mbox{} + {Z_{0, \,2}}^{2}\,{Z_{0, \,0}}\,{Z_{0, \,1}}\,{Z_{1, \,
2}}\,{Z_{2, \,2}} + {Z_{1, \,2}}^{2}\,{Z_{1, \,0}}\,{Z_{1, \,1}}
\,{Z_{2, \,2}}\,{Z_{0, \,2}} + {Z_{2, \,2}}^{2}\,{Z_{2, \,0}}\,{Z
_{2, \,1}}\,{Z_{0, \,2}}\,{Z_{1, \,2}} }
}
\end{maplelatex}
\end{maplegroup}
\emptyline
\begin{maplegroup}
\begin{mapleinput}
\mapleinline{active}{1d}{T[39]:=Gen(Z[0,0]^2*Z[0,1]*Z[0,2]*Z[1,1]*Z[2,2]);}{%
}
\end{mapleinput}

\mapleresult
\begin{maplelatex}
\mapleinline{inert}{2d}{T[39] :=
Z[0,0]^2*Z[0,1]*Z[0,2]*Z[1,1]*Z[2,2]+Z[1,0]^2*Z[1,1]*Z[1,2]*Z[2,1]*Z[0
,2]+Z[2,0]^2*Z[2,1]*Z[2,2]*Z[0,1]*Z[1,2]+Z[0,1]^2*Z[0,2]*Z[0,0]*Z[1,2]
*Z[2,0]+Z[1,1]^2*Z[1,2]*Z[1,0]*Z[2,2]*Z[0,0]+Z[2,1]^2*Z[2,2]*Z[2,0]*Z[
0,2]*Z[1,0]+Z[0,2]^2*Z[0,0]*Z[0,1]*Z[1,0]*Z[2,1]+Z[1,2]^2*Z[1,0]*Z[1,1
]*Z[2,0]*Z[0,1]+Z[2,2]^2*Z[2,0]*Z[2,1]*Z[0,0]*Z[1,1];}{%
\maplemultiline{
{T_{39}} := {Z_{0, \,0}}^{2}\,{Z_{0, \,1}}\,{Z_{0, \,2}}\,{Z_{1, 
\,1}}\,{Z_{2, \,2}} + {Z_{1, \,0}}^{2}\,{Z_{1, \,1}}\,{Z_{1, \,2}
}\,{Z_{2, \,1}}\,{Z_{0, \,2}} + {Z_{2, \,0}}^{2}\,{Z_{2, \,1}}\,{
Z_{2, \,2}}\,{Z_{0, \,1}}\,{Z_{1, \,2}} \\
\mbox{} + {Z_{0, \,1}}^{2}\,{Z_{0, \,2}}\,{Z_{0, \,0}}\,{Z_{1, \,
2}}\,{Z_{2, \,0}} + {Z_{1, \,1}}^{2}\,{Z_{1, \,2}}\,{Z_{1, \,0}}
\,{Z_{2, \,2}}\,{Z_{0, \,0}} + {Z_{2, \,1}}^{2}\,{Z_{2, \,2}}\,{Z
_{2, \,0}}\,{Z_{0, \,2}}\,{Z_{1, \,0}} \\
\mbox{} + {Z_{0, \,2}}^{2}\,{Z_{0, \,0}}\,{Z_{0, \,1}}\,{Z_{1, \,
0}}\,{Z_{2, \,1}} + {Z_{1, \,2}}^{2}\,{Z_{1, \,0}}\,{Z_{1, \,1}}
\,{Z_{2, \,0}}\,{Z_{0, \,1}} + {Z_{2, \,2}}^{2}\,{Z_{2, \,0}}\,{Z
_{2, \,1}}\,{Z_{0, \,0}}\,{Z_{1, \,1}} }
}
\end{maplelatex}
\end{maplegroup}
\emptyline
\begin{maplegroup}
\begin{mapleinput}
\mapleinline{active}{1d}{T[40]:=Gen(Z[0,0]^2*Z[0,1]*Z[0,2]*Z[1,2]*Z[2,1]);}{%
}
\end{mapleinput}

\mapleresult
\begin{maplelatex}
\mapleinline{inert}{2d}{T[40] :=
Z[0,0]^2*Z[0,1]*Z[0,2]*Z[1,2]*Z[2,1]+Z[1,0]^2*Z[1,1]*Z[1,2]*Z[2,2]*Z[0
,1]+Z[2,0]^2*Z[2,1]*Z[2,2]*Z[0,2]*Z[1,1]+Z[0,1]^2*Z[0,2]*Z[0,0]*Z[1,0]
*Z[2,2]+Z[1,1]^2*Z[1,2]*Z[1,0]*Z[2,0]*Z[0,2]+Z[2,1]^2*Z[2,2]*Z[2,0]*Z[
0,0]*Z[1,2]+Z[0,2]^2*Z[0,0]*Z[0,1]*Z[1,1]*Z[2,0]+Z[1,2]^2*Z[1,0]*Z[1,1
]*Z[2,1]*Z[0,0]+Z[2,2]^2*Z[2,0]*Z[2,1]*Z[0,1]*Z[1,0];}{%
\maplemultiline{
{T_{40}} := {Z_{0, \,0}}^{2}\,{Z_{0, \,1}}\,{Z_{0, \,2}}\,{Z_{1, 
\,2}}\,{Z_{2, \,1}} + {Z_{1, \,0}}^{2}\,{Z_{1, \,1}}\,{Z_{1, \,2}
}\,{Z_{2, \,2}}\,{Z_{0, \,1}} + {Z_{2, \,0}}^{2}\,{Z_{2, \,1}}\,{
Z_{2, \,2}}\,{Z_{0, \,2}}\,{Z_{1, \,1}} \\
\mbox{} + {Z_{0, \,1}}^{2}\,{Z_{0, \,2}}\,{Z_{0, \,0}}\,{Z_{1, \,
0}}\,{Z_{2, \,2}} + {Z_{1, \,1}}^{2}\,{Z_{1, \,2}}\,{Z_{1, \,0}}
\,{Z_{2, \,0}}\,{Z_{0, \,2}} + {Z_{2, \,1}}^{2}\,{Z_{2, \,2}}\,{Z
_{2, \,0}}\,{Z_{0, \,0}}\,{Z_{1, \,2}} \\
\mbox{} + {Z_{0, \,2}}^{2}\,{Z_{0, \,0}}\,{Z_{0, \,1}}\,{Z_{1, \,
1}}\,{Z_{2, \,0}} + {Z_{1, \,2}}^{2}\,{Z_{1, \,0}}\,{Z_{1, \,1}}
\,{Z_{2, \,1}}\,{Z_{0, \,0}} + {Z_{2, \,2}}^{2}\,{Z_{2, \,0}}\,{Z
_{2, \,1}}\,{Z_{0, \,1}}\,{Z_{1, \,0}} }
}
\end{maplelatex}
\end{maplegroup}
\emptyline
\begin{maplegroup}
\begin{mapleinput}
\mapleinline{active}{1d}{T[41]:=Gen(Z[0,0]^2*Z[1,0]*Z[2,0]*Z[1,1]*Z[2,2]);}{%
}
\end{mapleinput}

\mapleresult
\begin{maplelatex}
\mapleinline{inert}{2d}{T[41] :=
Z[0,0]^2*Z[1,0]*Z[2,0]*Z[1,1]*Z[2,2]+Z[1,0]^2*Z[2,0]*Z[0,0]*Z[2,1]*Z[0
,2]+Z[2,0]^2*Z[0,0]*Z[1,0]*Z[0,1]*Z[1,2]+Z[0,1]^2*Z[1,1]*Z[2,1]*Z[1,2]
*Z[2,0]+Z[1,1]^2*Z[2,1]*Z[0,1]*Z[2,2]*Z[0,0]+Z[2,1]^2*Z[0,1]*Z[1,1]*Z[
0,2]*Z[1,0]+Z[0,2]^2*Z[1,2]*Z[2,2]*Z[1,0]*Z[2,1]+Z[1,2]^2*Z[2,2]*Z[0,2
]*Z[2,0]*Z[0,1]+Z[2,2]^2*Z[0,2]*Z[1,2]*Z[0,0]*Z[1,1];}{%
\maplemultiline{
{T_{41}} := {Z_{0, \,0}}^{2}\,{Z_{1, \,0}}\,{Z_{2, \,0}}\,{Z_{1, 
\,1}}\,{Z_{2, \,2}} + {Z_{1, \,0}}^{2}\,{Z_{2, \,0}}\,{Z_{0, \,0}
}\,{Z_{2, \,1}}\,{Z_{0, \,2}} + {Z_{2, \,0}}^{2}\,{Z_{0, \,0}}\,{
Z_{1, \,0}}\,{Z_{0, \,1}}\,{Z_{1, \,2}} \\
\mbox{} + {Z_{0, \,1}}^{2}\,{Z_{1, \,1}}\,{Z_{2, \,1}}\,{Z_{1, \,
2}}\,{Z_{2, \,0}} + {Z_{1, \,1}}^{2}\,{Z_{2, \,1}}\,{Z_{0, \,1}}
\,{Z_{2, \,2}}\,{Z_{0, \,0}} + {Z_{2, \,1}}^{2}\,{Z_{0, \,1}}\,{Z
_{1, \,1}}\,{Z_{0, \,2}}\,{Z_{1, \,0}} \\
\mbox{} + {Z_{0, \,2}}^{2}\,{Z_{1, \,2}}\,{Z_{2, \,2}}\,{Z_{1, \,
0}}\,{Z_{2, \,1}} + {Z_{1, \,2}}^{2}\,{Z_{2, \,2}}\,{Z_{0, \,2}}
\,{Z_{2, \,0}}\,{Z_{0, \,1}} + {Z_{2, \,2}}^{2}\,{Z_{0, \,2}}\,{Z
_{1, \,2}}\,{Z_{0, \,0}}\,{Z_{1, \,1}} }
}
\end{maplelatex}
\end{maplegroup}
\emptyline
\begin{maplegroup}
\begin{mapleinput}
\mapleinline{active}{1d}{T[42]:=Gen(Z[0,0]^2*Z[1,0]*Z[2,0]*Z[1,2]*Z[2,1]);}{%
}
\end{mapleinput}

\mapleresult
\begin{maplelatex}
\mapleinline{inert}{2d}{T[42] :=
Z[0,0]^2*Z[1,0]*Z[2,0]*Z[1,2]*Z[2,1]+Z[1,0]^2*Z[2,0]*Z[0,0]*Z[2,2]*Z[0
,1]+Z[2,0]^2*Z[0,0]*Z[1,0]*Z[0,2]*Z[1,1]+Z[0,1]^2*Z[1,1]*Z[2,1]*Z[1,0]
*Z[2,2]+Z[1,1]^2*Z[2,1]*Z[0,1]*Z[2,0]*Z[0,2]+Z[2,1]^2*Z[0,1]*Z[1,1]*Z[
0,0]*Z[1,2]+Z[0,2]^2*Z[1,2]*Z[2,2]*Z[1,1]*Z[2,0]+Z[1,2]^2*Z[2,2]*Z[0,2
]*Z[2,1]*Z[0,0]+Z[2,2]^2*Z[0,2]*Z[1,2]*Z[0,1]*Z[1,0];}{%
\maplemultiline{
{T_{42}} := {Z_{0, \,0}}^{2}\,{Z_{1, \,0}}\,{Z_{2, \,0}}\,{Z_{1, 
\,2}}\,{Z_{2, \,1}} + {Z_{1, \,0}}^{2}\,{Z_{2, \,0}}\,{Z_{0, \,0}
}\,{Z_{2, \,2}}\,{Z_{0, \,1}} + {Z_{2, \,0}}^{2}\,{Z_{0, \,0}}\,{
Z_{1, \,0}}\,{Z_{0, \,2}}\,{Z_{1, \,1}} \\
\mbox{} + {Z_{0, \,1}}^{2}\,{Z_{1, \,1}}\,{Z_{2, \,1}}\,{Z_{1, \,
0}}\,{Z_{2, \,2}} + {Z_{1, \,1}}^{2}\,{Z_{2, \,1}}\,{Z_{0, \,1}}
\,{Z_{2, \,0}}\,{Z_{0, \,2}} + {Z_{2, \,1}}^{2}\,{Z_{0, \,1}}\,{Z
_{1, \,1}}\,{Z_{0, \,0}}\,{Z_{1, \,2}} \\
\mbox{} + {Z_{0, \,2}}^{2}\,{Z_{1, \,2}}\,{Z_{2, \,2}}\,{Z_{1, \,
1}}\,{Z_{2, \,0}} + {Z_{1, \,2}}^{2}\,{Z_{2, \,2}}\,{Z_{0, \,2}}
\,{Z_{2, \,1}}\,{Z_{0, \,0}} + {Z_{2, \,2}}^{2}\,{Z_{0, \,2}}\,{Z
_{1, \,2}}\,{Z_{0, \,1}}\,{Z_{1, \,0}} }
}
\end{maplelatex}
\end{maplegroup}
\emptyline
\begin{maplegroup}
\begin{mapleinput}
\mapleinline{active}{1d}{T[43]:=Gen(Z[0,0]^2*Z[1,1]*Z[2,2]*Z[1,2]*Z[2,1]);}{%
}
\end{mapleinput}

\mapleresult
\begin{maplelatex}
\mapleinline{inert}{2d}{T[43] :=
Z[0,0]^2*Z[1,1]*Z[2,2]*Z[1,2]*Z[2,1]+Z[1,0]^2*Z[2,1]*Z[0,2]*Z[2,2]*Z[0
,1]+Z[2,0]^2*Z[0,1]*Z[1,2]*Z[0,2]*Z[1,1]+Z[0,1]^2*Z[1,2]*Z[2,0]*Z[1,0]
*Z[2,2]+Z[1,1]^2*Z[2,2]*Z[0,0]*Z[2,0]*Z[0,2]+Z[2,1]^2*Z[0,2]*Z[1,0]*Z[
0,0]*Z[1,2]+Z[0,2]^2*Z[1,0]*Z[2,1]*Z[1,1]*Z[2,0]+Z[1,2]^2*Z[2,0]*Z[0,1
]*Z[2,1]*Z[0,0]+Z[2,2]^2*Z[0,0]*Z[1,1]*Z[0,1]*Z[1,0];}{%
\maplemultiline{
{T_{43}} := {Z_{0, \,0}}^{2}\,{Z_{1, \,1}}\,{Z_{2, \,2}}\,{Z_{1, 
\,2}}\,{Z_{2, \,1}} + {Z_{1, \,0}}^{2}\,{Z_{2, \,1}}\,{Z_{0, \,2}
}\,{Z_{2, \,2}}\,{Z_{0, \,1}} + {Z_{2, \,0}}^{2}\,{Z_{0, \,1}}\,{
Z_{1, \,2}}\,{Z_{0, \,2}}\,{Z_{1, \,1}} \\
\mbox{} + {Z_{0, \,1}}^{2}\,{Z_{1, \,2}}\,{Z_{2, \,0}}\,{Z_{1, \,
0}}\,{Z_{2, \,2}} + {Z_{1, \,1}}^{2}\,{Z_{2, \,2}}\,{Z_{0, \,0}}
\,{Z_{2, \,0}}\,{Z_{0, \,2}} + {Z_{2, \,1}}^{2}\,{Z_{0, \,2}}\,{Z
_{1, \,0}}\,{Z_{0, \,0}}\,{Z_{1, \,2}} \\
\mbox{} + {Z_{0, \,2}}^{2}\,{Z_{1, \,0}}\,{Z_{2, \,1}}\,{Z_{1, \,
1}}\,{Z_{2, \,0}} + {Z_{1, \,2}}^{2}\,{Z_{2, \,0}}\,{Z_{0, \,1}}
\,{Z_{2, \,1}}\,{Z_{0, \,0}} + {Z_{2, \,2}}^{2}\,{Z_{0, \,0}}\,{Z
_{1, \,1}}\,{Z_{0, \,1}}\,{Z_{1, \,0}} }
}
\end{maplelatex}
\end{maplegroup}
\emptyline
\begin{maplegroup}
Proc\'edure pour restreindre aux  4 plans correspondants aux points de
la forme (x, 00):
\end{maplegroup}
\begin{maplegroup}
\begin{mapleinput}
\mapleinline{active}{1d}{                                                                                                                                                                                                                                                                    
H:=proc(P,r,s)   x:=P; for j from 0 to 2  do for i from 0 to 2  do if 
q(r*i,s*j) <> 0 then      x:=subs(Z[i,j]=0,x) else end if od od;
RETURN(x) end;     }{%
}
\end{mapleinput}

\mapleresult

\begin{maplelatex}
\mapleinline{inert}{2d}{H := proc (P, r, s) local x, j, i; x := P; for j from 0 to 2 do for i
from 0 to 2 do if q(r*i,s*j) <> 0 then x := subs(Z[i,j] = 0,x) else
NULL end if end do end do; RETURN(x) end proc;}{%
\maplemultiline{
H := \textbf{proc} (P, \,r, \,s) \\
\textbf{local} \,x, \,j, \,i; \\
\mapleIndent{1} x := P\,; \\
\mapleIndent{1} \textbf{for} \,j\,\textbf{from} \,0\,\textbf{to} 
\,2\,\textbf{do} \textbf{for} \,i\,\textbf{from} \,0\,\textbf{to
} \,2\,\textbf{do}  \\
\mapleIndent{3} \textbf{if} \,\mathrm{q}(r \ast i, \,s \ast j)
\neq 0\,\textbf{then} \,x := \mathrm{subs}({Z_{i, \,j}}=0, \,x)\,
\textbf{else} \,\,\textbf{end if}  \\
\mapleIndent{2} \textbf{end do}  \\
\mapleIndent{1} \textbf{end do} ; \\
\mapleIndent{1} \mathrm{RETURN}(x) \\
\textbf{end proc}  }
}
\end{maplelatex}
\end{maplegroup}
\emptyline
\begin{maplegroup}
\begin{mapleinput}
\mapleinline{active}{1d}{M:= array(1..39);}{%
}
\end{mapleinput}
\mapleresult
\begin{maplelatex}
\mapleinline{inert}{2d}{M := array(1 .. 39,[]);}{%
\[
M := \mathrm{array}(1 .. 39, \,[])
\]
}
\end{maplelatex}
\end{maplegroup}
\emptyline
\begin{maplegroup}
\begin{mapleinput}
\mapleinline{active}{1d}{x:=1; for i from 1 to 43 do if H(T[i],0,1)=0 then M[x]:=T[i]; 
x:=x+1 else end if od: print(x); }{%
}
\end{mapleinput}
\end{maplegroup}
\begin{maplegroup}
\mapleresult
\begin{maplelatex}
\mapleinline{inert}{2d}{40;}{%
\[
40
\]
}
\end{maplelatex}
\end{maplegroup}
\begin{maplegroup}
\begin{mapleinput}
\mapleinline{active}{1d}{M2:= array(1..39);
}{%
}
\end{mapleinput}

\mapleresult
\begin{maplelatex}
\mapleinline{inert}{2d}{M2 := array(1 .. 39,[]);}{%
\[
\mathit{M2} := \mathrm{array}(1 .. 39, \,[])
\]
}
\end{maplelatex}
\end{maplegroup}
\emptyline
\begin{maplegroup}
\begin{mapleinput}
\mapleinline{active}{1d}{x2:=1; for i from 1 to 39 do if H(M[i],1,0)=0 then M2[x2]:=M[i];
x2:=x2+1 else end if od: }{%
}
\end{mapleinput}
\end{maplegroup}
\begin{maplegroup}
\begin{mapleinput}
\mapleinline{active}{1d}{print(x2);}{%
}
\end{mapleinput}

\mapleresult
\begin{maplelatex}
\mapleinline{inert}{2d}{37;}{%
\[
37
\]
}
\end{maplelatex}

\end{maplegroup}
\begin{maplegroup}
\begin{mapleinput}
\mapleinline{active}{1d}{M3:= array(1..39);
}{%
}
\end{mapleinput}

\mapleresult
\begin{maplelatex}
\mapleinline{inert}{2d}{M3 := array(1 .. 39,[]);}{%
\[
\mathit{M3} := \mathrm{array}(1 .. 39, \,[])
\]
}
\end{maplelatex}

\end{maplegroup}
\emptyline
\begin{maplegroup}
\begin{mapleinput}
\mapleinline{active}{1d}{x3:=1; for i from 1 to 39 do if H(M2[i],1,1)=0 then M3[x3]:=M2[i];
x3:=x3+1 else end if od: }{%
}
\end{mapleinput}
\end{maplegroup}

\begin{maplegroup}
\begin{mapleinput}
\mapleinline{active}{1d}{print(x3);}{%
}
\end{mapleinput}

\mapleresult
\begin{maplelatex}
\mapleinline{inert}{2d}{34;}{%
\[
34
\]
}
\end{maplelatex}
\end{maplegroup}
\emptyline
\begin{maplegroup}
\begin{mapleinput}
\mapleinline{active}{1d}{M4:= array(1..39);
}{%
}
\end{mapleinput}

\mapleresult
\begin{maplelatex}
\mapleinline{inert}{2d}{M4 := array(1 .. 39,[]);}{%
\[
\mathit{M4} := \mathrm{array}(1 .. 39, \,[])
\]
}
\end{maplelatex}
\end{maplegroup}

\begin{maplegroup}
\begin{mapleinput}
\mapleinline{active}{1d}{x4:=1; for i from 1 to 39 do if H(M3[i],1,2)=0 then M4[x4]:=M3[i];
x4:=x4+1 else end if od: }{%
}
\end{mapleinput}
\end{maplegroup}

\begin{maplegroup}
\begin{mapleinput}
\mapleinline{active}{1d}{print(x4);}{%
}
\end{mapleinput}

\mapleresult
\begin{maplelatex}
\mapleinline{inert}{2d}{31;}{%
\[
31
\]
}
\end{maplelatex}
\end{maplegroup}
\begin{maplegroup}
\emptyline
On d\'eduit que  le rang de l'image de l'application
$\nu$  est au moins 13.  \\ On r\'eduit \`a 30  la
dimension de  l'espace de d\'epart.\\
\\
Variables dans un plan de points fixes:
%\emptyline
\end{maplegroup}
\begin{maplegroup}
\begin{mapleinput}
\mapleinline{active}{1d}{Y:=array(0..2);}{%
}
\end{mapleinput}

\mapleresult
\begin{maplelatex}
\mapleinline{inert}{2d}{Y := array(0 .. 2,[]);}{%
\[
Y := \mathrm{array}(0 .. 2, \,[])
\]
}
\end{maplelatex}

\end{maplegroup}
\begin{maplegroup}
\begin{mapleinput}
\mapleinline{active}{1d}{with(linalg);}{%
}
\end{mapleinput}
\end{maplegroup}
\emptyline
\begin{maplegroup}

Matrice de coefficients:
%\emptyline
\end{maplegroup}
\begin{maplegroup}
\begin{mapleinput}
\mapleinline{active}{1d}{C:=matrix(144,30);}{%
}
\end{mapleinput}

\mapleresult
\begin{maplelatex}
\mapleinline{inert}{2d}{C := array(1 .. 144,1 .. 30,[]);}{%
\[
C := \mathrm{array}(1 .. 144, \,1 .. 30, \,[])
\]
}
\end{maplelatex}

\end{maplegroup}
\begin{maplegroup}
%\begin{maplelatex}
\emptyline
Il  y a  36 espaces de sextiques $K_{\eta}$-invariantes  $(S^6V_{\eta})^{K_{\eta}}$, chacun \\
de dimension 4.  
\\
Matrice de polyn\^omes restreints aux espaces
$V_{\eta}$ :
\end{maplegroup}
\begin{maplegroup}
\begin{mapleinput}
\mapleinline{active}{1d}{N:=matrix(36,30);}{%
}
\end{mapleinput}

\mapleresult
\begin{maplelatex}
\mapleinline{inert}{2d}{N := array(1 .. 36,1 .. 30,[]);}{%
\[
N := \mathrm{array}(1 .. 36, \,1 .. 30, \,[])
\]
}
\end{maplelatex}

\end{maplegroup}
\begin{maplegroup}
\begin{mapleinput}
\mapleinline{active}{1d}{for i from 1 to 36 do for j from 1 to 30 do N[i,j]:= M4[j] od od ; 

                  }{%
}
\end{mapleinput}
\end{maplegroup}
\emptyline
\begin{maplegroup}

Produit scalaire:

\end{maplegroup}
\begin{maplegroup}
\begin{mapleinput}
\mapleinline{active}{1d}{prod:= proc(u,v,x1,x2)  RETURN( (u*x1 + v*x2) mod 3) end;}{%
}
\end{mapleinput}

\mapleresult
\begin{maplelatex}
\mapleinline{inert}{2d}{prod := proc (u, v, x1, x2) RETURN(`mod`(u*x1+v*x2,3)) end proc;}{%
\[
\mathit{prod} := \textbf{proc} (u, \,v, \,\mathit{x1}, \,\mathit{
x2})\,\mathrm{RETURN}((u \ast \mathit{x1} + v \ast \mathit{x2})\,
\mathrm{mod}\,3)\,\textbf{end proc} 
\]
}
\end{maplelatex}
\end{maplegroup}
\emptyline
\begin{maplegroup}
Restriction aux sous-espaces de points fixes correspondants aux points
la forme (01, x*):
\emptyline
\end{maplegroup}
\begin{maplegroup}
\begin{mapleinput}
\mapleinline{active}{1d}{for u from 0 to 2 do 
for v from 0 to 2 do b:=3*u+v+1;                
for j from 1 to 30 do                               
N[b,j]:= subs(\{Z[0,0]=Y[0], Z[0,2]=w^prod(u,v,0,1)*Y[0],                      
Z[0,1]=Y[0], Z[1,0]=Y[1],
Z[1,2]=w^prod(u,v,2,1)*Y[1], Z[1,1]=w^prod(u,v,1,0)*Y[1], 
Z[2,0]=Y[2], Z[2,2]=w^prod(u,v,1,1)*Y[2], 
Z[2,1]=w^prod(u,v,2,0)*Y[2]\},N[b,j]);        
N[b,j]:=subs(\{w^4=w,w^5=w^2,w^6=1,w^3=1,w^7=w,w^8=w^2,
w^9=1,w^10=w \},N[b,j]);    simplify (N[b,j])                                      
od; od od;                                
 }{%
}
\end{mapleinput}
\end{maplegroup}

\emptyline
\mapleresult
%\begin{maplegroupe}
Remplissage de la matrice de coefficients:
\emptyline
%\end{maplegroupe}
\begin{maplegroup}
\begin{mapleinput}
\mapleinline{active}{1d}{for u from 0 to 2 do                                                 
                    for v from 0 to 2 do b:=3*u+v+1;                  
for j from 1 to 30 do                                                                    
if coeff(subs(\{Y[1]=1,Y[2]=1\},N[b,j]), Y[0]^2)=0 
	then C[4*(b-1)+1,j]:=0                        
else C[4*(b-1)+1,j]:=lcoeff(N[b,j],[Y[0],Y[1],Y[2]],'t') 
	end if;           
if coeff(subs(\{Y[1]=1,Y[2]=1\},N[b,j]), Y[0]^3)=0 
	then C[4*(b-1)+2,j]:=0
else  C[4*(b-1)+2,j]:=lcoeff(N[b,j],[Y[0],Y[1],Y[2]],'t')
	end if;                                                               
if coeff(subs(\{Y[1]=1,Y[2]=1\},N[b,j]), Y[0]^4)=0 
	then C[4*(b-1)+3,j]:=0       
else C[4*(b-1)+3,j]:=lcoeff(N[b,j], [Y[0],Y[1],Y[2]],'t') 
	end if;          
if coeff(subs(\{Y[1]=1,Y[2]=1\},N[b,j]), Y[0]^6)=0
	then C[4*(b-1)+4,j]:=0                    
else C[4*(b-1)+4,j]:=lcoeff(N[b,j],[Y[0],Y[1],Y[2]],'t') 
	end if;           
       od:  od od; 
}{%
}
\end{mapleinput}

\end{maplegroup}
\begin{maplegroup}
%\begin{maplelatex}
\emptyline
Restriction aux sous-espaces avec  $\eta=  (10,x^*)$:
\emptyline
%\end{maplelatex}
\end{maplegroup}
\begin{maplegroup}
\begin{mapleinput}
\mapleinline{active}{1d}{for u from 0 to 2 do 
for v from 0 to 2 do b:=3*u+v+1;
for j from 1 to 30 do                                               
N[b+9,j]:= subs(\{Z[0,0]=Y[0], Z[2,0]=w^prod(u,v,1,0)*Y[0], 
Z[1,0]=Y[0], Z[0,1]=Y[1], Z[2,1]=w^prod(u,v,1,2)*Y[1], 
Z[1,1]=w^prod(u,v,0,1)*Y[1], Z[0,2]=Y[2], 
Z[2,2]=w^prod(u,v,1,1)*Y[2], 
Z[1,2]=w^prod(u,v,0,2)*Y[2]\},N[b+9,j]); 
N[b+9,j]:=subs(\{w^4=w,w^5=w^2,w^6=1,w^3=1,w^7=w,w^8=w^2,
w^9=1,w^10=w \},  N[b+9,j]);    
simplify (N[b+9,j])                                
      od; od od;      }{%
}
\end{mapleinput}
\end{maplegroup}
\emptyline
\begin{maplegroup}
\begin{mapleinput}
\mapleinline{active}{1d}{for u from 0 to 2 do 
for v from 0 to 2 do b:=3*u+v+1;                
for j from 1 to 30 do                                                                    
 if coeff(subs(\{Y[1]=1,Y[2]=1\},N[b+9,j]), Y[0]^2)=0 
	then C[4*(b+8)+1,j]:=0             
else C[4*(b+8)+1,j]:=lcoeff(N[b+9,j],[Y[0],Y[1],Y[2]],'t') 
	end if;         
if coeff(subs(\{Y[1]=1,Y[2]=1\},N[b+9,j]), Y[0]^3)=0 
	then C[4*(b+8)+2,j]:=0                     
else C[4*(b+8)+2,j]:=lcoeff(N[b+9,j],[Y[0],Y[1],Y[2]],'t') 
	end if;        
if coeff(subs(\{Y[1]=1,Y[2]=1\},N[b+9,j]), Y[0]^4)=0
	then C[4*(b+8)+3,j]:=0                     
else C[4*(b+8)+3,j]:=lcoeff(N[b+9,j], [Y[0],Y[1],Y[2]],'t') 
	end if;        
if coeff(subs(\{Y[1]=1,Y[2]=1\},N[b+9,j]),Y[0]^6)=0 
	then C[4*(b+8)+4,j]:=0                
else C[4*(b+8)+4,j]:=lcoeff(N[b+9,j],[Y[0],Y[1],Y[2]],'t') 
end if;         
            od:  od od; 
}{%
}
\end{mapleinput}
\end{maplegroup}

\begin{maplegroup}
\emptyline
Restriction aux sous-espaces avec $\eta = (11,x^*):$
\end{maplegroup}
\emptyline
\begin{maplegroup}
\begin{mapleinput}
\mapleinline{active}{1d}{for u from 0 to 2 do 
for v from 0 to 2 do b:=3*u+v+1;                
  for j from 1 to 30 do                                     
N[b+18,j]:= subs(\{Z[0,0]=Y[0], Z[2,2]=w^prod(u,v,1,1)*Y[0], 
Z[1,1]=Y[0], Z[0,1]=Y[1],
Z[2,0]=w^prod(u,v,1,0)*Y[1], Z[1,2]=w^prod(u,v,0,1)*Y[1], 
Z[0,2]=Y[2], Z[2,1]=w^prod(u,v,1,2)*Y[2], 
Z[1,0]=w^prod(u,v,0,2)*Y[2]\},N[b+18,j]);
   N[b+18,j]:=subs(\{w^4=w,w^5=w^2,w^6=1,w^3=1,w^7=w,w^8=w^2,
w^9=1,w^10=w \},N[b+18,j]);    
simplify (N[b+18,j])                                
                  od; od od;      }{%
}
\end{mapleinput}
\end{maplegroup}
\emptyline
\begin{maplegroup}
\begin{mapleinput}
\mapleinline{active}{1d}{for u from 0 to 2 do                                                 
   for v from 0 to 2 do b:=3*u+v+1;                  
for j from 1 to 30 do                                                                   
if coeff(subs(\{Y[1]=1,Y[2]=1\},N[b+18,j]), Y[0]^2)=0 
then C[4*(b+17)+1,j]:=0                       
else C[4*(b+17)+1,j]:=lcoeff(N[b+18,j],[Y[0],Y[1],Y[2]],'t') 
end if;       
if coeff(subs(\{Y[1]=1,Y[2]=1\},N[b+18,j]), Y[0]^3)=0 
then C[4*(b+17)+2,j]:=0                   
else C[4*(b+17)+2,j]:=lcoeff(N[b+18,j],[Y[0],Y[1],Y[2]],'t') 
end if;      
if coeff(subs(\{Y[1]=1,Y[2]=1\},N[b+18,j]), Y[0]^4)=0
then C[4*(b+17)+3,j]:=0               
else C[4*(b+17)+3,j]:=lcoeff(N[b+18,j], [Y[0],Y[1],Y[2]],'t') 
end if;      
if coeff(subs(\{Y[1]=1,Y[2]=1\},N[b+18,j]), Y[0]^6)=0 
then C[4*(b+17)+4,j]:=0              
else C[4*(b+17)+4,j]:=lcoeff(N[b+18,j],[Y[0],Y[1],Y[2]],'t') 
end if;       
    od:  od od;              }{%
}
\end{mapleinput}
\end{maplegroup}   
\begin{maplegroup}
\emptyline
Restriction aux sous-espaces  avec $\eta = (12,x^*):$
\emptyline
\end{maplegroup}
\begin{maplegroup}
\begin{mapleinput}
\mapleinline{active}{1d}{for u from 0 to 2 do                                                 
for v from 0 to 2 do b:=3*u+v+1;                
for j from 1 to 30 do                                                           
N[b+27,j]:= subs(\{Z[0,0]=Y[0], Z[2,1]=w^prod(u,v,1,2)*Y[0],
Z[1,2]=Y[0], Z[0,1]=Y[1], Z[2,2]=w^prod(u,v,1,1)*Y[1],
Z[1,0]=w^prod(u,v,0,1)*Y[1], Z[0,2]=Y[2], 
Z[2,0]=w^prod(u,v,1,0)*Y[2],
Z[1,1]=w^prod(u,v,0,2)*Y[2]\},N[b+27,j]);    
N[b+27,j]:=subs(\{w^4=w,w^5=w^2,w^6=1,w^3=1,w^7=w,w^8=w^2,
w^9=1,w^10=w \},N[b+27,j]);    
simplify (N[b+27,j])                                
                     od; od od;      }{%
}
\end{mapleinput}
\end{maplegroup}
\emptyline
\begin{maplegroup}
\begin{mapleinput}
\mapleinline{active}{1d}{for u from 0 to 2 do                                                 
for v from 0 to 2 do b:=3*u+v+1;                       
for j from 1 to 30 do    
if coeff(subs(\{Y[1]=1,Y[2]=1\},N[b+27,j]), Y[0]^2)=0 
then C[4*(b+26)+1,j]:=0                   
else C[4*(b+26)+1,j]:=lcoeff(N[b+27,j],[Y[0],Y[1],Y[2]],'t')
 end if;       
if coeff(subs(\{Y[1]=1,Y[2]=1\},N[b+27,j]), Y[0]^3)=0 
then C[4*(b+26)+2,j]:=0             
else C[4*(b+26)+2,j]:=lcoeff(N[b+27,j],[Y[0],Y[1],Y[2]],'t')
end if;      
if coeff(subs(\{Y[1]=1,Y[2]=1\},N[b+27,j]), Y[0]^4)=0
then C[4*(b+26)+3,j]:=0            
else C[4*(b+26)+3,j]:=lcoeff(N[b+27,j], [Y[0],Y[1],Y[2]],'t') 
end if;      
if coeff(subs(\{Y[1]=1,Y[2]=1\},N[b+27,j]),
Y[0]^6)=0 then C[4*(b+26)+4,j]:=0            
else C[4*(b+26)+4,j]:=lcoeff(N[b+27,j],[Y[0],Y[1],Y[2]],'t') 
end if;       
            od:  od od;   }{%
}
\end{mapleinput}
\end{maplegroup}
\emptyline
\begin{maplegroup}
\begin{mapleinput}
\mapleinline{active}{1d}{e:=array(1..144);}{%
}
\end{mapleinput}

\mapleresult
\begin{maplelatex}
\mapleinline{inert}{2d}{e := array(1 .. 144,[]);}{%
\[
e := \mathrm{array}(1 .. 144, \,[])
\]
}
\end{maplelatex}
\end{maplegroup}
\begin{maplegroup}
\begin{mapleinput}
\mapleinline{active}{1d}{for i from 1 to 144 do e[i]:=0 od:
}{%
}
\end{mapleinput}

\end{maplegroup}
\begin{maplegroup}
\emptyline 

Calcul du rang de la matrice de coefficients et d'une base pour le
noyau de $\nu:$
\emptyline
\end{maplegroup}
\begin{maplegroup}
\begin{mapleinput}
\mapleinline{active}{1d}{linsolve(C,e,'r');
}{%
}
\end{mapleinput}
\mapleresult
\begin{maplelatex}
\mapleinline{inert}{2d}{vector([0, 0, -t[1], t[1], -t[2], t[2], -t[3], t[3], 
0, 0, 0, 0,
\\
0, 0, 0, 0, 0, 0, 0, 0, 0, 0, 0, 0, 0, 0, 0, 0, 0, 0]);}
{%
\maplemultiline{
[ \  0, \,0, \, - {\mathit{t}_{1}}, \,{\mathit{t}_{1}}
, \, - {\mathit{t}_{2}}, \,{\mathit{t}_{2}}, \, - {\mathit{t}_{3}}, \,{\mathit{t}_{3}}, \,0, \,0, \,0, \,0, \,0, \\ 
\,0, \,0, \,0, \,0, \,0, \,0, \,0, \,0, \,0, \,0, \,0, \,0, \,0, \,0, \,0, \,0, \,0 \ ] }
}
\end{maplelatex}
\end{maplegroup}
\emptyline
\begin{maplegroup}
\begin{mapleinput}
\mapleinline{active}{1d}{r;}{%
}
\end{mapleinput}
\mapleresult
\begin{maplelatex}
\mapleinline{inert}{2d}{27;}{%
\[
27
\]
}
\end{maplelatex}
\end{maplegroup}
\begin{maplegroup}
\emptyline 
Et donc une base pour le noyau de $\nu$ est donn\'ee par : 
\begin{eqnarray*}
\omega_1:= T[11]-T[10] &&\\
\omega_2:= T[14]-T[13] &&\\
\omega_3:= T[17]-T[16] 
\end{eqnarray*}
\emptyline
\end{maplegroup}

%\end{document}
%% End of Maple 8.00 Output

\backmatter
\end{document}